\documentclass[a4paper,10pt,reqno]{amsart}
\usepackage{amsmath}
\usepackage{amssymb}
\usepackage{amsthm}
\usepackage[all]{xy}
\usepackage{booktabs}
\usepackage{multirow}
\usepackage{multicol}
\usepackage{aliascnt}
\usepackage{hyperref}

\newlength{\realsidemargin}
\newlength{\sidemargin}
\newcommand{\NewTheorem}[1]{
	\newaliascnt{#1}{TheoremEnvironment}
	\newtheorem{#1}[#1]{#1}
	\aliascntresetthe{#1}
	\expandafter\newcommand\csname #1autorefname\endcsname{#1}
}

\theoremstyle{definition}

\NewTheorem{Definition}
\NewTheorem{Remark}
\NewTheorem{Axiom}
\NewTheorem{Example}
\NewTheorem{Conventions}
\NewTheorem{Notations}
\NewTheorem{Setting}
\NewTheorem{Question}
\NewTheorem{Conjecture}

\theoremstyle{plain}
\NewTheorem{Proposition}
\NewTheorem{Lemma}
\NewTheorem{Theorem}
\NewTheorem{Corollary}

\newcommand{\Sectionref}[1]{Section~\ref{#1}}
\def\itemautorefname~#1\null{(#1)\null}
\def\equationautorefname~#1\null{(#1)\null}

\newcommand{\SwapSymbols}[1]{
	\expandafter\let\expandafter\temporarysymbol\csname #1\endcsname
	\expandafter\let\csname #1\expandafter\endcsname\csname var#1\endcsname
	\expandafter\let\csname var#1\endcsname\temporarysymbol
}

\SwapSymbols{Gamma}
\SwapSymbols{Delta}
\SwapSymbols{Theta}
\SwapSymbols{Lambda}
\SwapSymbols{Xi}
\SwapSymbols{Pi}
\SwapSymbols{Sigma}
\SwapSymbols{Upsilon}
\SwapSymbols{Phi}
\SwapSymbols{Psi}
\SwapSymbols{Omega}

\newcommand{\mbP}{\mathbb{P}}

\newcommand{\mbZ}{\mathbb{Z}}

\newcommand{\mcA}{\mathcal{A}}
\newcommand{\mcB}{\mathcal{B}}
\newcommand{\mcC}{\mathcal{C}}

\newcommand{\mcF}{\mathcal{F}}

\newcommand{\mcN}{\mathcal{N}}
\newcommand{\mcO}{\mathcal{O}}

\newcommand{\mcS}{\mathcal{S}}

\newcommand{\mcU}{\mathcal{U}}

\newcommand{\mcX}{\mathcal{X}}
\newcommand{\mcY}{\mathcal{Y}}
\newcommand{\mcZ}{\mathcal{Z}}

\newcommand{\mfa}{\mathfrak{a}}

\newcommand{\mfm}{\mathfrak{m}}

\newcommand{\mfp}{\mathfrak{p}}
\newcommand{\mfq}{\mathfrak{q}}

\let\originalleft\left
\let\originalright\right
\renewcommand{\left}{\mathopen{}\mathclose\bgroup\originalleft}
\renewcommand{\right}{\aftergroup\egroup\originalright}

\newcommand{\set}[2][]{\mathopen{#1\{}#2\mathclose{#1\}}}
\newcommand{\setwithspace}[2][]{\mathopen{#1\{}\,#2\,\mathclose{#1\}}}
\newcommand{\conditionalset}[3][]{\mathopen{#1\{}\,#2\mathrel{#1|}#3\,\mathclose{#1\}}}
\newcommand{\generatedset}[3][]{\mathopen{#1\langle}#3\mathclose{#1\rangle}_{\mathrm{#2}}}
\newcommand{\conditionalgeneratedset}[4][]{\mathopen{#1\langle}\,#3\mathrel{#1|}#4\,\mathclose{#1\rangle}_{\mathrm{#2}}}

\newcommand{\preloc}[2][]{{\generatedset[#1]{preloc}{#2}}}
\newcommand{\loc}[2][]{{\generatedset[#1]{loc}{#2}}}

\newcommand{\locfilt}[2][]{{\generatedset[#1]{locfilt}{#2}}}

\DeclareMathOperator{\Ob}{Ob}
\DeclareMathOperator{\Mor}{Mor}

\DeclareMathOperator{\Hom}{Hom}
\DeclareMathOperator{\End}{End}

\newcommand{\into}{\hookrightarrow}

\newcommand{\onto}{\twoheadrightarrow}

\newcommand{\isoto}{\xrightarrow{\smash{\raisebox{-0.25em}{$\sim$}}}}

\newcommand{\definehookarrow}{\newdir^{ (}{!/-5pt/@^{(}}}

\DeclareMathOperator{\Mod}{Mod}

\DeclareMathOperator{\QCoh}{QCoh}

\DeclareMathOperator{\Ann}{Ann}

\DeclareMathOperator{\Ker}{Ker}
\DeclareMathOperator{\Cok}{Cok}
\let\Im\relax
\DeclareMathOperator{\Im}{Im}

\DeclareMathOperator{\Spec}{Spec}

\DeclareMathOperator{\Ass}{Ass}
\DeclareMathOperator{\Supp}{Supp}

\DeclareMathOperator{\ASpec}{ASpec}

\DeclareMathOperator{\AAss}{AAss}
\DeclareMathOperator{\ASupp}{ASupp}

\newcommand{\usp}[1]{{|{#1}|}}
\newcommand{\uspX}{{\usp{X}}}
\newcommand{\OX}{\mcO_{X}}
\newcommand{\OXx}{\mcO_{X,x}}
\newcommand{\OXy}{\mcO_{X,y}}

\makeatletter

\renewcommand{\p@enumii}{}
\renewcommand{\p@enumiii}{}

\makeatother

\title{Classification of categorical subspaces of locally noetherian schemes}

\author{Ryo Kanda}
\thanks{The author is a Research Fellow of Japan Society for the Promotion of Science. This work is supported by Grant-in-Aid for JSPS Fellows 25$\cdot$249.}
\address{Graduate School of Mathematics, Nagoya University, Furo-cho, Chikusa-ku, Nagoya-shi, Aichi-ken, 464-8602, Japan}
\email{ryo.kanda.math@gmail.com}

\begin{document}

\begin{abstract}
	We classify the prelocalizing subcategories of the category of quasi-coherent sheaves on a locally noetherian scheme. In order to give the classification, we introduce the notion of a local filter of subobjects of the structure sheaf. The essential part of the argument is given as results on a Grothendieck category with certain properties. We also classify the localizing subcategories, the closed subcategories, and the bilocalizing subcategories in terms of filters.
\end{abstract}

\maketitle
\tableofcontents

\section{Introduction}
\label{sec:Intro}

Gabriel \cite{Gabriel} introduced a classification theory of subcategories of the category of modules over a ring. The theory relates several classes of subcategories to collections of ideals, and it reveals that geometry of the prime spectrum of a commutative ring is reflected in the structure of subcategories of modules. In this paper, we extend Gabriel's result (\autoref{IntroBijectionBetweenPrelocSubcatsAndFiltsOfSubobjOfRing}) to an arbitrary locally noetherian scheme $X$ and give a systematic classification of subcategories of the category $\QCoh X$ of quasi-coherent sheaves on $X$.

We deal with the following classes of subcategories.

\begin{Definition}\label{IntroPrelocSubcatAndClosedSubcatAndLocSubcatAndBilocSubcat}
	Let $\mcA$ be a Grothendieck category.
	\begin{enumerate}
		\item\label{IntroPrelocSubcat} A \emph{prelocalizing subcategory} $\mcX$ of $\mcA$ is a full subcategory of $\mcA$ closed under subobjects, quotient objects, and arbitrary direct sums.
		\item\label{IntroLocSubcat} A \emph{localizing subcategory} of $\mcA$ is a prelocalizing subcategory of $\mcA$ closed under extensions.
		\item\label{IntroClosedSubcat} A \emph{closed subcategory} of $\mcA$ is a prelocalizing subcategory of $\mcA$ closed under arbitrary direct products.
		\item\label{IntroBilocSubcat} A \emph{bilocalizing subcategory} of $\mcA$ is a prelocalizing subcategory of $\mcA$ which is both localizing and closed.
	\end{enumerate}
\end{Definition}

It is known that a full subcategory $\mcX$ of a Grothendieck category $\mcA$ is prelocalizing (resp.\ closed) if and only if $\mcX$ is closed under subobjects and quotient objects, and the inclusion functor $\mcX\to\mcA$ has a right adjoint (resp.\ both a right and a left adjoint). See \autoref{CharacterizationOfPrelocSubcat} and \autoref{CharacterizationOfClosedSubcat}.

The notion of closed subcategories can be regarded as the categorical reformulation of closed subschemes of a given scheme. In fact, for every ring $\Lambda$, Rosenberg \cite{Rosenberg1} showed that there is a bijection
\begin{equation*}
	\setwithspace{\text{two-sided ideals of }\Lambda}\to\setwithspace{\text{closed subcategories of }\Mod\Lambda}
\end{equation*}
given by $I\mapsto\conditionalset{M\in\Mod\Lambda}{MI=0}$. It has been shown that the analogous results hold for every noetherian scheme with an ample line bundle (\cite[Theorem~4.1]{Smith}) and for every separated scheme (\cite[Proposition~3.18]{Brandenburg}).

One of the aims of this paper is to classify the closed subcategories of $\QCoh X$ for a locally noetherian scheme $X$. In more generality, we classify the prelocalizing subcategories of $\QCoh X$ by giving an analog of the following famous theorem by Gabriel \cite{Gabriel}.

\begin{Theorem}[{\cite[Lemma~V.2.1]{Gabriel}; \autoref{BijectionBetweenPrelocSubcatsAndFiltsOfSubobjOfRing}}]\label{IntroBijectionBetweenPrelocSubcatsAndFiltsOfSubobjOfRing}
	Let $\Lambda$ be a ring. There is a bijection
	\begin{equation*}
		\setwithspace{\text{prelocalizing subcategories of }\Mod\Lambda}\to\setwithspace{\text{prelocalizing filters of right ideals of }\Lambda}
	\end{equation*}
	given by $\mcY\mapsto\conditionalset{L\subset\Lambda\text{ in }\Mod\Lambda}{\Lambda/L\in\mcY}$.
\end{Theorem}

Note that the prelocalizing filters of right ideals of $\Lambda$ bijectively correspond to the right linear topologies on $\Lambda$ (see \cite[section~VI.4]{Stenstrom}).

For a locally noetherian scheme $X$, there exist too many filters of quasi-coherent subsheaves of $\OX$ compared with the prelocalizing subcategories of $\QCoh X$. Hence we need to consider a suitable class of filters, which we call \emph{local filters} (\autoref{LocFiltOfSubobjs}). By using local filters, we obtain a classification of prelocalizing subcategories, and as a consequence, we deduce classifications of localizing subcategories, closed subcategories, and bilocalizing subcategories.

\begin{Theorem}[\autoref{ClassficationOfPrelocSubcatsForLocNoethSch}, \autoref{ClassificationOfLocSubcatsForLocNoethSch}, \autoref{ClassificationOfClosedSubcatsForLocNoethSch}, \autoref{ClassificationOfBilocSubcatsForLocNoethSch}, \autoref{BijectionBetweenClosedSubcatsAndClosedSubschs}, and \autoref{BijectionBetweenBilocSubcatsAndIdempSubobjAndClopenSubsch}]\label{IntroClassficationOfPrelocSubcatsForLocNoethSch}
	Let $X$ be a locally noetherian scheme. There is a bijection
	\begin{equation*}
		\setwithspace{\text{prelocalizing subcategories of }\QCoh X}\to\setwithspace{\text{local filters of quasi-coherent subsheaves of }\OX}
	\end{equation*}
	given by
	\begin{equation*}
		\mcY\mapsto\conditionalset[\bigg]{I\subset\OX\text{ in }\QCoh X}{\frac{\OX}{I}\in\mcY}.
	\end{equation*}
	This bijection restricts to bijections
	\begin{equation*}
		\setwithspace{\text{localizing subcategories of }\QCoh X}
		\to
		\setwithspace[\bigg]{
			\begin{gathered}
				\text{local filters of quasi-coherent subsheaves of }\OX\\
				\text{ closed under products}
			\end{gathered}
		},
	\end{equation*}
	\begin{equation*}
		\setwithspace{\text{closed subcategories of }\QCoh X}
		\to
		\setwithspace{\text{principal filters of quasi-coherent subsheaves of }\OX},
	\end{equation*}
	and
	\begin{equation*}
		\setwithspace{\text{bilocalizing subcategories of }\QCoh X}
		\to
		\setwithspace[\bigg]{
			\begin{gathered}
				\text{principal filters of quasi-coherent subsheaves of }\OX\\
				\text{ closed under products}
			\end{gathered}
		}.
	\end{equation*}
	
	In particular, there exists a bijection between the closed subcategories of $\QCoh X$ and the closed subschemes of $X$, and it restricts to a bijection between the bilocalizing subcategories of $\QCoh X$ and the subsets of $X$ which are open and closed.
\end{Theorem}

The key of the proof of \autoref{IntroClassficationOfPrelocSubcatsForLocNoethSch} is to reduce the problem to open affine subschemes, and this part is in fact a consequence of the general theory of Grothendieck categories (\autoref{LocOfPrelocSubcat}). The notion of atom spectrum plays a crucial role in this process, and it clarifies the essential properties of the Grothendieck category $\QCoh X$.

The \emph{atom spectrum} $\ASpec\mcA$ of a Grothendieck category $\mcA$ is the set of \emph{atoms} in $\mcA$ which were introduced by Storrer \cite{Storrer} (\autoref{ASpec}). It is regarded as the collection of structural elements of the Grothendieck category in our previous studies \cite{Kanda1,Kanda3,Kanda2}. An atom is a generalization of a prime ideal of a commutative ring. Indeed, for every commutative ring $R$, there exists a canonical bijection between $\ASpec(\Mod R)$ and $\Spec R$ (\autoref{MonoformObjAndASpecOfCommRing}). For a locally noetherian scheme $X$, it is shown in this paper that there exists a canonical bijection between $\ASpec(\QCoh X)$ and the underlying space of $X$ (\autoref{DescriptionOfASpecOfLocNoethSch}). Several fundamental notions of commutative rings and locally noetherian schemes are generalized to Grothendieck categories in terms of atom spectrum as summarized in \autoref{CorrespondingNotions}.

\begin{table}\footnotesize\caption{Corresponding notions on $\ASpec\mcA$, $\Spec R$, and $X$}\label{CorrespondingNotions}
	\begin{tabular}{ccc}
		\toprule
		Grothendieck category $\mcA$ & Commutative ring $R$ & Locally noetherian scheme $X$\\
		\midrule
		Atom spectrum $\ASpec\mcA$ & Prime spectrum $\Spec R$ & Underlying space $\uspX$\\
		Atom $\alpha$ in $\mcA$ & Prime ideal $\mfp$ of $R$ & Point $x\in X$\\
		Associated atoms $\AAss M$ & Associated primes $\Ass M$ & Associated points $\Ass M$\\
		Atom support $\ASupp M$ & Support $\Supp M$ & Support $\Supp M$\\
		Open subsets of $\ASpec\mcA$ & Specialization-closed subsets of $\Spec R$ & Specialization-closed subsets of $X$\\
		$\overline{\set{\alpha}}$ for $\alpha\in\ASpec\mcA$ & $\conditionalset{\mfq\in\Spec R}{\mfq\subset\mfp}$ for $\mfp\in\Spec R$ & $\conditionalset{y\in X}{x\in\overline{\set{y}}}$ for $x\in X$\\
		$\alpha_{1}\leq\alpha_{2}$ & $\mfp_{1}\subset\mfp_{2}$ & $\overline{\set{x_{1}}}\ni x_{2}$\\
		Maximal atoms in $\mcA$ & Maximal ideals of $R$ & Closed points in $X$\\
		Open points in $\ASpec\mcA$ & Maximal ideals of $R$ & Closed points in $X$\\
		Minimal atoms in $\mcA$ & Minimal prime ideals of $R$ & Points in $X$ of height $0$\\
		($=$Closed points in $\ASpec\mcA$) &\\
		Generic point in $\ASpec\mcA$ & Unique maximal ideal of $R$ & Unique closed point in $X$\\
		Injective envelope $E(\alpha)$ & Injective envelope $E(R/\mfp)$ & ${j_{x}}_{*}E(x)$\\
		Residue field $k(\alpha)$ & Residue field $k(\mfp)$ & Residue field $k(x)$\\
		Atomic object $H(\alpha)$ & Residue field $k(\mfp)$ & ${j_{x}}_{*}k(x)$\\
		Localization $\mcA_{\alpha}$ & $\Mod R_{\mfp}$ & $\Mod\OXx$\\
		\bottomrule
	\end{tabular}
\end{table}

In this paper, we also generalize another kind of Gabriel's classification of localizing subcategories. Gabriel \cite{Gabriel} showed that for a noetherian scheme $X$, the localizing subcategories of $\QCoh X$ bijectively correspond to the specialization-closed subsets of the underlying space of $X$ (\cite[Proposition~VI.2.4 (b)]{Gabriel}). This result has been generalized by a number of authors. (For example, \cite{Hovey}, \cite{Krause3}, \cite{GarkushaPrest1}, \cite{GarkushaPrest2}, \cite{Takahashi1}, \cite{Takahashi2}, \cite{Herzog}, \cite{Krause1}, \cite{Kanda1}, and \cite{Kanda2} to some abelian categories. See \cite{GarkushaPrest1} or \cite{Takahashi2} for generalizations to derived categories.) By combining the theory of atom spectrum and the description of the atom spectrum of $\QCoh X$ for a locally noetherian scheme $X$ (\autoref{DescriptionOfASpecOfLocNoethSch}), we obtain the following result.

\begin{Theorem}[\autoref{BijectionBetweenLocSubcatsAndSpecializationClosedSubsetsForLocNoethSch}]\label{IntroBijectionBetweenLocSubcatsAndSpecializationClosedSubsetsForLocNoethSch}
	Let $X$ be a locally noetherian scheme. There is a bijection
	\begin{equation*}
		\setwithspace{\text{localizing subcategories of }\QCoh X}\to\setwithspace{\text{specialization-closed subsets of }X}
	\end{equation*}
	given by $\mcX\mapsto\Supp\mcX$. Its inverse is given by $\Phi\mapsto\Supp^{-1}\Phi$.
\end{Theorem}

This paper is organized as follows. In \autoref{sec:ASpec}, we recall the definition of the atom spectrum and fundamental notions and results on it. \Sectionref{sec:SubcatsAndQuotCats} is devoted to preliminary results on subcategories and quotient categories by localizing subcategories. In \autoref{sec:ASpecOfQuotCatsAndLoc}, we summarize results on the atom spectrum and the localization at an atom. In \autoref{sec:GrothCatsWithEnoughAtoms}, we introduce the class of Grothendieck categories with enough atoms and show that the localizing subcategories are classified in terms of the atom spectrum for a Grothendieck category with enough atoms (\autoref{BijectionBetweenLocSubcatsAndLocSubsets}). In \autoref{sec:ASpecOfLocNoethSch}, we describe the atom spectrum of the Grothendieck category $\QCoh X$ for a locally noetherian scheme $X$ and show that $\QCoh X$ has enough atoms (\autoref{DescriptionOfASpecOfLocNoethSch}). In \autoref{sec:LocOfPrelocSubcatsAndLocSubcats}, we investigate a Grothendieck category $\mcA$ with some properties and relate the prelocalizing subcategories (resp.\ localizing subcategories) of $\mcA$ with the prelocalizing subcategories (resp.\ localizing subcategories) of quotient categories of $\mcA$. For a locally noetherian scheme $X$, the prelocalizing subcategories, the localizing subcategories, the closed subcategories, and the bilocalizing subcategories of $\QCoh X$ are classified in \autoref{sec:ClassificationOfPrelocSubcats}, \autoref{sec:ClassificationOfLocSubcats}, \autoref{sec:ClassificationOfClosedSubcats}, and \autoref{sec:ClassificationOfBilocSubcats}, respectively.

\begin{Remark}\label{IntroRemarkOnTerminology}
	In this paper, we use the words ``prelocalizing'', ``localizing'', and ``bilocalizing'' subcategories in the same way as in \cite{Popescu}. Some authors use different terminology on these subcategories and also on ``closed'' subcategories, which are summarized below. Note that we always work inside a Grothendieck category.
	\begin{enumerate}
		\item Prelocalizing subcategories are often called \emph{weakly closed} subcategories. This terminology was introduced by Van den Bergh \cite{VandenBergh}, and followed by \cite{Smith} and \cite{Pappacena}, for example. Closed subcategories in \cite{VandenBergh} were defined in the same way as we do.
		\item Van den Bergh used different terminology in the preprint version \cite{VandenBerghX}. Weakly closed (resp.\ closed) subcategories in the published version \cite{VandenBergh} were called \emph{closed} (resp.\ \emph{biclosed}) subcategories in \cite{VandenBerghX}. This fits into Gabriel's terminology \cite{Gabriel}.
		\item In the context of torsion theory, such as in \cite[Chapter VI]{Stenstrom}, prelocalizing subcategories, localizing subcategories, and bilocalizing subcategories in this paper are called \emph{hereditary pretorsion class}, \emph{hereditary torsion class}, and \emph{TTF-class} (TTF indicates ``torsion torsionfree''), respectively.
		\item In \cite{Brandenburg}, our prelocalizing subcategories (resp.\ closed subcategories) are called \emph{topologizing} subcategories (resp.\ \emph{reflective topologizing} subcategories). This preprint is aimed at modifying a theory of Rosenberg \cite{Rosenberg2}, and the definition of topologizing subcategories was also changed. In Rosenberg's paper \cite{Rosenberg2}, our prelocalizing subcategories (resp.\ closed subcategories) are called \emph{coreflective topologizing} subcategories (resp.\ \emph{reflective topologizing} subcategories), and they are also called \emph{closed} subcategories (resp.\ \emph{left closed} subcategories) in \cite{Rosenberg1}.
	\end{enumerate}
\end{Remark}

\begin{Conventions}\label{Conventions}
	Throughout this paper, we fix a Grothendieck universe. A set is called \emph{small} if it is an element of the universe. For every category $\mcC$, the collection $\Ob\mcC$ (resp.\ $\Mor\mcC$) of objects (resp.\ morphisms) in $\mcC$ is a set, and $\Hom_{\mcC}(X,Y)$ is supposed to be small for all objects $X$ and $Y$ in $\mcC$. A category $\mcC$ is called \emph{skeletally small} if the set of isomorphism classes of objects in $\mcC$ is in bijection with a small set. The index set of each limit and colimit is assumed to be skeletally small.
	
	Rings, modules over rings, schemes, and sheaves on schemes are assumed to be small. Every ring is associative and has an identity element.
\end{Conventions}

\section{Acknowledgement}
\label{sec:Acknowledgement}

The author would like to express his deep gratitude to Osamu Iyama for his elaborate guidance. The author thanks Mitsuyasu Hashimoto, S.~Paul Smith, and Ryo Takahashi for their valuable comments.

\section{Atom spectrum}
\label{sec:ASpec}

In this section, we recall the definition of the atom spectrum of a Grothendieck category and fundamental results. We start with the definition of a Grothendieck category.

\begin{Definition}\label{GrothCatAndLocNoethGrothCat}\leavevmode
	\begin{enumerate}
		\item\label{GrothCat} An abelian category $\mcA$ is called a \emph{Grothendieck category} if it satisfies the following conditions.
		\begin{enumerate}
			\item $\mcA$ admits arbitrary direct sums (and hence arbitrary direct limits), and for every direct system of short exact sequences in $\mcA$, its direct limit is also a short exact sequence.
			\item $\mcA$ has a generator $G$, that is, every object in $\mcA$ is isomorphic to a quotient object of the direct sum of some copies of $G$.
		\end{enumerate}
		\item\label{LocNoethGrothCat} A Grothendieck category is called \emph{locally noetherian} if it admits a small generating set consisting of noetherian objects.
	\end{enumerate}
\end{Definition}

The exactness of direct limits has the following characterizations.

\begin{Proposition}\label{CharacterizationOfAb5}
	Let $\mcA$ be an abelian category with arbitrary direct sums. Then the following assertions are equivalent.
	\begin{enumerate}
		\item For every direct system of short exact sequences in $\mcA$, its direct limit is also a short exact sequence.
		\item Let $M$ be an object in $\mcA$. For each subobject $L$ of $M$ and each family $\mcN$ of subobjects of $M$ such that every finite subfamily of $\mcN$ has an upper bound in $\mcN$, we have
		\begin{equation*}
			L\cap\sum_{N\in\mcN}N=\sum_{N\in\mcN}(L\cap N).
		\end{equation*}
		\item For every family $\set{M_{\lambda}}_{\lambda\in\Lambda}$ of objects in $\mcA$ and every subobject $L$ of $\bigoplus_{\lambda\in\Lambda}M_{\lambda}$,
		\begin{equation*}
			L=\sum_{\Lambda'\in\mcS}\left(L\cap\bigoplus_{\lambda\in\Lambda'}M_{\lambda}\right),
		\end{equation*}
		where $\mcS$ is the set of finite subsets of $\Lambda$.
	\end{enumerate}
\end{Proposition}

\begin{proof}
	\cite[Theorem~2.8.6]{Popescu}.
\end{proof}

From now on, we deal with a Grothendieck category $\mcA$. The atom spectrum of a Grothendieck category is defined by using monoform objects defined as follows.

\begin{Definition}\label{MonoformObjAndAtomEquiv}\leavevmode
	\begin{enumerate}
		\item\label{MonoformObj} A nonzero object $H$ in $\mcA$ is called \emph{monoform} if for each nonzero subobject $L$ of $H$, no nonzero subobject of $H$ is isomorphic to a subobject of $H/L$.
		\item\label{AtomEquiv} For monoform objects $H_{1}$ and $H_{2}$ in $\mcA$, we say that $H_{1}$ is \emph{atom-equivalent} to $H_{2}$ if there exists a nonzero subobject of $H_{1}$ which is isomorphic to a subobject of $H_{2}$.
	\end{enumerate}
\end{Definition}

We recall the definitions of essential subobjects and uniform objects. These are also important notions in a Grothendieck category and related to monoform objects.

\begin{Definition}\label{EssSubobjAndUniformObj}\leavevmode
	\begin{enumerate}
		\item\label{EssSub} Let $M$ be an object in $\mcA$. A subobject $L$ of $M$ is called \emph{essential} if for every nonzero subobject $L'$ of $M$, we have $L\cap L'\neq 0$.
		\item\label{UniformObj} A nonzero object $U$ in $\mcA$ is called \emph{uniform} if every nonzero subobject of $U$ is essential.
	\end{enumerate}
\end{Definition}

In other words, a nonzero object $U$ in $\mcA$ is uniform if and only if for all nonzero subobjects $L_{1}$ and $L_{2}$ of $U$, we have $L_{1}\cap L_{2}\neq 0$.

It is easy to show that each nonzero subobject of a uniform object is uniform. This type of result also holds for monoform objects.

\begin{Proposition}\label{PropertiesOfMonoformObj}\leavevmode
	\begin{enumerate}
		\item\label{SubOfMonoformObjIsMonoform} Each nonzero subobject of a monoform object is monoform.
		\item\label{MonoformObjIsUniform} Every monoform object is uniform.
		\item\label{NoethObjHasMonoformSub} Every nonzero noetherian object has a monoform subobject.
	\end{enumerate}
\end{Proposition}

\begin{proof}
	\autoref{SubOfMonoformObjIsMonoform} \cite[Proposition~2.2]{Kanda1}.
	
	\autoref{MonoformObjIsUniform} \cite[Proposition~2.6]{Kanda1}.
	
	\autoref{NoethObjHasMonoformSub} \cite[Theorem~2.9]{Kanda1}.
\end{proof}

It follows from \autoref{PropertiesOfMonoformObj} \autoref{MonoformObjIsUniform} that atom equivalence is an equivalence relation on the set of monoform objects in $\mcA$ (\cite[Proposition~2.8]{Kanda1}). The atom spectrum is defined by using this relation.

\begin{Definition}\label{ASpec}
	Let $\mcA$ be a Grothendieck category. Denote by $\ASpec\mcA$ the quotient set of the set of monoform objects in $\mcA$ by atom equivalence. We call it the \emph{atom spectrum} of $\mcA$. Each element of $\ASpec\mcA$ is called an \emph{atom} in $\mcA$. For each monoform object $H$ in $\mcA$, the equivalence class of $H$ is denoted by $\overline{H}$.
\end{Definition}

It is shown in \cite[Proposition~2.7 (2)]{Kanda3} that the atom spectrum $\ASpec\mcA$ of a Grothendieck category $\mcA$ is in bijection with a small set.

The following result shows that the atom spectrum of a Grothendieck category is a generalization of the prime spectrum of a commutative ring.

\begin{Proposition}\label{MonoformObjAndASpecOfCommRing}
	Let $R$ be a commutative ring.
	\begin{enumerate}
		\item\label{MonoformObjAndPrimeIdeal} \textnormal{(\cite[Lemma~1.5]{Storrer})} Let $\mfa$ be an ideal of $R$. Then $R/\mfa$ is a monoform object in $\Mod R$ if and only if $\mfa$ is a prime ideal.
		\item\label{ASpecOfCommRing} \textnormal{(\cite[p.~631]{Storrer})} There is a bijection $\Spec R\to\ASpec(\Mod R)$ given by $\mfp\mapsto\overline{R/\mfp}$.
	\end{enumerate}
\end{Proposition}

We can also generalize the notions of supports and associated primes in commutative ring theory.

\begin{Definition}\label{AAssAndASupp}
	Let $M$ be an object in $\mcA$.
	\begin{enumerate}
		\item\label{AAss} Define the subset $\AAss M$ of $\ASpec\mcA$ by
		\begin{equation*}
			\AAss M=\conditionalset{\alpha\in\ASpec\mcA}{\alpha=\overline{H}\text{ for some monoform subobject }H\text{ of }M}.
		\end{equation*}
		We call each element of $\AAss M$ an \emph{associated atom} of $M$.
		\item\label{ASupp} Define the subset $\ASupp M$ of $\ASpec\mcA$ by
		\begin{equation*}
			\ASupp M=\conditionalset{\alpha\in\ASpec\mcA}{\alpha=\overline{H}\text{ for some monoform subquotient }H\text{ of }M}.
		\end{equation*}
		We call it the \emph{atom support} of $M$.
	\end{enumerate}
\end{Definition}

\begin{Proposition}\label{AAssAndASuppOfCommRing}
	Let $R$ be a commutative ring, and let $M$ be an $R$-module. Then the bijection $\Spec R\to\ASpec(\Mod R)$ in \autoref{MonoformObjAndASpecOfCommRing} \autoref{ASpecOfCommRing} restricts to bijections $\Ass M\to\AAss M$ and $\Supp M\to\ASupp M$.
\end{Proposition}

\begin{proof}
	\cite[Proposition~2.13]{Kanda3}.
\end{proof}

The following results are generalizations of fundamental results in commutative ring theory.

\begin{Proposition}\label{AAssAndASuppAndShortExactSeq}
	Let $0\to L\to M\to N\to 0$ be an exact sequence in $\mcA$. Then
	\begin{equation*}
		\AAss L\subset\AAss M\subset\AAss L\cup\AAss N,
	\end{equation*}
	and
	\begin{equation*}
		\ASupp M=\ASupp L\cup\ASupp N.
	\end{equation*}
\end{Proposition}

\begin{proof}
	\cite[Proposition~3.5]{Kanda1} and \cite[Proposition~3.3]{Kanda1}.
\end{proof}

\begin{Proposition}\label{AAssAndASuppAndDirectSumAndSum}\leavevmode
	\begin{enumerate}
		\item\label{AAssAndASuppAndDirectSum} Let $\set{M_{\lambda}}_{\lambda\in\Lambda}$ be a family of objects in $\mcA$. Then
		\begin{equation*}
			\AAss\bigoplus_{\lambda\in\Lambda}M_{\lambda}=\bigcup_{\lambda\in\Lambda}\AAss M_{\lambda},
		\end{equation*}
		and
		\begin{equation*}
			\ASupp\bigoplus_{\lambda\in\Lambda}M_{\lambda}=\bigcup_{\lambda\in\Lambda}\ASupp M_{\lambda}.
		\end{equation*}
		\item\label{ASuppAndSum} Let $M$ be an object in $\mcA$, and let $\set{L_{\lambda}}_{\lambda\in\Lambda}$ be a family of subobjects of $M$. Then
		\begin{equation*}
			\ASupp\sum_{\lambda\in\Lambda}L_{\lambda}=\bigcup_{\lambda\in\Lambda}\ASupp L_{\lambda}.
		\end{equation*}
	\end{enumerate}
\end{Proposition}

\begin{proof}
	\autoref{AAssAndASuppAndDirectSum} \cite[Proposition~2.12]{Kanda3}.
	
	\autoref{ASuppAndSum} Since we have the canonical epimorphism $\bigoplus_{\lambda\in\Lambda}L_{\lambda}\onto\sum_{\lambda\in\Lambda}L_{\lambda}$ and the inclusion $L_{\mu}\subset\sum_{\lambda\in\Lambda}L_{\lambda}$ for each $\mu\in\Lambda$, we obtain
	\begin{equation*}
		\ASupp L_{\mu}\subset\ASupp\sum_{\lambda\in\Lambda}L_{\lambda}\subset\bigcup_{\lambda\in\Lambda}\ASupp L_{\lambda}
	\end{equation*}
	by \autoref{AAssAndASuppAndDirectSum}. Hence the claim follows.
\end{proof}

Similarly to the case of commutative rings, we have the following results on the associated atoms of uniform objects and essential subobjects.

\begin{Proposition}\label{AAssOfUniformObjAndEssSubobj}\leavevmode
	\begin{enumerate}
		\item\label{AAssOfUniformObj} Let $U$ be a uniform object in $\mcA$. Then $\AAss U$ consists of at most one element. In particular, for every monoform object $H$ in $\mcA$, we have $\AAss H=\set{\overline{H}}$.
		\item\label{AAssOfEssSubobj} Let $M$ be an object in $\mcA$, and let $L$ be an essential subobject of $M$. Then $\AAss L=\AAss M$.
	\end{enumerate}
\end{Proposition}

\begin{proof}
	\autoref{AAssOfUniformObj} \cite[Proposition~2.15 (1)]{Kanda3}.
	
	\autoref{AAssOfEssSubobj} \cite[Proposition~2.16]{Kanda3}.
\end{proof}

We introduce a topology on the atom spectrum.

\begin{Definition}\label{LocSubset}
	We call a subset $\Phi$ of $\ASpec\mcA$ a \emph{localizing subset} if there exists an object $M$ in $\mcA$ such that $\Phi=\ASupp M$.
\end{Definition}

\begin{Proposition}\label{LocSubsetsFormOpenSubsets}
	The set of localizing subsets of $\ASpec\mcA$ satisfies the axioms of open subsets of $\ASpec\mcA$.
\end{Proposition}

\begin{proof}
	\cite[Proposition~3.8]{Kanda1}.
\end{proof}

We call the topology on $\ASpec\mcA$ defined by the set of localizing subsets of $\ASpec\mcA$ the \emph{localizing topology}. Throughout this paper, we regard $\ASpec\mcA$ as a topological space in this way. For a commutative ring $R$, the localizing subsets of $\ASpec(\Mod R)$ define a different topology from the Zariski topology on $\Spec R$. Recall that a subset $\Phi$ of $\Spec R$ is said to be \emph{closed under specialization} if for every $\mfp,\mfq\in\Spec R$, the conditions $\mfp\in\Phi$ and $\mfp\subset\mfq$ imply $\mfq\in\Phi$.

\begin{Proposition}\label{LocSubsetOfCommRing}
	Let $R$ be a commutative ring, and let $\Phi$ be a subset of $\Spec R$. Then the corresponding subset
	\begin{equation*}
		\conditionalset[\Bigg]{\overline{\left(\frac{R}{\mfp}\right)}\in\ASpec(\Mod R)}{\mfp\in\Phi}
	\end{equation*}
	of $\ASpec(\Mod R)$ is localizing if and only if $\Phi$ is closed under specialization.
\end{Proposition}

\begin{proof}
	\cite[Proposition~7.2 (2)]{Kanda1}.
\end{proof}

For each $\alpha\in\ASpec\mcA$, let $\Lambda(\alpha)$ be the topological closure of $\set{\alpha}$ in $\ASpec\mcA$. We introduce a partial order on the atom spectrum.

\begin{Definition}\label{SpecializationOrderOnASpec}
	For $\alpha,\beta\in\ASpec\mcA$, we write $\alpha\leq\beta$ if $\alpha\in\Lambda(\beta)$.
\end{Definition}

The relation $\leq$ is called the \emph{specialization order} on the topological space $\ASpec\mcA$ with respect to the localizing topology. This is in fact a partial order on $\ASpec\mcA$ since the topological space $\ASpec\mcA$ is a Kolmogorov space (\cite[Proposition~3.5]{Kanda3}).

By definition, $\Lambda(\beta)=\conditionalset{\alpha\in\ASpec\mcA}{\alpha\leq\beta}$ for each $\beta\in\ASpec\mcA$. The partial order has the following descriptions.

\begin{Proposition}\label{CharacterizationOfSpecializationOrder}
	Let $\alpha,\beta\in\ASpec\mcA$. Then the following assertions are equivalent.
	\begin{enumerate}
		\item $\alpha\leq\beta$, that is, $\alpha\in\Lambda(\beta)$.
		\item For every object $M$ in $\mcA$, the condition $\alpha\in\ASupp M$ implies $\beta\in\ASupp M$.
		\item For every monoform object $H$ in $\mcA$ with $\overline{H}=\alpha$, we have $\beta\in\ASupp H$.
	\end{enumerate}
\end{Proposition}

\begin{proof}
	\cite[Proposition~4.2]{Kanda3}.
\end{proof}

The following result claims that the partial order $\leq$ on $\ASpec\mcA$ is a generalization of the inclusion relation between prime ideals of a commutative ring.

\begin{Proposition}\label{SpecializationOrderOfCommRing}
	Let $R$ be a commutative ring and $\mfp,\mfq\in\Spec R$. Then $\overline{R/\mfp}\leq\overline{R/\mfq}$ in $\ASpec(\Mod R)$ if and only if $\mfp\subset\mfq$. In other words, the bijection $\Spec R\to\ASpec(\Mod R)$ in \autoref{MonoformObjAndASpecOfCommRing} \autoref{ASpecOfCommRing} is an isomorphism between the partially ordered sets $(\Spec R,\subset)$ and $(\ASpec(\Mod R),\leq)$.
\end{Proposition}

\begin{proof}
	\cite[Proposition~4.3]{Kanda3}.
\end{proof}

\section{Subcategories and quotient categories}
\label{sec:SubcatsAndQuotCats}

In this section, we show preliminary results on subcategories and quotient categories of a Grothendieck category $\mcA$. We start with defining some classes of subcategories, which are the main objects in this paper.

\begin{Definition}\label{SubcatClosedUnderExtAndLocSubcat}\leavevmode
	\begin{enumerate}
		\item\label{ExtOfFullSubcats} For full subcategories $\mcX_{1}$ and $\mcX_{2}$ of $\mcA$, we denote by $\mcX_{1}*\mcX_{2}$ the full subcategory of $\mcA$ consisting of all objects $M$ admitting an exact sequence
		\begin{equation*}
			0\to M_{1}\to M\to M_{2}\to 0
		\end{equation*}
		in $\mcA$, where $M_{i}$ belongs to $\mcX_{i}$ for each $i=1,2$.
		\item\label{SubcatClosedUnderExt} We say that a full subcategory $\mcX$ of $\mcA$ is \emph{closed under extension} if $\mcX*\mcX\subset\mcX$, that is, for every exact sequence $0\to L\to M\to N\to 0$ in $\mcA$, the condition $L,N\in\mcX$ implies $M\in\mcX$.
		\item\label{PrelocSubcat} A full subcategory $\mcX$ of $\mcA$ is called a \emph{prelocalizing subcategory} (or \emph{weakly closed subcategory} in \cite{VandenBergh}) of $\mcA$ if $\mcX$ is closed under subobjects, quotient objects, and arbitrary direct sums.
		\item\label{LocSubcat} A prelocalizing subcategory $\mcX$ of $\mcA$ is called a \emph{localizing subcategory} of $\mcA$ if $\mcX$ is also closed under extensions.
		\item\label{PrelocClosureAndLocClosure} For a full subcategory $\mcX$ of $\mcA$, denote by $\preloc{\mcX}$ (resp.\ $\loc{\mcX}$) the smallest prelocalizing (resp.\ localizing) subcategory of $\mcA$ containing $\mcX$. For an object $M$ in $\mcA$, let $\preloc{M}=\preloc{\set{M}}$ and $\loc{M}=\loc{\set{M}}$.
	\end{enumerate}
\end{Definition}

\begin{Proposition}\label{PropertiesOfExt}\leavevmode
	\begin{enumerate}
		\item\label{ExtIsAssociative} Let $\mcX_{1}$, $\mcX_{2}$, and $\mcX_{3}$ be full subcategories of $\mcA$. Then
		\begin{equation*}
			(\mcX_{1}*\mcX_{2})*\mcX_{3}=\mcX_{1}*(\mcX_{2}*\mcX_{3}).
		\end{equation*}
		\item\label{ExtOfPrelocSubcatsIsPreloc} Let $\mcX_{1}$ and $\mcX_{2}$ be prelocalizing subcategories of $\mcA$. Then $\mcX_{1}*\mcX_{2}$ is also a prelocalizing subcategory of $\mcA$.
	\end{enumerate}
\end{Proposition}

\begin{proof}
	\autoref{ExtIsAssociative} \cite[Proposition~2.4 (2)]{Kanda1}.
	
	\autoref{ExtOfPrelocSubcatsIsPreloc} \cite[Lemma~4.8.11]{Popescu}.
\end{proof}

Prelocalizing subcategories are characterized as follows.

\begin{Proposition}\label{CharacterizationOfPrelocSubcat}
	Let $\mcA$ be a Grothendieck category (or more generally, an abelian category admitting arbitrary direct sums), and let $\mcX$ be a full subcategory of $\mcA$ closed under subobjects and quotient objects. Then the following assertions are equivalent.
	\begin{enumerate}
		\item\label{CharacterizationOfPrelocSubcat:PrelocSubcat} $\mcX$ is closed under arbitrary direct sums, that is, $\mcX$ is a prelocalizing subcategory of $\mcA$.
		\item\label{CharacterizationOfPrelocSubcat:RightAdjoint} The inclusion functor $\mcX\into\mcA$ has a right adjoint.
		\item\label{CharacterizationOfPrelocSubcat:LargestSubobj} For each object $M$ in $\mcA$, there exists a largest subobject $L$ of $M$ which belongs to $\mcX$.
	\end{enumerate}
\end{Proposition}

\begin{proof}
	Assume \autoref{CharacterizationOfPrelocSubcat:LargestSubobj}. Then the functor $\mcA\to\mcX$ which sends each object $M$ to its largest subobject belonging to $\mcX$ and each morphism to the induced morphism is a right adjoint of the inclusion functor $\mcX\into\mcA$.
	
	The dual statement of \autoref{CharacterizationOfPrelocSubcat:RightAdjoint}$\Rightarrow$\autoref{CharacterizationOfPrelocSubcat:PrelocSubcat} is essentially shown in \cite[Proposition~X.1.2]{Stenstrom}.
	
	\autoref{CharacterizationOfPrelocSubcat:PrelocSubcat}$\Rightarrow$\autoref{CharacterizationOfPrelocSubcat:LargestSubobj} follows from the next remark.
\end{proof}

\begin{Remark}\label{LargestSubobjWrtFullSubcatClosedUnderQuotsAndDirectSums}
	Let $\mcX$ be a full subcategory of $\mcA$ closed under quotient objects and arbitrary direct sums, and let $M$ be an object in $\mcA$. Since the sum $L=\sum_{\lambda\in\Lambda}L_{\lambda}$ of all subobjects of $M$ which belong to $\mcX$ is a quotient object of the direct sum $\bigoplus_{\lambda\in\Lambda}L_{\lambda}$, the subobject $L$ of $M$ also belongs to $\mcX$. Hence $L$ is the largest subobject of $M$ which belongs to $\mcX$.
\end{Remark}

The operation in \autoref{LargestSubobjWrtFullSubcatClosedUnderQuotsAndDirectSums} of taking the subobject $L$ from $M$ is used throughout this paper. The following result shows that this operation commutes with taking arbitrary direct sums.

\begin{Proposition}\label{LargestSubobjOfDirectSumWrtPrelocSubcat}
	Let $\mcA$ be a Grothendieck category, and let $\mcX$ be a full subcategory of $\mcA$ closed under quotient objects and arbitrary direct sums. Let $\set{M_{\lambda}}_{\lambda\in\Lambda}$ be a family of objects in $\mcA$, and take $L_{\lambda}$ to be the largest subobject of $M_{\lambda}$ which belongs to $\mcX$ for each $\lambda\in\Lambda$. Then $\bigoplus_{\lambda\in\Lambda}L_{\lambda}$ is the largest subobject of $\bigoplus_{\lambda\in\Lambda}M_{\lambda}$ which belongs to $\mcX$.
\end{Proposition}

\begin{proof}
	Let $N$ be the largest subobject of $\bigoplus_{\lambda\in\Lambda}M_{\lambda}$ which belongs to $\mcX$. It suffices to show that $N\subset\bigoplus_{\lambda\in\Lambda}L_{\lambda}$.
	
	We show the claim in the case where $\Lambda=\set{1,\ldots,n}$ for some $n\in\mbZ_{\geq 1}$. Let $\pi_{i}\colon M_{1}\oplus\cdots\oplus M_{n}\onto M_{i}$ be the projection for each $i\in\set{1,\ldots,n}$. Since $\pi_{i}(N)$ is a quotient object of $N$, it belongs to $\mcX$. By the maximality of $L_{i}$, we have $\pi_{i}(N)\subset L_{i}$. Hence
	\begin{equation*}
		N\subset\pi_{1}(N)\oplus\cdots\oplus\pi_{n}(N)\subset L_{1}\oplus\cdots\oplus L_{n}
	\end{equation*}
	as subobjects of $M_{1}\oplus\cdots\oplus M_{n}$.
	
	In the general case, let $\mcS$ be the set of finite subsets of $\Lambda$. Then by \autoref{CharacterizationOfAb5},
	\begin{equation*}
		N=\sum_{\Lambda'\in\mcS}\left(N\cap\bigoplus_{\lambda\in\Lambda'}M_{\lambda}\right)\subset\sum_{\Lambda'\in\mcS}\bigoplus_{\lambda\in\Lambda'}L_{\lambda}=\bigoplus_{\lambda\in\Lambda}L_{\lambda}.\qedhere
	\end{equation*}
\end{proof}

For a localizing subcategory $\mcX$ of $\mcA$, we have the \emph{quotient category} $\mcA/\mcX$ of $\mcA$ by $\mcX$. It is a Grothendieck category together with a canonical (covariant) functor $\mcA\to\mcA/\mcX$ (\cite[Corollary~4.6.2]{Popescu}). We refer the reader to \cite[Definition~5.2]{Kanda3} for the explicit definition of the quotient category. Instead, we state a universal property of the quotient category.

\begin{Theorem}\label{FundamentalAndUniversalPropertyOfQuotCat}
	Let $\mcA$ be a Grothendieck category, and let $\mcX$ be a localizing subcategory of $\mcA$. The canonical functor is denoted by $F\colon\mcA\to\mcA/\mcX$.
	\begin{enumerate}
		\item\label{FundamentalPropertyOfQuotCat} The functor $F\colon\mcA\to\mcA/\mcX$ is exact and has a right adjoint $\mcA/\mcX\to\mcA$. For every object $M$ in $\mcA$, we have $F(M)=0$ if and only if $M$ belongs to $\mcX$.
		\item\label{UniversalPropertyOfQuotCat} Let $\mcB$ be an abelian category together with an exact functor $Q\colon\mcA\to\mcB$ with $Q(M)=0$ for each object $M$ in $\mcX$. Then there exists a unique functor $\overline{Q}\colon\mcA/\mcX\to\mcB$ such that $\overline{Q}F=Q$. Moreover, the functor $\overline{Q}$ is exact.
	\end{enumerate}
\end{Theorem}

\begin{proof}
	\autoref{FundamentalPropertyOfQuotCat} \cite[Proposition~4.6.3]{Popescu}, \cite[Theorem~4.3.8]{Popescu}, and \cite[Lemma~4.3.4]{Popescu}.
	
	\autoref{UniversalPropertyOfQuotCat} \cite[Corollary~4.3.11]{Popescu} and \cite[Corollary~4.3.12]{Popescu}.
\end{proof}

Every object $M$ in a Grothendieck category $\mcA$ has an \emph{injective envelope} $E(M)$ (\cite[Theorem~II.6.2]{Gabriel}, see also \cite[Theorem~3.10.10]{Popescu}). By definition, the object $M$ is an essential subobject of the injective object $E(M)$. The object $E(M)$ is also denoted by $E_{\mcA}(M)$ in order to specify the category explicitly.

Let $\mcX$ be a localizing subcategory of $\mcA$. An object $M$ in $\mcA$ is called $\mcX$\emph{-torsionfree} if $M$ has no nonzero subobject belonging to $\mcX$. Note that every subobject of an $\mcX$-torsionfree object is $\mcX$-torsionfree.

\begin{Proposition}\label{TorsionfreenessOfQuotObj}
	Let $\mcX$ be a localizing subcategory of $\mcA$. Let $M$ be an object in $\mcA$, and let $L$ be the largest subobject of $M$ which belongs to $\mcX$. Then $M/L$ is $\mcX$-torsionfree.
\end{Proposition}

\begin{proof}
	Assume that $M/L$ is not $\mcX$-torsionfree. Then there exists a subobject $L'$ of $M$ such that $L\subsetneq L'$, and $L'/L$ belongs to $\mcX$. The subobject $L'$ of $M$ also belongs to $\mcX$. This contradicts the maximality of $L$.
\end{proof}

For an object $M$ in $\mcA$, it is also important to consider the torsionfreeness of $E(M)/M$.

\begin{Proposition}\label{ExactSeqAndTorsionfreeness}
	Let $\mcX$ be a localizing subcategory of $\mcA$, and let
	\begin{equation*}
		0\to L\to M\to N\to 0
	\end{equation*}
	be an exact sequence in $\mcA$. If $M$ and $E(L)/L$ are $\mcX$-torsionfree, then $N$ is $\mcX$-torsionfree.
\end{Proposition}

\begin{proof}
	This can be shown similarly to the proof of \cite[Proposition~4.5.5]{Popescu}.
\end{proof}

We state important properties of the canonical functor to a quotient category and its right adjoint by using the notion of torsionfreeness.

\begin{Proposition}\label{PropertiesOfQuotCat}
	Let $\mcX$ be a localizing subcategory of $\mcA$. Denote the canonical functor by $F\colon\mcA\to\mcA/\mcX$ and its right adjoint by $G\colon\mcA/\mcX\to\mcA$.
	\begin{enumerate}
		\item\label{QuotFunctorIsSurj} $F$ is surjective, that is, each object in $\mcA/\mcX$ is of the form $F(M)$, where $M$ is some object in $\mcA$.
		\item\label{CounitIsIso} The counit morphism $\varepsilon\colon FG\to 1_{\mcA/\mcX}$ is an isomorphism. Hence $G$ is fully faithful.
		\item\label{DescriptionOfUnit} Let $\eta\colon 1_{\mcA}\to GF$ be the unit morphism. Then for each object $M$ in $\mcA$, the subobject $\Ker\eta_{M}$ of $M$ is the largest subobject belonging to $\mcX$, the subobject $\Im\eta_{M}$ of $GF(M)$ is essential, and $\Cok\eta_{M}$ belongs to $\mcX$. The objects $GF(M)$ and $E(GF(M))/GF(M)$ are $\mcX$-torsionfree.
		\item\label{TorsionfreenessOfClosedObj} Let $M'$ be an object in $\mcA/\mcX$. Then $G(M')$ and $E(G(M'))/G(M')$ are $\mcX$-torsionfree.
	\end{enumerate}
\end{Proposition}

\begin{proof}
	\autoref{QuotFunctorIsSurj} This is obvious from the definition of the canonical functor $F$. It also follows from \autoref{FundamentalAndUniversalPropertyOfQuotCat}.
	
	\autoref{CounitIsIso} \cite[Proposition~4.4.3 (1)]{Popescu}.
	
	\autoref{DescriptionOfUnit} This follows from \cite[Proposition~4.4.3 (2)]{Popescu} and the proof of \cite[Proposition~4.4.5]{Popescu}.
	
	\autoref{TorsionfreenessOfClosedObj} This follows from \autoref{QuotFunctorIsSurj} and \autoref{DescriptionOfUnit}.
\end{proof}

The next result is necessary to describe subobjects of an object in a quotient category.

\begin{Proposition}\label{SubobjOfObjInQuotCat}
	Let $\mcX$ be a localizing subcategory of $\mcA$. Denote the canonical functor by $F\colon\mcA\to\mcA/\mcX$ and its right adjoint by $G\colon\mcA/\mcX\to\mcA$. Let $M$ be an object in $\mcA$. For each subobject $L'$ of $F(M)$, there exists a largest subobject $L$ of $M$ satisfying $F(L)\subset L'$ as a subobject of $F(M)$. Moreover, it holds that $F(L)=L'$, and the quotient object $M/L$ is $\mcX$-torsionfree. The quotient object $F(M)/L'$ of $F(M)$ is equal to $F(M/L)$.
\end{Proposition}

\begin{proof}
	Since $G$ is left exact, the object $G(L')$ can be regarded as a subobject of $GF(M)$. Let $\eta\colon 1_{\mcA}\to GF$ be the unit morphism. There is a commutative diagram
	\begin{equation*}
		\definehookarrow
		\xymatrix{
			0\ar[r] & \eta_{M}^{-1}(G(L'))\ar[d]\ar[r] & M\ar[d]^-{\eta_{M}}\ar[r] & \dfrac{M}{\eta_{M}^{-1}(G(L'))}\ar@{^{ (}->}[d]\ar[r] & 0\\
			0\ar[r] & G(L')\ar[r] & GF(M)\ar[r] & \dfrac{GF(M)}{G(L')}\ar[r] & 0
		}.
	\end{equation*}
	By applying $F$ to this diagram, we obtain the commutative diagram
	\begin{equation*}
		\definehookarrow
		\xymatrix{
			0\ar[r] & F(\eta_{M}^{-1}(G(L')))\ar[d]^-{\cong}\ar[r] & F(M)\ar[d]^-{\cong}\ar[r] & F\left(\dfrac{M}{\eta_{M}^{-1}(G(L'))}\right)\ar[d]^-{\cong}\ar[r] & 0\\
			0\ar[r] & FG(L')\ar[d]^-{\cong}\ar[r] & FGF(M)\ar[d]^-{\cong}\ar[r] & F\left(\dfrac{GF(M)}{G(L')}\right)\ar@{=}[d]\ar[r] & 0\\
			0\ar[r] & L'\ar[r] & F(M)\ar[r] & F\left(\dfrac{GF(M)}{G(L')}\right)\ar[r] & 0
		}
	\end{equation*}
	by \autoref{PropertiesOfQuotCat} \autoref{CounitIsIso} and \autoref{PropertiesOfQuotCat} \autoref{DescriptionOfUnit}. Hence the subobject $L:=\eta_{M}^{-1}(G(L'))$ of $M$ satisfies $F(L)=L'$, and $F(M)/L'=F(M/L)$. By \autoref{ExactSeqAndTorsionfreeness}, the object $GF(M)/G(L')$ is $\mcX$-torsionfree, and hence $M/L$ is also $\mcX$-torsionfree.
	
	Let $\widetilde{L}$ be a subobject of $M$ such that $F(\widetilde{L})\subset L'$. Since we have the commutative diagram
	\begin{equation*}
		\definehookarrow
		\xymatrix{
			\widetilde{L}\ar@{^{ (}->}[r]\ar[d]^-{\eta_{\widetilde{L}}} & M\ar[d]^-{\eta_{M}}\\
			GF(\widetilde{L})\ar@{^{ (}->}[r] & GF(M)
		},
	\end{equation*}
	it holds that $\eta_{M}(\widetilde{L})\subset GF(\widetilde{L})$. Therefore
	\begin{equation*}
		\widetilde{L}\subset\eta_{M}^{-1}(GF(\widetilde{L}))\subset\eta_{M}^{-1}(G(L'))=L.\qedhere
	\end{equation*}
\end{proof}

Several properties of objects are preserved by the canonical functor to a quotient category and its right adjoint as in the following results.

\begin{Proposition}\label{SubobjAndInjectivityInQuotCat}
	Let $\mcX$ be a localizing subcategory of $\mcA$. Denote the canonical functor by $F\colon\mcA\to\mcA/\mcX$ and its right adjoint by $G\colon\mcA/\mcX\to\mcA$.
	\begin{enumerate}
		\item\label{ImageOfEssSubobjInQuotCat} Let $M'$ be an object in $\mcA/\mcX$, and let $L'$ be an essential subobject of $M'$. Then $G(L')$ is an essential subobject of $G(M')$.
		\item\label{ImageOfUniformObjInQuotCat} Let $U'$ be a uniform object in $\mcA/\mcX$. Then $G(U')$ is a uniform object in $\mcA$.
		\item\label{ImageOfMonoformObjInQuotCat} Let $H'$ be a monoform object in $\mcA/\mcX$. Then $G(H')$ is a monoform object in $\mcA$.
		\item\label{ImageOfInjObjInQuotCat} Let $I'$ be an injective object in $\mcA/\mcX$. Then $G(I')$ is an injective object in $\mcA$.
		\item\label{ImageOfIndecObjInQuotCat} Let $M'$ be an indecomposable object in $\mcA/\mcX$. Then $G(M')$ is an indecomposable object in $\mcA$.
	\end{enumerate}
\end{Proposition}

\begin{proof}
	\autoref{ImageOfEssSubobjInQuotCat} \cite[Corollary~4.4.7]{Popescu}.
	
	\autoref{ImageOfUniformObjInQuotCat} Let $L$ be a nonzero subobject of $G(U')$. We have a commutative diagram
	\begin{equation*}
		\definehookarrow
		\xymatrix{
			L\ar@{^{ (}->}[r]\ar[d] & G(U')\ar[d]\\
			GF(L)\ar[r] & GFG(U')
		},
	\end{equation*}
	and the morphism $G(U')\to GFG(U')$ is an isomorphism by \autoref{PropertiesOfQuotCat} \autoref{CounitIsIso}. Hence the morphism $L\to GF(L)$ is a monomorphism, and in particular $F(L)$ is a nonzero subobject of $FG(U')\cong U'$. By the uniformness of $U'$ and \autoref{ImageOfEssSubobjInQuotCat}, $GF(L)$ is an essential subobject of $GFG(U')$. Since $L$ is essential as a subobject of $GF(L)$ by \autoref{PropertiesOfQuotCat} \autoref{DescriptionOfUnit}, $L$ is an essential subobject of $G(U')$.
	
	\autoref{ImageOfMonoformObjInQuotCat} \cite[Lemma~5.14 (1)]{Kanda3}.
	
	\autoref{ImageOfInjObjInQuotCat} \cite[Corollary~4.4.7]{Popescu}.
	
	\autoref{ImageOfIndecObjInQuotCat} This follows from \autoref{PropertiesOfQuotCat} \autoref{CounitIsIso}.
\end{proof}

\begin{Proposition}\label{PropertiesOfImageOfObj}
	Let $\mcX$ be a localizing subcategory of $\mcA$. Denote the canonical functor by $F\colon\mcA\to\mcA/\mcX$ and its right adjoint by $G\colon\mcA/\mcX\to\mcA$.
	\begin{enumerate}
		\item\label{EssentialityOfImageOfObj} Let $M$ be an $\mcX$-torsionfree object in $\mcA$, and let $L$ be an essential subobject of $M$. Then $F(L)$ is an essential subobject of $F(M)$.
		\item\label{UniformnessOfImageOfObj} Let $U$ be a uniform $\mcX$-torsionfree object in $\mcA$. Then $F(U)$ is a uniform object in $\mcA/\mcX$.
		\item\label{MonoformnessOfImageOfObj} Let $H$ be a monoform $\mcX$-torsionfree object in $\mcA$. Then $F(H)$ is a monoform object in $\mcA/\mcX$.
		\item\label{InjectivityOfImageOfObj} Let $I$ be an injective $\mcX$-torsionfree object in $\mcA$. Then $F(I)$ is an injective object in $\mcA/\mcX$.
	\end{enumerate}
\end{Proposition}

\begin{proof}
	\autoref{EssentialityOfImageOfObj} \cite[Lemma~4.4.6 (3)]{Popescu}.
	
	\autoref{UniformnessOfImageOfObj} This follows from \autoref{SubobjOfObjInQuotCat} and \autoref{EssentialityOfImageOfObj}.
	
	\autoref{MonoformnessOfImageOfObj} \cite[Lemma~5.14 (2)]{Kanda3}.
	
	\autoref{InjectivityOfImageOfObj} \cite[Lemma~4.5.1 (2)]{Popescu}.
\end{proof}

The prelocalizing subcategories of $\mcA$ and those of quotient categories are related by the following operations.

\begin{Proposition}\label{QuotOfPrelocSubcat}
	Let $\mcX$ be a localizing subcategory of $\mcA$. Denote the canonical functor by $F\colon\mcA\to\mcA/\mcX$ and its right adjoint by $G\colon\mcA/\mcX\to\mcA$.
	\begin{enumerate}
		\item\label{InverseImageOfPrelocSubcatIsPrelocSubcat} For each prelocalizing subcategory $\mcY'$ of $\mcA/\mcX$, the full subcategory
		\begin{equation*}
			F^{-1}(\mcY'):=\conditionalset{M\in\mcA}{F(M)\in\mcY'}
		\end{equation*}
		of $\mcA$ is a prelocalizing subcategory, and
		\begin{equation*}
			\mcX*F^{-1}(\mcY')*\mcX=F^{-1}(\mcY').
		\end{equation*}
		\item\label{QuotOfPrelocSubcatIsPrelocSubcat} For each prelocalizing subcategory $\mcY$ of $\mcA$, the full subcategory
	\begin{equation*}
		F(\mcY):=\conditionalset[\bigg]{N\in\frac{\mcA}{\mcX}}{N\cong F(M)\text{ for some }M\in\mcY}
	\end{equation*}
	of $\mcA/\mcX$ is a prelocalizing subcategory.
		\item\label{QuotOfExtOfSubcatsIsExtOfQuotOfSubcats} Let $\mcY_{1}$ and $\mcY_{2}$ be prelocalizing subcategories of $\mcA$. Then
		\begin{equation*}
			F(\mcY_{1}*\mcX*\mcY_{2})=F(\mcY_{1})*F(\mcY_{2}).
		\end{equation*}
	\end{enumerate}
\end{Proposition}

\begin{proof}
	\autoref{InverseImageOfPrelocSubcatIsPrelocSubcat} Since $F$ is exact and commutes with arbitrary direct sums, the full subcategory $F^{-1}(\mcY')$ is a prelocalizing subcategory. The inclusion $F^{-1}(\mcY')\subset\mcX*F^{-1}(\mcY')*\mcX$ is obvious. By \autoref{FundamentalAndUniversalPropertyOfQuotCat} \autoref{FundamentalPropertyOfQuotCat},
	\begin{equation*}
		F(\mcX*F^{-1}(\mcY')*\mcX)\subset F(\mcX)*F(F^{-1}(\mcY'))*F(\mcX)\subset\mcY'.
	\end{equation*}
	Hence $\mcX*F^{-1}(\mcY')*\mcX\subset F^{-1}(\mcY')$.
	
	\autoref{QuotOfPrelocSubcatIsPrelocSubcat} By \autoref{SubobjOfObjInQuotCat}, the full subcategory $F(\mcY)$ of $\mcA/\mcX$ is closed under subobjects and quotient objects. It is also closed under arbitrary direct sums since $F$ commutes with arbitrary direct sums.
	
	\autoref{QuotOfExtOfSubcatsIsExtOfQuotOfSubcats} Since $F$ is exact, $F(\mcY_{1}*\mcX*\mcY_{2})\subset F(\mcY_{1})*F(\mcY_{2})$ by \autoref{FundamentalAndUniversalPropertyOfQuotCat} \autoref{FundamentalPropertyOfQuotCat}. Let $M'$ be an object in $\mcA/\mcX$ which belongs to $F(\mcY_{1})*F(\mcY_{2})$. Then there exists an exact sequence
	\begin{equation*}
		0\to F(M_{1})\to M'\to F(M_{2})\to 0
	\end{equation*}
	where $M_{i}$ is an object in $\mcA$ which belongs to $\mcY_{i}$ for each $i=1,2$. Since $G$ is left exact, we have the exact sequence
	\begin{equation*}
		0\to GF(M_{1})\to G(M')\to GF(M_{2}).
	\end{equation*}
	Let $\eta\colon 1_{\mcA}\to GF$ be the unit morphism, and let $B$ be the image of the morphism $G(M')\to GF(M_{2})$. Then we obtain a commutative diagram
	\begin{equation*}
		\definehookarrow
		\xymatrix{
			0\ar[r] & GF(M_{1})\ar@{=}[d]\ar[r] & M\ar@{^{ (}->}[d]\ar[r] & B\cap\Im\eta_{M_{2}}\ar@{^{ (}->}[d]\ar[r] & 0\\
			0\ar[r] & GF(M_{1})\ar[r] & G(M')\ar[r] & B\ar[r] & 0
		},
	\end{equation*}
	where $M$ is an object in $\mcA$. Let $N$ be the cokernel of the composite $\Im\eta_{M_{1}}\into GF(M_{1})\into G(M')$. There is a commutative diagram
	\begin{equation*}
		\definehookarrow
		\xymatrix{
			0\ar[r] & \Im\eta_{M_{1}}\ar@{^{ (}->}[d]\ar[r] & M\ar@{=}[d]\ar[r] & N\ar@{>>}[d]\ar[r] & 0\\
			0\ar[r] & GF(M_{1})\ar[r] & M\ar[r] & B\cap\Im\eta_{M_{2}}\ar[r] & 0
		}.
	\end{equation*}
	By the snake lemma, we have an exact sequence
	\begin{equation*}
		0\to\Cok\eta_{M_{1}}\to N\to B\cap\Im\eta_{M_{2}}\to 0.
	\end{equation*}
	By \autoref{PropertiesOfQuotCat} \autoref{DescriptionOfUnit}, the object $\Cok\eta_{M_{i}}$ belongs to $\mcX$ for each $i=1,2$. Hence $F(\Cok\eta_{M_{1}})=0$, and
	\begin{equation*}
		F\left(\frac{B}{B\cap\Im\eta_{M_{2}}}\right)\cong F\left(\frac{B+\Im\eta_{M_{2}}}{\Im\eta_{M_{2}}}\right)\subset F\left(\frac{GF(M_{2})}{\Im\eta_{M_{2}}}\right)=0.
	\end{equation*}
	By applying $F$ to the morphisms $B\cap\Im\eta_{M_{2}}\into B$ and $\Im\eta_{M_{1}}\into GF(M_{1})$, we obtain $F(B\cap\Im\eta_{M_{2}})\isoto F(B)$ and $F(\Im\eta_{M_{1}})\isoto FGF(M_{1})\isoto F(M_{1})$. Hence we have the commutative diagram
	\begin{equation*}
		\definehookarrow
		\xymatrix{
			0\ar[r] & F(\Im\eta_{M_{1}})\ar[d]^-{\cong}\ar[r] & F(M)\ar@{=}[d]\ar[r] & F(N)\ar[d]^-{\cong}\ar[r] & 0\\
			0\ar[r] & FGF(M_{1})\ar@{=}[d]\ar[r] & F(M)\ar[d]^-{\cong}\ar[r] & F(B\cap\Im\eta_{M_{2}})\ar[d]^-{\cong}\ar[r] & 0\\
			0\ar[r] & FGF(M_{1})\ar[d]^-{\cong}\ar[r] & FG(M')\ar[d]^-{\cong}\ar[r] & F(B)\ar[d]^-{\cong}\ar[r] & 0\\
			0\ar[r] & F(M_{1})\ar[r] & M'\ar[r] & F(M_{2})\ar[r] & 0
		}.
	\end{equation*}
	For each $i=1,2$, the quotient object $\Im\eta_{M_{i}}$ of $M_{i}$ belongs to $\mcY_{i}$, and hence $N$ belongs to $\mcX*\mcY_{2}$. Therefore $M'$ belongs to $F(\mcY_{1}*\mcX*\mcY_{2})$.
\end{proof}

\begin{Proposition}\label{PrelocSubcatAndQuotOfPrelocSubcat}
	Let $\mcX$ be a localizing subcategory of $\mcA$. Denote the canonical functor by $F\colon\mcA\to\mcA/\mcX$ and its right adjoint by $G\colon\mcA/\mcX\to\mcA$.
	\begin{enumerate}
		\item\label{BijectionBetweenPrelocSubcatAndPrelocSubcatOfQuotSubcat} There is a bijection
		\begin{equation*}
			\setwithspace[\bigg]{
			\begin{gathered}
				\text{prelocalizing subcategories }\mcY\text{ of }\mcA\\
				\text{ satisfying }\mcX*\mcY*\mcX=\mcY
			\end{gathered}
			}
			\to\setwithspace[\bigg]{\text{prelocalizing subcategories of }\frac{\mcA}{\mcX}}
		\end{equation*}
		given by $\mcY\mapsto F(\mcY)$. Its inverse is given by $\mcY'\mapsto F^{-1}(\mcY')$.
		\item\label{ExtOfQuotOfPrelocSubcat} For each $i=1,2$, let $\mcY_{i}$ be a prelocalizing subcategory of $\mcA$ such that $\mcX*\mcY_{i}*\mcX=\mcY_{i}$. Then
		\begin{equation*}
			F(\mcY_{1}*\mcY_{2})=F(\mcY_{1})*F(\mcY_{2}).
		\end{equation*}
		\item\label{BijectionBetweenLocSubcatAndLocSubcatOfQuotSubcat} The bijection in \autoref{BijectionBetweenPrelocSubcatAndPrelocSubcatOfQuotSubcat} restricts to a bijection
		\begin{equation*}
			\setwithspace[\bigg]{
			\begin{gathered}
				\text{localizing subcategories }\mcY\text{ of }\mcA\\
				\text{ satisfying }\mcX\subset\mcY
			\end{gathered}
			}
			\to\setwithspace[\bigg]{\text{localizing subcategories of }\frac{\mcA}{\mcX}}.
		\end{equation*}
	\end{enumerate}
\end{Proposition}

\begin{proof}
	\autoref{BijectionBetweenPrelocSubcatAndPrelocSubcatOfQuotSubcat} By \autoref{QuotOfPrelocSubcat} \autoref{InverseImageOfPrelocSubcatIsPrelocSubcat} and \autoref{QuotOfPrelocSubcat} \autoref{QuotOfPrelocSubcatIsPrelocSubcat}, these maps are well-defined. Let $\eta\colon 1_{\mcA}\to GF$ be the unit morphism.
	
	Let $\mcY$ be a prelocalizing subcategory of $\mcA$ satisfying $\mcX*\mcY*\mcX=\mcY$. It is obvious that $\mcY\subset F^{-1}F(\mcY)$. Let $M$ be an object in $\mcA$ which belongs to $F^{-1}F(\mcY)$. Then there exists an object $N$ in $\mcA$ which belongs to $\mcY$ such that $F(M)\cong F(N)$. There is an exact sequence
	\begin{equation*}
		0\to\Im\eta_{N}\to GF(N)\to\Cok\eta_{N}\to 0.
	\end{equation*}
	The quotient object $\Im\eta_{N}$ of $N$ belongs to $\mcY$. By \autoref{PropertiesOfQuotCat} \autoref{DescriptionOfUnit}, the object $\Cok\eta_{N}$ belongs to $\mcX$. Hence $GF(M)\cong GF(N)$ belongs to $\mcY*\mcX$. By \autoref{PropertiesOfExt} \autoref{ExtOfPrelocSubcatsIsPreloc}, the subobject $\Im\eta_{M}$ of $GF(M)$ belongs to $\mcY*\mcX$. There is an exact sequence
	\begin{equation*}
		0\to\Ker\eta_{M}\to M\to\Im\eta_{M}\to 0,
	\end{equation*}
	where $\Ker\eta_{M}$ belongs to $\mcX$. Therefore $M$ belongs to $\mcX*\mcY*\mcX=\mcY$. This shows that $F^{-1}F(\mcY)\subset\mcY$.
	
	Let $\mcY'$ be a prelocalizing subcategory of $\mcA/\mcX$. It is obvious that $FF^{-1}(\mcY')\subset\mcY'$. Let $M'$ be an object in $\mcA/\mcX$ which belongs to $\mcY'$. Then by \autoref{PropertiesOfQuotCat} \autoref{QuotFunctorIsSurj}, there exists an object $M$ in $\mcA$ such that $F(M)=M'$. Since $M$ belongs to $F^{-1}(\mcY')$, the object $M'=F(M)$ belongs to $FF^{-1}(\mcY')$. This shows that $\mcY'\subset FF^{-1}(\mcY')$.
	
	\autoref{ExtOfQuotOfPrelocSubcat} By \autoref{QuotOfPrelocSubcat} \autoref{QuotOfExtOfSubcatsIsExtOfQuotOfSubcats},
	\begin{equation*}
		F(\mcY_{1}*\mcY_{2})=F(\mcY_{1}*\mcX*\mcY_{2})=F(\mcY_{1})*F(\mcY_{2}).
	\end{equation*}
	
	\autoref{BijectionBetweenLocSubcatAndLocSubcatOfQuotSubcat} This follows from \autoref{ExtOfQuotOfPrelocSubcat}.
\end{proof}

\begin{Remark}\label{ExtOfQuotIsNotNecQuotOfExt}
	In the setting of \autoref{QuotOfPrelocSubcat} \autoref{QuotOfExtOfSubcatsIsExtOfQuotOfSubcats}, the assertion $F(\mcY_{1}*\mcX*\mcY_{2})=F(\mcY_{1}*\mcY_{2})$ does not necessarily hold. The next example gives a counter-example.
\end{Remark}

\begin{Example}\label{ExOfExtOfQuotIsNotNecQuotOfExt}
	Let $K$ be a field, and let $\Lambda$ be the ring
	\begin{equation*}
		\Lambda=
		\begin{bmatrix}
			K & 0 & 0\\
			K & K & 0\\
			K & K & K
		\end{bmatrix}
	\end{equation*}
	of $3\times 3$ lower triangular matrices. Define simple $\Lambda$-modules $S_{i}$ for each $i=1,2,3$ by
	\begin{align*}
		S_{1}&=
		\begin{bmatrix}
			K & 0 & 0
		\end{bmatrix}
		,\\
		S_{2}&=
		\frac{
		\begin{bmatrix}
			K & K & 0
		\end{bmatrix}
		}{
		\begin{bmatrix}
			K & 0 & 0
		\end{bmatrix}
		},\\
		S_{3}&=
		\frac{
		\begin{bmatrix}
			K & K & K
		\end{bmatrix}
		}{
		\begin{bmatrix}
			K & K & 0
		\end{bmatrix}
		},
	\end{align*}
	and let $\mcX_{i}$ be the localizing subcategory of $\Mod\Lambda$ consisting of arbitrary direct sums of copies of $S_{i}$. Let $F\colon\mcA\to\mcA/\mcX_{2}$ and $G\colon\mcA/\mcX_{2}\to\mcA$ denote the canonical functors. Since the $\Lambda$-module
	\begin{equation*}
		M=
		\begin{bmatrix}
			K & K & K
		\end{bmatrix}
	\end{equation*}
	belongs to $\mcX_{1}*\mcX_{2}*\mcX_{3}$, it follows that $M\cong GF(M)$ belongs to $GF(\mcX_{1}*\mcX_{2}*\mcX_{3})$.
	
	On the other hand, every $\Lambda$-module belonging to $\mcX_{1}*\mcX_{3}$ is the direct sum of some object in $\mcX_{1}$ and some object in $\mcX_{3}$. Since $\Mod\Lambda$ is a locally noetherian Grothendieck category, by \cite[Proposition~5.8.12]{Popescu}, the functor $G$ commutes with arbitrary direct sums. Hence every $\Lambda$-module belonging to $GF(\mcX_{1}*\mcX_{3})$ is the direct sum of some object in $GF(\mcX_{1})=\mcX_{1}*\mcX_{2}$ and some object in $GF(\mcX_{3})=\mcX_{3}$. Since $M$ is indecomposable and belongs to neither $\mcX_{1}*\mcX_{2}$ nor $\mcX_{3}$, the $\Lambda$-module $M$ does not belong to $GF(\mcX_{1}*\mcX_{3})$. This shows that $F(\mcX_{1}*\mcX_{2}*\mcX_{3})\not\subset F(\mcX_{1}*\mcX_{3})$.
\end{Example}

The following result gives a characterization of a quotient category.

\begin{Proposition}\label{CharacterizationOfQuotCat}
	Let $\mcA$ and $\mcB$ be Grothendieck categories, and let $Q\colon\mcA\to\mcB$ be an exact functor with a fully faithful right adjoint $\mcB\to\mcA$. Then the full subcategory
	\begin{equation*}
		\mcX=\conditionalset{M\in\mcA}{Q(M)=0}
	\end{equation*}
	of $\mcA$ is a localizing subcategory, and there exists a unique equivalence $\overline{Q}\colon\mcA/\mcX\isoto\mcB$ such that $\overline{Q}F=Q$, where $F\colon\mcA\to\mcA/\mcX$ is the canonical functor.
\end{Proposition}

\begin{proof}
	\cite[Theorem~4.4.9]{Popescu}.
\end{proof}

We state some facts on the image of a localizing subcategory in a quotient category.

\begin{Proposition}\label{QuotOfLocSubcat}
	Let $\mcX$ and $\mcY$ be localizing subcategories of $\mcA$. Denote the canonical functor by $F\colon\mcA\to\mcA/\mcX$ and its right adjoint by $G\colon\mcA/\mcX\to\mcA$.
	\begin{enumerate}
		\item\label{DescriptionOfQuotOfLocSubcat} It holds that $\loc{F(\mcY)}=F(\loc{\mcX\cup\mcY})$.
		\item\label{LocClosureOfImageOfLocSubcat} If $\mcX\subset\mcY$, then the composite $\mcY\to\mcA\to\mcA/\mcX$ induces an equivalence
		\begin{equation*}
			\frac{\mcY}{\mcX}\isoto F(\mcY).
		\end{equation*}
		\item\label{QuotOfQuotCat} If $\mcX\subset\mcY$, then the composite
		\begin{equation*}
			\mcA\to\mcA/\mcX\to\frac{\mcA/\mcX}{F(\mcY)}
		\end{equation*}
		induces an equivalence
		\begin{equation*}
			\frac{\mcA}{\mcY}\isoto\frac{\mcA/\mcX}{F(\mcY)}.
		\end{equation*}
	\end{enumerate}
\end{Proposition}

\begin{proof}
	\autoref{DescriptionOfQuotOfLocSubcat} It is obvious that $\loc{F(\mcY)}\subset F(\loc{\mcX\cup\mcY})$. Since $F$ is exact and commutes with arbitrary direct sums,
	\begin{equation*}
		F(\loc{\mcX\cup\mcY})\subset\loc{F(\mcX\cup\mcY)}=\loc{F(\mcX)\cup F(\mcY)}=\loc{F(\mcY)}
	\end{equation*}
	by \autoref{FundamentalAndUniversalPropertyOfQuotCat} \autoref{FundamentalPropertyOfQuotCat}.
	
	\autoref{LocClosureOfImageOfLocSubcat} The equivalence follows from the construction of $\mcA/\mcX$ (see \cite[p.~365]{Gabriel} or \cite[Definition~5.2]{Kanda3}).
	
	\autoref{QuotOfQuotCat} By \autoref{PrelocSubcatAndQuotOfPrelocSubcat} \autoref{BijectionBetweenLocSubcatAndLocSubcatOfQuotSubcat}, the full subcategory $F(\mcY)$ of $\mcA/\mcX$ is a localizing subcategory, and $F^{-1}F(\mcY)=\mcY$. By \autoref{PropertiesOfQuotCat} \autoref{CounitIsIso}, the composite is an exact functor with a fully faithful right adjoint. Hence by \autoref{CharacterizationOfQuotCat}, it induces an equivalence
	\begin{equation*}
		\frac{\mcA}{\mcY}=\frac{\mcA}{F^{-1}F(\mcY)}\isoto\frac{\mcA/\mcX}{F(\mcY)}.\qedhere
	\end{equation*}
\end{proof}

\section{Atom spectra of quotient categories and localization}
\label{sec:ASpecOfQuotCatsAndLoc}

Throughout this section, let $\mcA$ be a Grothendieck category. We recall a description of the atom spectrum of a quotient category of $\mcA$ and fundamental results on the localization of $\mcA$ at an atom. We start with relating localizing subcategories of $\mcA$ and localizing subsets of $\ASpec\mcA$.

\begin{Definition}\label{ASuppOfSubcatAndASuppInverse}\leavevmode
	\begin{enumerate}
		\item\label{ASuppOfSubcat} For a full subcategory $\mcX$ of $\mcA$, define the subset $\ASupp\mcX$ of $\ASpec\mcA$ by
		\begin{equation*}
			\ASupp\mcX=\bigcup_{M\in\mcX}\ASupp M.
		\end{equation*}
		\item\label{ASuppInverse} For a subset $\Phi$ of $\ASpec\mcA$, define the full subcategory $\ASupp^{-1}\Phi$ of $\mcA$ by
		\begin{equation*}
			\ASupp^{-1}\Phi=\conditionalset{M\in\mcA}{\ASupp M\subset\Phi}.
		\end{equation*}
	\end{enumerate}
\end{Definition}

\begin{Proposition}\label{ASuppOfSubcatIsLocalizingAndASuppInverseIsLocalizing}\leavevmode
	\begin{enumerate}
		\item\label{ASuppOfSubcatIsLocalizing} For every full subcategory $\mcX$ of $\mcA$, the subset $\ASupp\mcX$ of $\ASpec\mcA$ is a localizing subset.
		\item\label{ASuppInverseIsLocalizing} For every subset $\Phi$ of $\ASpec\mcA$, the full subcategory $\ASupp^{-1}\Phi$ of $\mcA$ is a localizing subcategory.
	\end{enumerate}
\end{Proposition}

\begin{proof}
	\autoref{ASuppOfSubcatIsLocalizing} Recall that $\ASpec\mcA$ is in bijection with a small set. For each $\alpha\in\ASupp\mcX$, choose an object $M(\alpha)$ in $\mcA$ which belongs to $\mcX$ such that $\alpha\in\ASupp M(\alpha)$. Then
	\begin{equation*}
		\ASupp\bigoplus_{\alpha\in\ASupp\mcX}M(\alpha)=\bigcup_{\alpha\in\ASupp\mcX}\ASupp M(\alpha)=\ASupp\mcX
	\end{equation*}
	by \autoref{AAssAndASuppAndDirectSumAndSum} \autoref{AAssAndASuppAndDirectSum}.
	
	\autoref{ASuppInverseIsLocalizing} This follows from \autoref{AAssAndASuppAndShortExactSeq} and \autoref{AAssAndASuppAndDirectSumAndSum} \autoref{AAssAndASuppAndDirectSum}.
\end{proof}

The following result shows that a localizing subset of $\ASpec\mcA$ is determined by the corresponding localizing subcategory of $\mcA$.

\begin{Proposition}\label{ASuppOfASuppInverseIsIdentity}
	For every localizing subset $\Phi$ of $\ASpec\mcA$,
	\begin{equation*}
		\ASupp(\ASupp^{-1}\Phi)=\Phi.
	\end{equation*}
\end{Proposition}

\begin{proof}
	This follows from the proof of \cite[Theorem~4.3]{Kanda1}.
\end{proof}

If $\mcA$ is a locally noetherian Grothendieck category, we also have $\ASupp^{-1}(\ASupp\mcX)=\mcX$ for every localizing subcategory $\mcX$ of $\mcA$, and these correspondences establish a bijection between the localizing subcategories of $\mcA$ and the localizing subsets of $\ASpec\mcA$ (\cite[Theorem~5.5]{Kanda1}). We generalize this result later as \autoref{BijectionBetweenLocSubcatsAndLocSubsets}.

We describe the atom spectrum of the quotient category by a localizing subcategory.

\begin{Theorem}\label{ASpecOfQuotCat}
	Let $\mcA$ be a Grothendieck category, and let $\mcX$ be a localizing subcategory of $\mcA$. Denote the canonical functor by $F\colon\mcA\to\mcA/\mcX$ and its right adjoint by $G\colon\mcA/\mcX\to\mcA$. Then the map $\ASpec\mcA\setminus\ASupp\mcX\to\ASpec(\mcA/\mcX)$ given by $\overline{H}\mapsto\overline{F(H)}$ is a homeomorphism. Its inverse is given by $\overline{H'}\mapsto\overline{G(H')}$.
\end{Theorem}

\begin{proof}
	\cite[Theorem~5.17]{Kanda3}.
\end{proof}

\begin{Remark}\label{ASpecOfLocSubcatAndQuotCat}
	Every localizing subcategory $\mcX$ of $\mcA$ is a Grothendieck category, and $\ASpec\mcX$ is homeomorphic to the localizing subset $\ASupp\mcX$ of $\ASpec\mcA$ by the correspondence $\overline{H}\mapsto\overline{H}$ (\cite[Proposition~5.12]{Kanda3}). We identify $\ASpec\mcX$ with $\ASupp\mcX$, and $\ASpec(\mcA/\mcX)$ with $\ASpec\mcA\setminus\ASupp\mcX$ via the homeomorphism in \autoref{ASpecOfQuotCat}. Then
	\begin{equation*}
		\ASpec\mcA=\ASpec\mcX\cup\ASpec\frac{\mcA}{\mcX},
	\end{equation*}
	and
	\begin{equation*}
		\ASpec\mcX\cap\ASpec\frac{\mcA}{\mcX}=\emptyset.
	\end{equation*}
\end{Remark}

We describe atom supports and associated atoms in a quotient category.

\begin{Proposition}\label{AAssAndASuppAndQuotCat}
	Let $\mcX$ be a localizing subcategory of $\mcA$. Denote the canonical functor by $F\colon\mcA\to\mcA/\mcX$ and its right adjoint by $G\colon\mcA/\mcX\to\mcA$.
	\begin{enumerate}
		\item\label{AAssAndASuppAndSectionFunctor} For every object $M'$ in $\mcA/\mcX$,
		\begin{equation*}
			\AAss G(M')=\AAss M',
		\end{equation*}
		and
		\begin{equation*}
			\ASupp G(M')\setminus\ASupp\mcX=\ASupp M'.
		\end{equation*}
		\item\label{AAssAndASuppAndCanonicalFunctor} For every object $M$ in $\mcA$,
		\begin{equation*}
			\AAss F(M)\supset\AAss M\setminus\ASupp\mcX,
		\end{equation*}
		and
		\begin{equation*}
			\ASupp F(M)=\ASupp M\setminus\ASupp\mcX.
		\end{equation*}
	\end{enumerate}
\end{Proposition}

\begin{proof}
	These follow from \cite[Lemma~5.16]{Kanda3}. By considering \autoref{PropertiesOfQuotCat} \autoref{TorsionfreenessOfClosedObj}, the assertion $\AAss G(M')=\AAss M'$ also follows.
\end{proof}

The atom spectrum of the image of a localizing subcategory in a quotient category is described as follows.

\begin{Proposition}\label{ASpecOfQuotOfLocSubcat}
	Let $\mcX$ and $\mcY$ be localizing subcategories of $\mcA$. Denote the canonical functor by $F\colon\mcA\to\mcA/\mcX$ and its right adjoint by $G\colon\mcA/\mcX\to\mcA$. Then
	\begin{align*}
		\ASpec\loc{F(\mcY)}&=\ASpec F(\loc{\mcX\cup\mcY})\\
		&=\ASpec\mcY\cap\ASpec\frac{\mcA}{\mcX}\\
		&=\ASpec\mcY\setminus\ASpec\mcX,
	\end{align*}
	and
	\begin{align*}
		\ASpec\frac{\mcA/\mcX}{\loc{F(\mcY)}}&=\ASpec\frac{\mcA}{\loc{\mcX\cup\mcY}}\\
		&=\ASpec\frac{\mcA}{\mcX}\cap\ASpec\frac{\mcA}{\mcY}\\
		&=\ASpec\mcA\setminus(\ASpec\mcX\cup\ASpec\mcY).
	\end{align*}
\end{Proposition}

\begin{proof}
	This follows from $\ASupp\loc{\mcX\cup\mcY}=\ASupp\mcX\cup\ASupp\mcY$ and \autoref{QuotOfLocSubcat}.
\end{proof}

\begin{Definition}\label{LocalizedCat}
	Let $\mcA$ be a Grothendieck category and $\alpha\in\ASpec\mcA$. Define a localizing subcategory $\mcX(\alpha)$ of $\mcA$ by $\mcX(\alpha)=\ASupp^{-1}(\ASpec\mcA\setminus\Lambda(\alpha))$. Define the \emph{localization} $\mcA_{\alpha}$ of $\mcA$ at $\alpha$ by $\mcA_{\alpha}=\mcA/\mcX(\alpha)$. The canonical functor $\mcA\to\mcA_{\alpha}$ is denoted by $(-)_{\alpha}$.
\end{Definition}

In \autoref{LocalizedCat}, the subset $\ASpec\mcA\setminus\Lambda(\alpha)$ of $\ASpec\mcA$ is localizing. By \autoref{ASuppOfASuppInverseIsIdentity}, $\ASupp\mcX(\alpha)=\ASpec\mcA\setminus\Lambda(\alpha)$. Therefore we have the following result.

\begin{Theorem}\label{ASpecOfLocalizedCat}
	Let $\mcA$ be a Grothendieck category and $\alpha\in\ASpec\mcA$. Then $\ASpec\mcA_{\alpha}=\Lambda(\alpha)$. In particular, the partially ordered set $\ASpec\mcA$ has the largest element $\alpha$.
\end{Theorem}

\begin{proof}
	\cite[Proposition~6.6 (1)]{Kanda3}.
\end{proof}

We obtain the following description of atom supports.

\begin{Proposition}\label{ASuppAndLoc}\leavevmode
	\begin{enumerate}
		\item\label{DescriptionPrimeLocSubcatByASupp} For every $\alpha\in\ASpec\mcA$,
		\begin{equation*}
			\mcX(\alpha)=\conditionalset{M\in\mcA}{\alpha\notin\ASupp M}.
		\end{equation*}
		\item\label{DescriptionOfASuppByLoc} For every object $M$ in $\mcA$,
	\begin{equation*}
		\ASupp M=\conditionalset{\alpha\in\ASpec\mcA}{M_{\alpha}\neq 0}.
	\end{equation*}
	\end{enumerate}
\end{Proposition}

\begin{proof}
	\cite[Proposition~6.2]{Kanda3}.
\end{proof}

We show that the localization of a Grothendieck category at an atom is ``local'' in the following sense.

\begin{Definition}\label{LocalGrothCatAndPrimeLocSubcat}
	Let $\mcA$ be a Grothendieck category.
	\begin{enumerate}
		\item\label{LocalGrothCat} We say that $\mcA$ is \emph{local} if there exists a simple object in $\mcA$ such that $E(S)$ is a cogenerator of $\mcA$.
		\item\label{PrimeLocSubcat} A localizing subcategory $\mcX$ of $\mcA$ is called \emph{prime} if $\mcA/\mcX$ is a local Grothendieck category.
	\end{enumerate}
\end{Definition}

\begin{Theorem}\label{CharacterizationOfLocalGrothCat}
	Let $\mcA$ be a Grothendieck category. Then the following assertions are equivalent.
	\begin{enumerate}
		\item $\mcA$ is local.
		\item There exists $\alpha\in\ASpec\mcA$ such that for every nonzero object $M$ in $\mcA$, we have $\alpha\in\ASupp M$.
		\item There exists $\alpha\in\ASpec\mcA$ such that the canonical functor $\mcA\to\mcA_{\alpha}$ is an equivalence.
	\end{enumerate}
\end{Theorem}

\begin{proof}
	\cite[Proposition~6.4 (1)]{Kanda3} and \cite[Proposition~6.6 (2)]{Kanda3}.
\end{proof}

In the case of where $\mcA$ is a locally noetherian Grothendieck category, the localness of $\mcA$ is characterized as follows.

\begin{Proposition}\label{CharacterizationOfLocalLocNoethGrothCat}
	Let $\mcA$ be a Grothendieck category. If $\mcA$ is local, then all simple objects in $\mcA$ are isomorphic to each other. In the case where $\mcA$ is a nonzero locally noetherian Grothendieck category, the converse also holds.
\end{Proposition}

\begin{proof}
	\cite[Proposition~6.4 (2)]{Kanda3}.
\end{proof}

\autoref{CharacterizationOfLocalGrothCat} shows that the localizing subcategory $\mcX(\alpha)$ is prime for every $\alpha\in\ASpec\mcA$. This correspondence gives the following bijection.

\begin{Theorem}\label{BijectionBetweenAtomsAndPrimeLocSubcat}
	Let $\mcA$ be a Grothendieck category. There is a bijection
	\begin{equation*}
		\ASpec\mcA\to\setwithspace{\text{prime localizing subcategories of }\mcA}
	\end{equation*}
	given by $\alpha\mapsto\mcX(\alpha)$. For each $\alpha,\beta\in\ASpec\mcA$, we have $\alpha\leq\beta$ if and only if $\mcX(\alpha)\supset\mcX(\beta)$.
\end{Theorem}

\begin{proof}
	\cite[Theorem~6.8]{Kanda3}.
\end{proof}

We consider the localization of a quotient category.

\begin{Proposition}\label{LocOfQuotCat}
	Let $\mcX$ be a localizing subcategory of $\mcA$ and $\alpha\in\ASpec\mcA\setminus\ASupp\mcX$. Then the composite of the canonical functors $\mcA\to\mcA/\mcX$ and $\mcA/\mcX\to(\mcA/\mcX)_{\alpha}$ induces an equivalence $\mcA_{\alpha}\isoto(\mcA/\mcX)_{\alpha}$.
\end{Proposition}

\begin{proof}
	By \autoref{ASuppAndLoc} \autoref{DescriptionPrimeLocSubcatByASupp}, $\mcX\subset\mcX(\alpha)$. Hence the claim follows from \autoref{AAssAndASuppAndQuotCat} \autoref{AAssAndASuppAndCanonicalFunctor} and \autoref{QuotOfLocSubcat} \autoref{QuotOfQuotCat}.
\end{proof}

In the setting of \autoref{LocOfQuotCat}, we identify $\mcA_{\alpha}$ and $(\mcA/\mcX)_{\alpha}$.

The following result shows that the localization of a Grothendieck category at an atom is a generalization of the localization a commutative ring at a prime ideal.

\begin{Proposition}\label{LocOfCommRingAndLocCommRing}
	Let $R$ be a commutative ring.
	\begin{enumerate}
		\item\label{LocOfCommRing} Let $\mfp\in\Spec R$. Denote by $\alpha$ the corresponding atom $\overline{R/\mfp}$ in $\Mod R$. Then the functor $-\otimes_{R}R_{\mfp}\colon\Mod R\to\Mod R_{\mfp}$ induces an equivalence $(\Mod R)_{\alpha}\isoto\Mod R_{\mfp}$.
		\item\label{LocalnessOfModCatOverCommRing} The Grothendieck category $\Mod R$ is local if and only if the commutative ring $R$ is local.
	\end{enumerate}
\end{Proposition}

\begin{proof}
	\autoref{LocOfCommRing} \cite[Proposition~6.9]{Kanda3}.
	
	\autoref{LocalnessOfModCatOverCommRing} This follows from \autoref{CharacterizationOfLocalGrothCat} and \autoref{LocOfCommRing}.
\end{proof}

\section{Grothendieck categories with enough atoms}
\label{sec:GrothCatsWithEnoughAtoms}

The purpose of this paper is to investigate the category $\QCoh X$ of quasi-coherent sheaves on a locally noetherian scheme $X$. In general, the category $\QCoh X$ is a Grothendieck category but not necessarily locally noetherian (see \autoref{CohSheafIsNotNecFinGenObj}). In this section, we introduce the notion of a Grothendieck category with enough atoms and investigate its properties. It is shown later that $\QCoh X$ is a Grothendieck category with enough atoms.

Let $\mcA$ be a Grothendieck category. Recall that every monoform object in $\mcA$ is uniform (\autoref{PropertiesOfMonoformObj} \autoref{MonoformObjIsUniform}). We say that uniform objects $U_{1}$ and $U_{2}$ in $\mcA$ are \emph{equivalent} (denoted by $U_{1}\sim U_{2}$) if there exists a nonzero subobject of $U_{1}$ which is isomorphic to a subobject of $U_{2}$. The equivalence between monoform objects is exactly the same as the atom-equivalence defined in \autoref{MonoformObjAndAtomEquiv} \autoref{AtomEquiv}.

\begin{Proposition}\label{UniformEquivAndIndecInjIso}
	Let $U_{1}$ and $U_{2}$ be uniform objects in $\mcA$. Then $U_{1}$ is equivalent to $U_{2}$ if and only if $E(U_{1})$ is isomorphic to $E(U_{2})$.
\end{Proposition}

\begin{proof}
	\cite[Lemma~2]{Krause2}.
\end{proof}

Since every indecomposable injective object in $\mcA$ is uniform (\cite[Proposition~V.2.8]{Stenstrom}), the map
\begin{equation*}
	\frac{\setwithspace{\text{uniform objects in }\mcA}}{\sim}\to\frac{\setwithspace{\text{indecomposable injective objects in }\mcA}}{\cong}
\end{equation*}
induced by the correspondence $U\mapsto E(U)$ is bijective. We consider the restriction of this bijection to $\ASpec\mcA$.

\begin{Definition}\label{InjEnvOfAtom}
	Let $\mcA$ be a Grothendieck category. For $\alpha\in\ASpec\mcA$, define the \emph{injective envelope} $E(\alpha)$ of $\alpha$ by $E(\alpha)=E(H)$, where $H$ is a monoform object in $\mcA$ satisfying $\overline{H}=\alpha$.
\end{Definition}

\autoref{UniformEquivAndIndecInjIso} implies that the isomorphism class of $E(\alpha)$ in \autoref{InjEnvOfAtom} does not depend on the choice of the representative $H$.

\begin{Definition}\label{GrothCatWithEnoughAtoms}
	We say that a Grothendieck category $\mcA$ has \emph{enough atoms} if $\mcA$ satisfies the following conditions.
	\begin{enumerate}
		\item\label{ExistenceOfIndecDecomp} Every injective object in $\mcA$ has an indecomposable decomposition.
		\item\label{HavingEnoughAtoms} Each indecomposable injective object in $\mcA$ is isomorphic to $E(\alpha)$ for some $\alpha\in\ASpec\mcA$.
	\end{enumerate}
\end{Definition}

Note that an indecomposable decomposition of an injective object is unique in the following sense.

\begin{Theorem}\label{UniquenessOfIndecDecompIntoIndecInjObjs}
	Let $\mcA$ be a Grothendieck category, and let $I$ be an injective object with
		\begin{equation*}
			I\cong\bigoplus_{\lambda\in\Lambda}I_{\lambda}\cong\bigoplus_{\mu\in\Lambda'}I'_{\mu},
		\end{equation*}
		where $I_{\lambda}$ and $I'_{\mu}$ are indecomposable for each $\lambda\in\Lambda$ and $\mu\in\Lambda'$. Then there exists a bijection $\varphi\colon\Lambda\to\Lambda'$ such that $I_{\lambda}$ is isomorphic to $I'_{\varphi(\lambda)}$ for each $\lambda\in\Lambda$.
\end{Theorem}

\begin{proof}
	This follows from Krull--Remak--Schmidt--Azumaya's theorem (\cite[Theorem~5.1.3]{Popescu}) and the fact that the endomorphism ring of each indecomposable injective object in $\mcA$ is local (\cite[Lemma~4.20.3]{Popescu}).
\end{proof}

The following result shows that a Grothendieck category with enough atoms is a generalization of a locally noetherian Grothendieck category.

\begin{Proposition}\label{LocNoethGrothCatHasEnoughAtoms}
	Every locally noetherian Grothendieck category has enough atoms.
\end{Proposition}

\begin{proof}
	This follows from \autoref{PropertiesOfMonoformObj} \autoref{NoethObjHasMonoformSub} and \cite[Proposition~V.4.5]{Stenstrom} since every nonzero object in a locally noetherian Grothendieck category has a nonzero noetherian subobject.
\end{proof}

We show that every quotient category of a Grothendieck category with enough atoms has enough atoms.

\begin{Proposition}\label{QuotCatOfGrothCatWithEnoughAtoms}
	Let $\mcA$ be a Grothendieck category, and let $\mcX$ be a localizing subcategory of $\mcA$.
	\begin{enumerate}
		\item\label{IndecDecompInQuotCat} If every injective object in $\mcA$ has an indecomposable decomposition, then every injective object in $\mcA/\mcX$ has an indecomposable decomposition.
		\item\label{QuotCatOfGrothCatWithEnoughAtomsHasEnoughAtoms} If $\mcA$ has enough atoms, then $\mcA/\mcX$ has enough atoms.
	\end{enumerate}
\end{Proposition}

\begin{proof}
	Denote the canonical functor by $F\colon\mcA\to\mcA/\mcX$ and its right adjoint by $G\colon\mcA/\mcX\to\mcA$.
	
	\autoref{IndecDecompInQuotCat} Let $I'$ be an injective object in $\mcA/\mcX$. By \autoref{SubobjAndInjectivityInQuotCat} \autoref{ImageOfInjObjInQuotCat}, the object $G(I')$ in $\mcA$ is injective. Hence $G(I')$ has an indecomposable decomposition
	\begin{equation*}
		G(I')=\bigoplus_{\lambda\in\Lambda}I_{\lambda}.
	\end{equation*}
	We obtain
	\begin{equation*}
		I'\cong FG(I')\cong\bigoplus_{\lambda\in\Lambda}F(I_{\lambda}).
	\end{equation*}
	By \autoref{PropertiesOfQuotCat} \autoref{TorsionfreenessOfClosedObj}, \autoref{PropertiesOfImageOfObj} \autoref{UniformnessOfImageOfObj} and \autoref{PropertiesOfImageOfObj} \autoref{InjectivityOfImageOfObj}, the object $F(I_{\lambda})$ is an indecomposable injective object in $\mcA/\mcX$ for each $\lambda\in\Lambda$.
	
	\autoref{QuotCatOfGrothCatWithEnoughAtomsHasEnoughAtoms} Let $I'$ be an indecomposable injective object in $\mcA/\mcX$. Then by \autoref{SubobjAndInjectivityInQuotCat} \autoref{ImageOfIndecObjInQuotCat} and \autoref{SubobjAndInjectivityInQuotCat} \autoref{ImageOfInjObjInQuotCat}, the object $G(I')$ in $\mcA$ is indecomposable and injective. Hence there exists $\alpha\in\ASpec\mcA$ such that $G(I')\cong E_{\mcA}(\alpha)$. We obtain $I'\cong FG(I')\cong F(E_{\mcA}(\alpha))$. Let $H$ be a monoform subobject of $E_{\mcA}(\alpha)$. By \autoref{PropertiesOfQuotCat} \autoref{TorsionfreenessOfClosedObj}, the object $H$ is $\mcX$-torsionfree. By \autoref{PropertiesOfImageOfObj} \autoref{MonoformnessOfImageOfObj}, the object $I'$ has the monoform subobject $F(H)$. This implies that $I'=E(F(H))=E_{\mcA/\mcX}(\alpha)$.
\end{proof}

A Grothendieck category $\mcA$ is called \emph{locally uniform}\footnote{In \cite[p.~330]{Popescu}, it is called \emph{locally coirreducible} since a uniform object is called a \emph{coirreducible object}.} if every nonzero object in $\mcA$ has a uniform subobject. It is shown that this holds whenever $\mcA$ has enough atoms.

\begin{Proposition}\label{GrothCatWithEnoughAtomsIsLocMonoform}
	Let $\mcA$ be a Grothendieck category with enough atoms. Then every nonzero object in $\mcA$ has a monoform subobject. In particular, the Grothendieck category $\mcA$ is locally uniform.
\end{Proposition}

\begin{proof}
	Let $M$ be a nonzero object in $\mcA$. Then there exists a family $\set{\alpha_{\lambda}}_{\lambda\in\Lambda}$ of atoms in $\mcA$ such that
	\begin{equation*}
		E(M)\cong\bigoplus_{\lambda\in\Lambda}E(\alpha_{\lambda}).
	\end{equation*}
	Hence $E(M)$ has a monoform subobject $H$. Since $M$ is an essential subobject of $E(M)$, the subobject $H\cap M$ of $M$ is monoform by \autoref{PropertiesOfMonoformObj} \autoref{SubOfMonoformObjIsMonoform}. The last assertion follows from \autoref{PropertiesOfMonoformObj} \autoref{MonoformObjIsUniform}.
\end{proof}

The classification of the localizing subcategories by the atom spectrum we mentioned after \autoref{ASuppOfASuppInverseIsIdentity} is generalized to a Grothendieck category with enough atoms.

\begin{Theorem}\label{BijectionBetweenLocSubcatsAndLocSubsets}
	Let $\mcA$ be a Grothendieck category with enough atoms. There is a bijection
	\begin{equation*}
		\setwithspace{\text{localizing subcategories of }\mcA}\to\setwithspace{\text{localizing subsets of }\ASpec\mcA}
	\end{equation*}
	given by $\mcX\mapsto\ASupp\mcX$. Its inverse is given by $\Phi\mapsto\ASupp^{-1}\Phi$.
\end{Theorem}

\begin{proof}
	By \autoref{ASuppOfSubcatIsLocalizingAndASuppInverseIsLocalizing} and \autoref{ASuppOfASuppInverseIsIdentity}, it suffices to show that $\ASupp^{-1}(\ASupp\mcX)=\mcX$ for each localizing subcategory $\mcX$ of $\mcA$. The inclusion $\mcX\subset\ASupp^{-1}(\ASupp\mcX)$ holds obviously. Let $M$ be an object in $\mcA$ which belongs to $\ASupp^{-1}(\ASupp\mcX)$, and let $L$ be the largest subobject of $M$ which belongs to $\mcX$. If $M/L$ is nonzero, then by \autoref{GrothCatWithEnoughAtomsIsLocMonoform}, there exists a monoform subobject $H$ of $M/L$. Since $\overline{H}\in\ASupp M\subset\ASupp\mcX$, there exists a nonzero subobject $H'$ of $H$ which belongs to $\mcX$. Let $H'=L'/L\subset M/L$. Since $L$ and $L'/L$ belongs to $\mcX$, the subobject $L'$ of $M$ also belongs to $\mcX$. This contradicts the maximality of $L$. Therefore $\ASupp^{-1}(\ASupp\mcX)=\mcX$.
\end{proof}

We show that every localizing subcategory is the intersection of some family of prime localizing subcategories.

\begin{Corollary}\label{LocSubcatIsIntersectionOfPrimeLocSubcat}
	Let $\mcA$ be a Grothendieck category with enough atoms. For every localizing subcategory $\mcX$ of $\mcA$,
	\begin{equation*}
		\mcX=\bigcap_{\alpha\in\ASpec\mcA\setminus\ASupp\mcX}\mcX(\alpha).
	\end{equation*}
\end{Corollary}

\begin{proof}
	By \autoref{ASuppAndLoc} \autoref{DescriptionPrimeLocSubcatByASupp} and \autoref{BijectionBetweenLocSubcatsAndLocSubsets},
	\begin{align*}
		\bigcap_{\alpha\in\ASpec\mcA\setminus\ASupp\mcX}\mcX(\alpha)&=\conditionalset{M\in\mcA}{\alpha\notin\ASupp M\text{ for each }\alpha\in\ASpec\mcA\setminus\ASupp\mcX}\\
		&=\conditionalset{M\in\mcA}{\ASupp M\subset\ASupp\mcX}\\
		&=\ASupp^{-1}(\ASupp\mcX)\\
		&=\mcX.\qedhere
	\end{align*}
\end{proof}

Let $\mcA$ be a Grothendieck category and $\alpha\in\ASpec\mcA$. It is shown in the proof of \cite[Theorem~2.5]{Kanda2} that the injective envelope $E(\alpha)$ has a largest monoform subobject $H(\alpha)$. The object $H(\alpha)$ is called the \emph{atomic object} corresponding to $\alpha$. It is straightforward to show that no monoform object in $\mcA$ has a proper essential subobject isomorphic to $H(\alpha)$.

The atomic objects correspond to the simple objects in the localizations.

\begin{Proposition}\label{SimpleObjAndAtomicObj}
	Let $\mcA$ be a Grothendieck category and $\alpha\in\ASpec\mcA$. Denote the canonical functor by $F_{\alpha}\colon\mcA\to\mcA_{\alpha}$ and its right adjoint by $G_{\alpha}\colon\mcA_{\alpha}\to\mcA$. Let $S'$ be the simple object in $\mcA_{\alpha}$.
	\begin{enumerate}
		\item\label{SimpleObjIsAtomicInLocalizedCat} $S'$ is the atomic object corresponding to the atom $\overline{S'}$ in $\mcA_{\alpha}$.
		\item\label{AtomicObjIsSimpleObjInLocalizedCat} $G_{\alpha}(S')$ is isomorphic to the atomic object $H(\alpha)$.
		\item\label{EndOfAtomicObjIsSkewField} The ring $\End_{\mcA}(H(\alpha))$ is isomorphic to the skew field $\End_{\mcA_{\alpha}}(S')$.
	\end{enumerate}
\end{Proposition}

\begin{proof}
	\autoref{SimpleObjIsAtomicInLocalizedCat} It holds that $S'\subset H(\overline{S'})\subset E(\overline{S'})=E(S')$. If $S'\subsetneq H(\overline{S'})$, then by \autoref{CharacterizationOfLocalGrothCat}, $\overline{S'}\in\ASupp(H(\overline{S'})/S')$, and hence there exist a subobject $L$ of $H(\overline{S'})$ with $S'\subset L$ and a subobject of $H(\overline{S'})/L$ which is isomorphic to $S'$. This contradicts the monoformness of $H(\overline{S'})$. Therefore $S'=H(\overline{S'})$.
	
	\autoref{AtomicObjIsSimpleObjInLocalizedCat} By \autoref{ASpecOfQuotCat}, the object $F_{\alpha}(H(\alpha))$ is a monoform object in $\mcA_{\alpha}$, and $G_{\alpha}F_{\alpha}(H(\alpha))$ is a monoform object in $\mcA$. By \autoref{SimpleObjIsAtomicInLocalizedCat}, $F_{\alpha}(H(\alpha))\cong S'$. Since $H(\alpha)$ is $\mcX(\alpha)$-torsionfree, by \autoref{PropertiesOfQuotCat} \autoref{DescriptionOfUnit}, the canonical morphism $H(\alpha)\to G_{\alpha}F_{\alpha}(H(\alpha))$ is a monomorphism, and $H(\alpha)$ is essential as a subobject of $G_{\alpha}F_{\alpha}(H(\alpha))$. Therefore the morphism $H(\alpha)\to G_{\alpha}F_{\alpha}(H(\alpha))$ is an isomorphism, and $G_{\alpha}(S')\cong G_{\alpha}F_{\alpha}(H(\alpha))\cong H(\alpha)$.
	
	\autoref{EndOfAtomicObjIsSkewField} By \autoref{AtomicObjIsSimpleObjInLocalizedCat} and \autoref{PropertiesOfQuotCat} \autoref{CounitIsIso},
	\begin{align*}
		\End_{\mcA}(H(\alpha))&\cong\End_{\mcA}(G_{\alpha}(S'))\\
		&\cong\Hom_{\mcA_{\alpha}}(F_{\alpha}G_{\alpha}(S'),S')\\
		&\cong\End_{\mcA_{\alpha}}(S').
	\end{align*}
	This gives a ring isomorphism $\End_{\mcA}(H(\alpha))\cong\End_{\mcA_{\alpha}}(S')$.
\end{proof}

The skew field $\End_{\mcA}(H(\alpha))$ is called the \emph{residue field} of $\alpha$ and denoted by $k(\alpha)$.

\section{The atom spectra of locally noetherian schemes}
\label{sec:ASpecOfLocNoethSch}

In this section, we describe the atom spectrum of the category of quasi-coherent sheaves on a locally noetherian scheme. Let $X$ be a locally noetherian scheme with the underlying topological space $\uspX$ and the structure sheaf $\OX$. It is known that the category $\Mod X$ of $\OX$-modules and the category $\QCoh X$ of quasi-coherent sheaves on $X$ are Grothendieck categories (see \cite[Theorem~II.7.8]{Hartshorne1} and \cite[Lemma~2.1.7]{Conrad}). For a commutative ring $R$, we identify $\QCoh(\Spec R)$ with $\Mod R$.

\begin{Proposition}\label{FunctorsBetweenQCoh}
	Let $U$ be an open affine subscheme of $X$, and let $i\colon U\into X$ be the immersion. Then the functor $i_{*}\colon\Mod U\to\Mod X$ and its left adjoint $i^{*}\colon\Mod X\to\Mod U$ induce the functor $i_{*}\colon\QCoh U\to\QCoh X$ and its left adjoint $i^{*}\colon\QCoh X\to\QCoh U$.
\end{Proposition}

\begin{proof}
	\cite[0.4.4.3.1]{Grothendieck1} and \cite[Proposition~I.9.4.2 (i)]{Grothendieck1}.
\end{proof}

In the rest of this paper, every quasi-coherent sheaf $M$ on $X$ is always regarded as an object in $\QCoh X$, not in $\Mod X$. Hence a subobject of $M$ means a quasi-coherent subsheaf of $M$.

For an open affine subscheme $U$ of $X$ with the immersion $i\colon U\into X$, the functor $i^{*}\colon\QCoh X\to\QCoh U$ is also denoted by $(-)|_{U}$. The category $\QCoh U$ is realized as a quotient category of $\QCoh X$ through this functor.

\begin{Proposition}\label{OpenSubschIsQuotCat}
	Let $U$ be an open affine subscheme of $X$. Then the functor $(-)|_{U}\colon\QCoh X\to\QCoh U$ induces an equivalence $(\QCoh X)/\mcX_{U}\isoto\QCoh U$, where $\mcX_{U}$ is a localizing subcategory of $\QCoh X$ defined by
	\begin{equation*}
		\mcX_{U}=\conditionalset{M\in\QCoh X}{M|_{U}=0}.
	\end{equation*}
\end{Proposition}

\begin{proof}
	Let $i\colon U\into X$ be the immersion. Since the counit functor $i^{*}i_{*}\to 1_{\QCoh U}$ is an isomorphism, the functor $i_{*}$ is fully faithful. The functor $i^{*}$ is exact. Hence the claim follows from \autoref{CharacterizationOfQuotCat}.
\end{proof}

For each object $M$ in $\QCoh X$, the subset $\Supp M$ of $X$ is defined by
\begin{equation*}
	\Supp M=\conditionalset{x\in X}{M_{x}\neq 0}.
\end{equation*}
For each $x\in X$, let $j_{x}\colon\Spec\OXx\to X$ be the canonical morphism. Note that $j_{x}^{*}$ is equal to the localization $(-)_{x}\colon\QCoh X\to\Mod\OXx$. The category $\Mod\OXx$ is realized as a quotient category of $\QCoh X$ through this morphism.

\begin{Proposition}\label{LocOfSchIsQuotCat}
	For every $x\in X$, the full subcategory
	\begin{equation*}
		\mcX(x):=\conditionalset{M\in\QCoh X}{x\notin\Supp M}=\conditionalset{M\in\QCoh X}{M_{x}=0}
	\end{equation*}
	of $\QCoh X$ is a prime localizing subcategory. The functor $(-)_{x}\colon\QCoh X\to\Mod\OXx$ induces an equivalence $(\QCoh X)/\mcX(x)\isoto\Mod\OXx$.
\end{Proposition}

\begin{proof}
	Let $i\colon U\into X$ be the immersion of an open affine subscheme with $x\in U$. Then the functor $(-)_{x}\colon\QCoh X\to\Mod\OXx$ is equal to the composite of $(-)|_{U}\colon\QCoh X\to\QCoh U$ and $(-)_{x}\colon\QCoh U\to\Mod\OXx$. By \autoref{OpenSubschIsQuotCat} and \autoref{LocOfCommRingAndLocCommRing} \autoref{LocOfCommRing}, these two functors are exact functors with fully faithful right adjoints. Hence we obtain the equivalence by \autoref{CharacterizationOfQuotCat}. By \autoref{LocOfCommRingAndLocCommRing} \autoref{LocalnessOfModCatOverCommRing}, the localizing subcategory $\mcX(x)$ is prime.
\end{proof}

For each $x\in X$, denote the unique maximal ideal of $\OXx$ by $\mfm_{x}$, the residue field of $x$ by $k(x)=\OXx/\mfm_{x}$, and an injective envelope of $k(x)$ in $\Mod\OXx$ by $E(x)=E_{\OXx}(k(x))$. We state that every injective object in $\QCoh X$ is a direct sum of indecomposable injective objects of this form.

\begin{Theorem}[{Hartshorne \cite{Hartshorne1}}]\label{InjObjInQCoh}
	Let $X=(\uspX,\OX)$ be a locally noetherian scheme.
	\begin{enumerate}
		\item\label{DirectSumOfInjObjsInQCohIsInj} For every family $\set{I_{\lambda}}_{\lambda\in\Lambda}$ of injective objects in $\QCoh X$, the direct sum $\bigoplus_{\lambda\in\Lambda}I_{\lambda}$ is also injective.
		\item\label{InjObjInQCohHasIndecDecomp} Every injective object in $\QCoh X$ has an indecomposable decomposition.
		\item\label{DescriptionOfIndecInjObj} There is a bijection
		\begin{equation*}
			\uspX\to\frac{\setwithspace{\text{indecomposable injective objects in }\QCoh X}}{\cong}
		\end{equation*}
		given by $x\mapsto {j_{x}}_{*}E(x)$.
	\end{enumerate}
\end{Theorem}

\begin{proof}
	\cite[Lemma~2.1.5]{Conrad}.
\end{proof}

\begin{Remark}\label{CohSheafIsNotNecFinGenObj}
	In \cite[p.~135]{Hartshorne1}, it is shown that there exists a locally noetherian scheme $X$ such that the Grothendieck category $\QCoh X$ is not locally noetherian. By combining \autoref{InjObjInQCoh} \autoref{DirectSumOfInjObjsInQCohIsInj} and \cite[Theorem~5.8.7]{Popescu}, we deduce that $\QCoh X$ is not even (categorically) locally finitely generated. On the other hand, the set of coherent sheaves on $X$ generates $\QCoh X$ \cite[Corollary~I.9.4.9]{Grothendieck1}. Consequently, a coherent sheaf on $X$ is not necessarily a finitely generated object in $\QCoh X$.
\end{Remark}

We give a description of the atom spectrum of $\QCoh X$.

\begin{Theorem}\label{DescriptionOfASpecOfLocNoethSch}
	Let $X=(\uspX,\OX)$ be a locally noetherian scheme.
	\begin{enumerate}
		\item\label{AtomCorrespondingToPointInLocNoethSch} For each $x\in X$, the set $\AAss{j_{x}}_{*}E(x)$ consists of one element, say $\alpha_{x}$. The injective envelope of $\alpha_{x}$ is $E(\alpha_{x})={j_{x}}_{*}E(x)$. The atomic object is $H(\alpha_{x})={j_{x}}_{*}k(x)$. The residue field is $k(\alpha_{x})\cong k(x)$.
		\item\label{ASpecOfLocNoethSch} There is a bijection $\uspX\to\ASpec(\QCoh X)$ given by $x\mapsto\alpha_{x}$. Moreover, the Grothendieck category $\QCoh X$ has enough atoms.
	\end{enumerate}
\end{Theorem}

\begin{proof}
	\autoref{AtomCorrespondingToPointInLocNoethSch} By \autoref{LocOfSchIsQuotCat} and \autoref{AAssAndASuppAndQuotCat} \autoref{AAssAndASuppAndSectionFunctor},
	\begin{equation*}
		\AAss{j_{x}}_{*}E(x)=\AAss E(x)=\set{\overline{k(x)}}.
	\end{equation*}
	Since ${j_{x}}_{*}E(x)$ is an indecomposable injective object by \autoref{InjObjInQCoh} \autoref{DescriptionOfIndecInjObj}, it is an injective envelope of each of its nonzero subobjects. Hence $E(\alpha_{x})={j_{x}}_{*}E(x)$.
	
	By \autoref{SimpleObjAndAtomicObj} \autoref{AtomicObjIsSimpleObjInLocalizedCat}, $H(\alpha_{x})={j_{x}}_{*}k(x)$. By \autoref{SimpleObjAndAtomicObj} \autoref{EndOfAtomicObjIsSkewField}, $k(\alpha_{x})\cong\End_{\OXx}(k(x))\cong k(x)$.
	
	\autoref{ASpecOfLocNoethSch} The bijection in \autoref{InjObjInQCoh} \autoref{DescriptionOfIndecInjObj} is the composite of the map
	\begin{equation*}
		\uspX\to\ASpec(\QCoh X)
	\end{equation*}
	given by $x\mapsto\alpha_{x}$ and the injection
	\begin{equation*}
		\ASpec(\QCoh X)\to\frac{\setwithspace{\text{indecomposable injective objects in }\QCoh X}}{\cong}
	\end{equation*}
	given by $\alpha\mapsto E(\alpha)$. Hence these maps are also bijective. By \autoref{InjObjInQCoh} \autoref{InjObjInQCohHasIndecDecomp}, the Grothendieck category $\QCoh X$ has enough atoms.
\end{proof}

A subset $\Phi$ of $X$ is said to be \emph{closed under specialization} if for every $x\in\Phi$, we have $\overline{\set{x}}\subset\Phi$. Atom supports and related notions in $\QCoh X$ are described as follows.

\begin{Corollary}\label{DescriptionOfASuppForLocNoethSch}\leavevmode
	\begin{enumerate}
		\item\label{ASuppIsSuppForLocNoethSch} Let $M$ be an object in $\QCoh X$. Then the bijection $\uspX\to\ASpec(\QCoh X)$ in \autoref{DescriptionOfASpecOfLocNoethSch} \autoref{ASpecOfLocNoethSch} restricts to a bijection $\Supp M\to\ASupp M$.
		\item\label{LocOfLocNoethSch} For each $x\in X$, we have $\mcX(\alpha_{x})=\mcX(x)$. The canonical functor $\QCoh X\to\Mod\OXx$ induces an equivalence $(\QCoh X)_{\alpha_{x}}\isoto\Mod\OXx$.
		\item\label{LocSubsetForLocNoethSch} For each subset $\Phi$ of $X$, the corresponding subset
		\begin{equation*}
			\conditionalset{\alpha_{x}\in\ASpec(\QCoh X)}{x\in\Phi}
		\end{equation*}
		of $\ASpec(\QCoh X)$ is localizing if and only if $\Phi$ is closed under specialization.
		\item\label{SpecializationOrderOfLocNoethSch} Let $x,y\in X$. Then $\alpha_{x}\leq\alpha_{y}$ if and only if $y\in\overline{\set{x}}$.
	\end{enumerate}
\end{Corollary}

\begin{proof}
	\autoref{ASuppIsSuppForLocNoethSch} For each $x\in X$, by \autoref{LocOfSchIsQuotCat} and \autoref{AAssAndASuppAndQuotCat} \autoref{AAssAndASuppAndCanonicalFunctor}, $\alpha_{x}\in\ASupp M$ if and only if $\overline{k(x)}\in\ASupp j_{x}^{*}M$. By \autoref{AAssAndASuppOfCommRing}, this is equivalent to $\mfm_{x}\in\Supp j_{x}^{*}M$, which means $M_{x}=j_{x}^{*}M\neq 0$.
	
	\autoref{LocOfLocNoethSch} By \autoref{ASuppIsSuppForLocNoethSch} and \autoref{ASuppAndLoc} \autoref{DescriptionPrimeLocSubcatByASupp}, $\mcX(\alpha_{x})=\mcX(x)$. The equivalence follows from \autoref{LocOfSchIsQuotCat}.
	
	\autoref{LocSubsetForLocNoethSch} By \autoref{LocOfLocNoethSch}, it suffices to show that $\Phi$ is closed under specialization if and only if there exists an object $M$ in $\QCoh X$ satisfying $\Phi=\Supp M$. For every object $M$ in $\QCoh X$, it is straightforward to show that $\Supp M$ is closed under specialization.
	
	Assume that $\Phi$ is closed under specialization. For each $x\in\Phi$, we have $\Supp{j_{x}}_{*}k(x)=\overline{\set{x}}$. Hence
	\begin{equation*}
		\Supp\bigoplus_{x\in\Phi}{j_{x}}_{*}k(x)=\bigcup_{x\in\Phi}\Supp{j_{x}}_{*}k(x)=\Phi.
	\end{equation*}
	
	\autoref{SpecializationOrderOfLocNoethSch} This follows from \autoref{LocSubsetForLocNoethSch}.
\end{proof}

We specialize \autoref{BijectionBetweenLocSubcatsAndLocSubsets} to the case of $\QCoh X$. For a full subcategory $\mcX$ of $\QCoh X$, define the specialization-closed subset $\Supp\mcX$ of $X$ by
\begin{equation*}
	\Supp\mcX=\bigcup_{M\in\mcX}\Supp M.
\end{equation*}
For a subset $\Phi$ of $X$, define the localizing subcategory $\Supp^{-1}\Phi$ of $\QCoh X$ by
\begin{equation*}
	\Supp^{-1}\Phi=\conditionalset{M\in\QCoh X}{\Supp M\subset\Phi}.
\end{equation*}

\begin{Theorem}\label{BijectionBetweenLocSubcatsAndSpecializationClosedSubsetsForLocNoethSch}
	Let $X$ be a locally noetherian scheme. There is a bijection
	\begin{equation*}
		\setwithspace{\text{localizing subcategories of }\QCoh X}\to\setwithspace{\text{specialization-closed subsets of }X}
	\end{equation*}
	given by $\mcX\mapsto\Supp\mcX$. Its inverse is given by $\Phi\mapsto\Supp^{-1}\Phi$.
\end{Theorem}

\begin{proof}
	In \autoref{DescriptionOfASpecOfLocNoethSch} \autoref{ASpecOfLocNoethSch}, we showed that the Grothendieck category $\QCoh X$ has enough atoms and described $\ASpec(\QCoh X)$. Hence the claim follows from \autoref{BijectionBetweenLocSubcatsAndLocSubsets} and \autoref{DescriptionOfASuppForLocNoethSch} \autoref{LocSubsetForLocNoethSch}.
\end{proof}

\begin{Definition}\label{AssForSch}
	Let $X$ be a locally noetherian scheme, and let $M$ be an object in $\QCoh X$. The subset $\Ass M$ of $X$ is defined by
	\begin{equation*}
		\Ass M=\conditionalset{x\in X}{\mfm_{x}\in\Ass_{\OXx} M_{x}}.
	\end{equation*}
	Each element of $\Ass M$ is called an \emph{associated point} of $M$.
\end{Definition}

In order to show that associated atoms are generalizations of associated points defined in \autoref{AssForSch}, we need the following results.

\begin{Proposition}\label{AssAndOpenAffineImmersion}
	Let $\Spec R$ be an open affine subscheme of $X$, and let $i\colon\Spec R\into X$ be the immersion. For every $R$-module $M$, we have $\Ass i_{*}M=i(\Ass_{R}M)$.
\end{Proposition}

\begin{proof}
	\cite[Proposition~3.1.13]{Grothendieck2} and \cite[Proposition~3.1.2]{Grothendieck2}.
\end{proof}

\begin{Lemma}\label{AssOfIndecInjInQCoh}
	For each $x\in X$, we have $\Ass{j_{x}}_{*}E(x)=\set{x}$ and $\Supp{j_{x}}_{*}E(x)=\overline{\set{x}}$.
\end{Lemma}

\begin{proof}
	Let $i\colon\Spec R\into X$ be the immersion of an open affine subscheme such that $x=i(\mfp)$ for some $\mfp\in\Spec R$. Then the morphism $j_{x}$ is the composite of $j\colon\Spec\OXx\cong\Spec R_{\mfp}\to\Spec R$ and $i\colon\Spec R\into X$. By \cite[Theorem~18.4 (vi)]{Matsumura}, $j_{*}E(x)=E_{R}(R/\mfp)$. By \autoref{AssAndOpenAffineImmersion},
	\begin{equation*}
		\Ass {j_{x}}_{*}E(x)=\Ass i_{*}E_{R}\left(\frac{R}{\mfp}\right)=i\left(\Ass_{R}E_{R}\left(\frac{R}{\mfp}\right)\right)=i(\set{\mfp})=\set{x}.
	\end{equation*}
	By the argument in \cite[p.~150]{Matsumura}, for each $\mfq\in\Spec R$, we have $E_{R}(R/\mfp)_{\mfq}=E_{R_{\mfq}}((R/\mfp)_{\mfq})$. Hence we obtain
	\begin{equation*}
		\Supp E_{R}\left(\frac{R}{\mfp}\right)=\conditionalset{\mfq\in\Spec R}{\mfp\subset\mfq}
	\end{equation*}
	and
	\begin{equation*}
		\Supp{j_{x}}_{*}E(x)=\Supp i_{*}\left(E_{R}\left(\frac{R}{\mfp}\right)\right)=\overline{\set{x}}.\qedhere
	\end{equation*}
\end{proof}

\begin{Proposition}\label{AAssIsAssForLocNoethSch}
	Let $M$ be an object in $\QCoh X$. Then the bijection $\uspX\to\ASpec(\QCoh X)$ in \autoref{DescriptionOfASpecOfLocNoethSch} \autoref{ASpecOfLocNoethSch} restricts to a bijection $\Ass M\to\AAss M$.
\end{Proposition}

\begin{proof}
	Assume that $\alpha_{x}\in\AAss M$, and let $i\colon U\into X$ be the immersion of an open affine subscheme with $x\in U$. By \autoref{OpenSubschIsQuotCat} and \autoref{AAssAndASuppAndQuotCat} \autoref{AAssAndASuppAndCanonicalFunctor}, $\alpha_{x}\in\AAss i^{*}M$. By \autoref{AAssAndASuppOfCommRing} and \autoref{AssAndOpenAffineImmersion}, we obtain $x\in\Ass i_{*}i^{*}M$. Since the canonical morphism $M\to i_{*}i^{*}M$ induces an isomorphism $M_{x}\isoto(i_{*}i^{*}M)_{x}$, we deduce that $x\in\Ass M$.
	
	Conversely, assume that $x\in\Ass M$. By \autoref{InjObjInQCoh} \autoref{InjObjInQCohHasIndecDecomp} and \autoref{InjObjInQCoh} \autoref{DescriptionOfIndecInjObj}, there exists a family $\set{x_{\lambda}}_{\lambda\in\Lambda}$ of points in $X$ such that
	\begin{equation*}
		E(M)\cong\bigoplus_{\lambda\in\Lambda}{j_{x_{\lambda}}}_{*}E(x_{\lambda}).
	\end{equation*}
	By \cite[Proposition~3.1.7]{Grothendieck2},
	\begin{equation*}
		x\in\Ass M\subset\Ass E(M)=\bigcup_{\lambda\in\Lambda}\Ass{j_{x_{\lambda}}}_{*}E(x_{\lambda}).
	\end{equation*}
	Hence there exists $\lambda\in\Lambda$ such that $x\in\Ass{j_{x_{\lambda}}}_{*}E(x_{\lambda})$. By \autoref{AssOfIndecInjInQCoh}, $x_{\lambda}=x$. By \autoref{AAssOfUniformObjAndEssSubobj} \autoref{AAssOfEssSubobj}, we deduce that
	\begin{equation*}
		\alpha_{x}\in\AAss{j_{x}}_{*}E(x)\subset\AAss E(M)=\AAss M.\qedhere
	\end{equation*}
\end{proof}

\section{Localization of prelocalizing subcategories and localizing subcategories}
\label{sec:LocOfPrelocSubcatsAndLocSubcats}

In order to classify the prelocalizing subcategories of $\QCoh X$ for a locally noetherian scheme $X$, we show that they are determined by their restrictions to open affine subschemes of $X$. In this section, we prove this claim in a categorical setting (\autoref{SettingOfLocOfPrelocSubcat}). We start with two lemmas, which show the setting includes the case of $\QCoh X$.

\begin{Lemma}\label{AssAndSuppAndSpecializationOrderInQCoh}
	Let $X$ be a locally noetherian scheme, and let $M$ be an object in $\QCoh X$. Then for each $y\in\Supp M$, there exists $x\in\Ass M$ with $y\in\overline{\set{x}}$.
\end{Lemma}

\begin{proof}
	By \autoref{InjObjInQCoh} \autoref{InjObjInQCohHasIndecDecomp} and \autoref{InjObjInQCoh} \autoref{DescriptionOfIndecInjObj}, there exists a family $\set{x_{\lambda}}_{\lambda\in\Lambda}$ of points in $X$ such that
	\begin{equation*}
		E(M)\cong\bigoplus_{\lambda\in\Lambda}{j_{x_{\lambda}}}_{*}E(x_{\lambda}).
	\end{equation*}
	Then
	\begin{equation*}
		y\in\Supp M\subset\Supp E(M)=\bigcup_{\lambda\in\Lambda}\Supp{j_{x_{\lambda}}}_{*}E(x_{\lambda}).
	\end{equation*}
	By \autoref{AssOfIndecInjInQCoh}, $y\in\Supp{j_{x_{\lambda}}}_{*}E(x_{\lambda})=\overline{\set{x_{\lambda}}}$ for some $\lambda\in\Lambda$. By \autoref{AAssIsAssForLocNoethSch} and \autoref{AAssOfUniformObjAndEssSubobj} \autoref{AAssOfEssSubobj}, we obtain
	\begin{equation*}
		\Ass M=\Ass E(M)=\bigcup_{\lambda\in\Lambda}\Ass{j_{x_{\lambda}}}_{*}E(x_{\lambda})=\conditionalset{x_{\lambda}}{\lambda\in\Lambda}.
	\end{equation*}
	Therefore the claim follows.
\end{proof}

\begin{Lemma}\label{LocOfModOverCommRingAndPrelocSubcat}
	Let $R$ be a commutative ring, and let $S$ be a multiplicatively closed subset of $R$. Let $M$ be an $R$-module. Then the $R$-module $M_{S}$ is a quotient object of the direct sum of some copies of $M$. In particular, for every $\mfp\in\Spec R$, the $R$-module $M_{\mfp}$ belongs to the prelocalizing subcategory $\preloc{M}$ of $\Mod R$.
\end{Lemma}

\begin{proof}
	For each $s\in S$, the image of the $R$-homomorphism $M\to M_{S}$ given by $x\mapsto xs^{-1}$ is $Ms^{-1}$. Hence the $R$-submodule $Ms^{-1}$ of $M_{S}$ is a quotient $R$-module of $M$. Since
	\begin{equation*}
		\bigoplus_{s\in S}Ms^{-1}\onto\sum_{s\in S}Ms^{-1}=M_{S},
	\end{equation*}
	the claim follows.
\end{proof}

In the rest of this section, we investigate localizations of prelocalizing subcategories in the following setting.

\begin{Setting}\label{SettingOfLocOfPrelocSubcat}
	Let $\mcA$ be a Grothendieck category with enough atoms, and let $\set{\mcX_{\lambda}}_{\lambda\in\Lambda}$ be a family of localizing subcategories of $\mcA$. For each $\lambda\in\Lambda$, let $\mcU_{\lambda}=\mcA/\mcX_{\lambda}$. Denote the canonical functors and their right adjoints by
	\begin{itemize}
		\item $F_{\lambda}\colon\mcA\to\mcU_{\lambda}$ and $G_{\lambda}\colon\mcU_{\lambda}\to\mcA$ for each $\lambda\in\Lambda$,
		\item $F^{\lambda}_{\mu}\colon\mcU_{\lambda}\to\mcU_{\mu}$ and $G^{\lambda}_{\mu}\colon\mcU_{\mu}\to\mcU_{\lambda}$ for each $\lambda,\mu\in\Lambda$ with $U_{\mu}\subset U_{\lambda}$,
		\item $F_{\alpha}\colon\mcA\to\mcA_{\alpha}$ and $G_{\alpha}\colon\mcA_{\alpha}\to\mcA$ for each $\alpha\in\ASpec\mcA$,
		\item $F^{\lambda}_{\alpha}\colon\mcU_{\lambda}\to(\mcU_{\lambda})_{\alpha}$ and $G^{\lambda}_{\alpha}\colon(\mcU_{\lambda})_{\alpha}\to\mcU_{\lambda}$ for each $\lambda\in\Lambda$ and $\alpha\in\ASpec\,\mcU_{\lambda}$. (Note that $(\mcU_{\lambda})_{\alpha}=\mcA_{\alpha}$.)
	\end{itemize}
	
	We assume the following properties.
	\begin{enumerate}
		\item\label{SettingOfLocOfPrelocSubcat:Covering} It holds that
		\begin{equation*}
			\ASpec\mcA=\bigcup_{\lambda\in\Lambda}\ASpec\,\mcU_{\lambda}.
		\end{equation*}
		Moreover, for each $\lambda_{1},\lambda_{2}\in\Lambda$ and $\alpha\in\ASpec\,\mcU_{\lambda_{1}}\cap\ASpec\,\mcU_{\lambda_{2}}$, there exists $\mu\in\Lambda$ such that
		\begin{equation*}
			\alpha\in\ASpec\,\mcU_{\mu}\subset\ASpec\,\mcU_{\lambda_{1}}\cap\ASpec\,\mcU_{\lambda_{2}}.
		\end{equation*}
		In other words, the family $\set{\ASpec\,\mcU_{\lambda}}_{\lambda\in\Lambda}$ satisfies the axiom of open basis of $\ASpec\mcA$.\footnote{However, we regard $\ASpec\mcA$ as a topological space only by the localizing topology. (See \autoref{LocSubsetsFormOpenSubsets}.)}
		\item\label{SettingOfLocOfPrelocSubcat:AAssAndASupp} For each object $M$ in $\mcA$ and $\beta\in\ASupp M$, there exists $\alpha\in\AAss M$ with $\alpha\leq\beta$.
		\item\label{SettingOfLocOfPrelocSubcat:AAssAndPrelocSubcatLocally} Let $\lambda\in\Lambda$, and let $M'$ be an object in $\mcU_{\lambda}$ and $\alpha\in\ASpec\,\mcU_{\lambda}$. Then the object $G^{\lambda}_{\alpha}F^{\lambda}_{\alpha}(M')$ belongs to $\preloc{M'}$.
	\end{enumerate}
\end{Setting}

For a locally noetherian scheme $X$, let $\set{U_{\lambda}}_{\lambda\in\Lambda}$ be an open affine basis of $X$ (that is, an open basis of $X$ consisting of affine subsets). Then \autoref{AssAndSuppAndSpecializationOrderInQCoh} and \autoref{LocOfModOverCommRingAndPrelocSubcat} show that the Grothendieck category $\QCoh X$ together with $\set{\QCoh U_{\lambda}}_{\lambda\in\Lambda}$ satisfies the conditions in \autoref{SettingOfLocOfPrelocSubcat}.

We assume \autoref{SettingOfLocOfPrelocSubcat} in the rest of this section.

We show that every quotient category of $\mcA$ also satisfies the same conditions.

\begin{Proposition}\label{QuotCatSatisfiesSetting}
	Let $\mcX$ be a localizing subcategory of $\mcA$. Denote the canonical functor by $F\colon\mcA\to\mcA/\mcX$ and its right adjoint by $G\colon\mcA/\mcX\to\mcA$. Then the Grothendieck category $\mcA/\mcX$ together with the family $\set{\loc{F(\mcX_{\lambda})}}_{\lambda\in\Lambda}$ of localizing subcategories of $\mcA/\mcX$ also satisfies the conditions in \autoref{SettingOfLocOfPrelocSubcat}. In particular, for every $\alpha\in\ASpec\mcA$, the Grothendieck category $\mcA_{\alpha}$ together with $\set{\loc{(\mcX_{\lambda})_{\alpha}}}_{\lambda\in\Lambda}$ satisfies the conditions in \autoref{SettingOfLocOfPrelocSubcat}.\footnote{It is shown in \autoref{CommutativityAndImageOfExtOfPrelocSubcats} \autoref{ImageOfLocSubcat} that $\loc{F(\mcX_{\lambda})}=F(\mcX_{\lambda})$. In particular, $\loc{(\mcX_{\lambda})_{\alpha}}=(\mcX_{\lambda})_{\alpha}$.}
\end{Proposition}

\begin{proof}
	By \autoref{QuotCatOfGrothCatWithEnoughAtoms} \autoref{QuotCatOfGrothCatWithEnoughAtomsHasEnoughAtoms}, the Grothendieck category $\mcA/\mcX$ has enough atoms.
	
	\autoref{SettingOfLocOfPrelocSubcat:Covering} By \autoref{ASpecOfQuotOfLocSubcat},
	\begin{equation*}
		\ASpec\frac{\mcA/\mcX}{\loc{F(\mcX_{\lambda})}}=\ASpec\frac{\mcA}{\mcX_{\lambda}}\cap\ASpec\frac{\mcA}{\mcX}.
	\end{equation*}
	
	\autoref{SettingOfLocOfPrelocSubcat:AAssAndASupp} Let $M'$ be an object in $\mcA/\mcX$, and let $\beta\in\ASupp M'$. By \autoref{AAssAndASuppAndQuotCat} \autoref{AAssAndASuppAndSectionFunctor}, $\beta\in\ASupp G(M')$. Hence there exists $\alpha\in\AAss G(M')=\AAss M'$ with $\alpha\leq\beta$.
	
	\autoref{SettingOfLocOfPrelocSubcat:AAssAndPrelocSubcatLocally} Let $\lambda\in\Lambda$. By \autoref{QuotOfLocSubcat} \autoref{DescriptionOfQuotOfLocSubcat} and \autoref{QuotOfLocSubcat} \autoref{QuotOfQuotCat},
	\begin{equation*}
		\frac{\mcA/\mcX}{\loc{F(\mcX_{\lambda})}}\cong\frac{\mcA}{\loc{\mcX_{\lambda}\cup\mcX}}\cong\frac{\mcA/\mcX_{\lambda}}{\loc{F_{\lambda}(\mcX)}}.
	\end{equation*}
	Let $\mcU'_{\lambda}:=\mcU_{\lambda}/\loc{F_{\lambda}(\mcX)}$. Denote the canonical functors by $F'\colon\mcU_{\lambda}\to\mcU'_{\lambda}$, $F'_{\alpha}\colon\mcU'_{\lambda}\to\mcA_{\alpha}$, and their right adjoints by $G'\colon\mcU'_{\lambda}\to\mcU_{\lambda}$, $G'_{\alpha}\colon\mcA_{\alpha}\to\mcU'_{\lambda}$, respectively. Let $M''$ be an object in $\mcU'_{\lambda}$, and let $\alpha\in\ASpec\,\mcU'_{\lambda}$. Then by the assumption, the object $G^{\lambda}_{\alpha}F^{\lambda}_{\alpha}G'(M'')$ belongs to $\preloc{G'(M'')}$. Since $F'$ is exact, the object $F'G^{\lambda}_{\alpha}F^{\lambda}_{\alpha}G'(M'')$ belongs to $\preloc{F'G'(M'')}=\preloc{M''}$. Since
	\begin{equation*}
		F'G^{\lambda}_{\alpha}F^{\lambda}_{\alpha}G'(M'')\cong F'G'G'_{\alpha}F'_{\alpha}F'G'(M'')\cong G'_{\alpha}F'_{\alpha}(M''),
	\end{equation*}
	the claim follows.
\end{proof}

Under the assumptions of \autoref{SettingOfLocOfPrelocSubcat}, we can show a complemental fact on associated atoms in a quotient category.

\begin{Lemma}\label{AAssOfObjInQuotCatInSetting}
	Let $\mcX$ be a localizing subcategory of $\mcA$. Denote the canonical functor by $F\colon\mcA\to\mcA/\mcX$ and its right adjoint by $G\colon\mcA/\mcX\to\mcA$. For every object $M$ in $\mcA$,
	\begin{equation*}
		\AAss F(M)=\AAss M\setminus\ASupp\mcX.
	\end{equation*}
	In particular, for every $\alpha\in\ASpec\mcA$,
	\begin{equation*}
		\AAss M_{\alpha}=\AAss M\cap\Lambda(\alpha).
	\end{equation*}
\end{Lemma}

\begin{proof}
	By \autoref{AAssAndASuppAndQuotCat},
	\begin{equation*}
		\AAss GF(M)=\AAss F(M)\supset\AAss M\setminus\ASupp\mcX.
	\end{equation*}
	Let $\eta\colon 1_{\mcA}\to GF$ be the unit morphism and $\beta\in\AAss GF(M)$. Note that $\beta\notin\ASupp\mcX$. By \autoref{PropertiesOfQuotCat} \autoref{DescriptionOfUnit}, the subobject $L:=\Ker\eta_{M}$ of $M$ belongs to $\mcX$, and $\Im\eta_{M}$ is an essential subobject of $GF(M)$. By \autoref{AAssOfUniformObjAndEssSubobj} \autoref{AAssOfEssSubobj}, $\beta\in\AAss(\Im\eta_{M})=\AAss(M/L)$. Hence there exists a subobject $L'$ of $M$ with $L\subset L'$ such that $L'/L$ is a monoform object representing $\beta$. Since $\beta\in\ASupp L'$, by \autoref{SettingOfLocOfPrelocSubcat} \autoref{SettingOfLocOfPrelocSubcat:AAssAndASupp}, there exists $\alpha\in\AAss L'$ with $\alpha\leq\beta$. Since $\beta\notin\ASupp\mcX$, it holds that $\alpha\notin\ASupp L$ by \autoref{CharacterizationOfSpecializationOrder}. By \autoref{AAssAndASuppAndShortExactSeq} and \autoref{AAssOfUniformObjAndEssSubobj} \autoref{AAssOfUniformObj}, $\alpha\in\AAss(L'/L)=\set{\beta}$. Therefore $\beta=\alpha\in\AAss L'\subset\AAss M$.
\end{proof}

We show two lemmas as parts of the proof of \autoref{LocOfObjAtAAssAndPrelocSubcat}. It is useful to determine whether an object belongs to a given prelocalizing subcategory.

\begin{Lemma}\label{LocOfUniformObjAtAAssAndPrelocSubcatLocally}
	Let $\lambda\in\Lambda$, and let $\mcY'$ be a prelocalizing subcategory of $\mcU_{\lambda}$. Let $U'$ be a uniform object in $\mcU_{\lambda}$ with $\AAss U'=\set{\alpha}$. If $U'_{\alpha}$ belongs to $\mcY'_{\alpha}$, then $U'$ belongs to $\mcY'$.
\end{Lemma}

\begin{proof}
	There exists an object $N'$ in $\mcU_{\lambda}$ which belongs to $\mcY'$ such that $U'_{\alpha}\cong N'_{\alpha}$. By \autoref{SettingOfLocOfPrelocSubcat} \autoref{SettingOfLocOfPrelocSubcat:AAssAndPrelocSubcatLocally}, the object $G^{\lambda}_{\alpha}F^{\lambda}_{\alpha}(U')\cong G^{\lambda}_{\alpha}F^{\lambda}_{\alpha}(N')$ belongs to $\mcY'$. Let $\eta\colon 1_{\mcU_{\lambda}}\to G^{\lambda}_{\alpha}F^{\lambda}_{\alpha}$ be the unit morphism. Then by \autoref{PropertiesOfQuotCat} \autoref{DescriptionOfUnit}, $\alpha\notin\ASupp(\Ker\eta_{U'})$. If $\Ker\eta_{U'}\neq 0$, then by \autoref{AAssOfUniformObjAndEssSubobj} \autoref{AAssOfEssSubobj}, $\alpha\in\AAss(\Ker\eta_{U'})\subset\ASupp(\Ker\eta_{U'})$. This is a contradiction. Hence $\eta_{U'}$ is a monomorphism. The object $U'$ belongs to $\mcY'$.
\end{proof}

\begin{Lemma}\label{LocOfUniformObjAtAAssAndPrelocSubcat}
	Let $\mcY$ be a prelocalizing subcategory of $\mcA$, and let $U$ be a uniform object in $\mcA$ with $\AAss U=\set{\alpha}$. If $U_{\alpha}$ belongs to $\mcY_{\alpha}$, then $U$ belongs to $\mcY$.
\end{Lemma}

\begin{proof}
	Let $L$ be the largest subobject of $U$ which belongs to $\mcY$. Assume that $L\subsetneq U$. Then by \autoref{GrothCatWithEnoughAtomsIsLocMonoform}, there exists $\beta\in\AAss(U/L)$. By \autoref{SettingOfLocOfPrelocSubcat} \autoref{SettingOfLocOfPrelocSubcat:AAssAndASupp}, $\alpha\leq\beta$. By \autoref{SettingOfLocOfPrelocSubcat} \autoref{SettingOfLocOfPrelocSubcat:Covering}, there exists $\lambda\in\Lambda$ such that $\beta\in\ASpec\,\mcU_{\lambda}=\ASpec\mcA\setminus\ASupp\mcX_{\lambda}$. Then by \autoref{CharacterizationOfSpecializationOrder}, we also have $\alpha\in\ASpec\,\mcU_{\lambda}$. By a similar argument to that in the proof of \autoref{LocOfUniformObjAtAAssAndPrelocSubcatLocally}, the canonical morphism $U\to G_{\lambda}F_{\lambda}(U)$ is a monomorphism, and $U$ is $\mcX_{\lambda}$-torsionfree. By \autoref{PropertiesOfImageOfObj} \autoref{UniformnessOfImageOfObj}, the object $F_{\lambda}(U)$ is uniform, and $\AAss F_{\lambda}(U)=\set{\alpha}$ by \autoref{AAssOfUniformObjAndEssSubobj} \autoref{AAssOfUniformObj}. Since $F_{\lambda}(U)_{\alpha}=U_{\alpha}$ belongs to $\mcY_{\alpha}=F_{\lambda}(\mcY)_{\alpha}$, by \autoref{LocOfUniformObjAtAAssAndPrelocSubcatLocally}, the object $F_{\lambda}(U)$ belongs to $F_{\lambda}(\mcY)$. We obtain an object $N$ in $\mcA$ which belongs to $\mcY$ such that $F_{\lambda}(N)\cong F_{\lambda}(U)$. Let $V$ be the image of the composite of the canonical morphism $N\to G_{\lambda}F_{\lambda}(N)$ and $G_{\lambda}F_{\lambda}(N)\isoto G_{\lambda}F_{\lambda}(U)$. By \autoref{PropertiesOfQuotCat} \autoref{DescriptionOfUnit}, the object $G_{\lambda}F_{\lambda}(U)/V$ belongs to $\mcX_{\lambda}$. Hence
	\begin{equation*}
		\frac{G_{\lambda}F_{\lambda}(U)}{U\cap V}\into\frac{G_{\lambda}F_{\lambda}(U)}{U}\oplus\frac{G_{\lambda}F_{\lambda}(U)}{V}\in\mcX_{\lambda}.
	\end{equation*}
	Since $U\cap V$ belongs to $\mcY$, we have $U\cap V\subset L$ by the maximality of $L$. Hence $U/L$ also belongs to $\mcX_{\lambda}$, and $\beta\in\ASupp(U/L)\subset\ASupp\mcX_{\lambda}$. This is a contradiction. Therefore $L=U$.
\end{proof}

\begin{Theorem}\label{LocOfObjAtAAssAndPrelocSubcat}
	Assume \autoref{SettingOfLocOfPrelocSubcat}. Let $\mcY$ be a prelocalizing subcategory of $\mcA$, and let $M$ be an object in $\mcA$. If $M_{\alpha}$ belongs to $\mcY_{\alpha}$ for every $\alpha\in\AAss M$, then $M$ belongs to $\mcY$.
\end{Theorem}

\begin{proof}
	Since $\mcA$ has enough atoms, there exists a family $\set{\alpha_{\omega}}_{\omega\in\Omega}$ of elements of $\ASpec\mcA$ such that
	\begin{equation*}
		E(M)\cong\bigoplus_{\omega\in\Omega}E(\alpha_{\omega}).
	\end{equation*}
	Let $\mcZ=\preloc{M}$. For each $\omega\in\Omega$, let $L_{\omega}$ be the largest subobject of $E(\alpha_{\omega})$ which belongs to $\mcZ$. Then by \autoref{LargestSubobjOfDirectSumWrtPrelocSubcat}, $M\subset\bigoplus_{\omega\in\Omega}L_{\omega}$. Since $L_{\omega}$ is uniform for each $\omega\in\Omega$, by \autoref{AAssOfUniformObjAndEssSubobj},
	\begin{equation*}
		\set{\alpha_{\omega}}=\AAss L_{\omega}\subset\AAss E(M)=\AAss M.
	\end{equation*}
	By \autoref{QuotOfPrelocSubcat} \autoref{QuotOfPrelocSubcatIsPrelocSubcat}, it is straightforward to show that $\mcZ_{\alpha_{\omega}}=\preloc{M_{\alpha_{\omega}}}$. Hence by the assumption, $\mcZ_{\alpha_{\omega}}\subset\mcY_{\alpha_{\omega}}$. Since $L_{\omega}$ belongs to $\mcZ$, the object $(L_{\omega})_{\alpha_{\omega}}$ belongs to $\mcY_{\alpha_{\omega}}$. By \autoref{LocOfUniformObjAtAAssAndPrelocSubcat}, we deduce that $L_{\omega}$ belongs to $\mcY$. Therefore the subobject $M$ of $\bigoplus_{\omega\in\Omega}L_{\omega}$ also belongs to $\mcY$.
\end{proof}

The following results are consequences of \autoref{LocOfObjAtAAssAndPrelocSubcat}.

\begin{Proposition}\label{LocAndPrelocSubcat}
	Let $\mcX$ be a localizing subcategory of $\mcA$. Denote the canonical functor by $F\colon\mcA\to\mcA/\mcX$ and its right adjoint by $G\colon\mcA/\mcX\to\mcA$. Then for every object $M$ in $\mcA$, the object $GF(M)$ belongs to $\preloc{M}$.
\end{Proposition}

\begin{proof}
	Let $\eta\colon 1_{\mcA}\to GF$ be the unit morphism. Let $\alpha\in\AAss GF(M)$. By \autoref{AAssAndASuppAndQuotCat} \autoref{AAssAndASuppAndSectionFunctor},
	\begin{equation*}
		\alpha\in\AAss GF(M)=\AAss F(M)\subset\ASpec\frac{\mcA}{\mcX}=\ASpec\mcA\setminus\ASupp\mcX.
	\end{equation*}
	By \autoref{PropertiesOfQuotCat} \autoref{DescriptionOfUnit}, the objects $\Ker\eta_{M}$ and $\Cok\eta_{M}$ belong to $\mcX$. By applying $(-)_{\alpha}$ to the exact sequence
	\begin{equation*}
		0\to\Ker\eta_{M}\to M\to GF(M)\to\Cok\eta_{M}\to 0,
	\end{equation*}
	we obtain the isomorphism $M_{\alpha}\isoto GF(M)_{\alpha}$. Hence $GF(M)_{\alpha}$ belongs to $(\preloc{M})_{\alpha}$. By \autoref{LocOfObjAtAAssAndPrelocSubcat}, we deduce that $GF(M)$ belongs to $\preloc{M}$.
\end{proof}

\begin{Proposition}\label{AtomAndAtomicObjBelongToPrelocSubcat}
	Let $\mcY$ be a prelocalizing subcategory of $\mcA$ and $\alpha\in\ASpec\mcA$. Then $\alpha\in\ASupp\mcY$ if and only if $H(\alpha)$ belongs to $\mcY$.
\end{Proposition}

\begin{proof}
	If $H(\alpha)$ belongs to $\mcY$, then $\alpha=\overline{H(\alpha)}\in\ASupp\mcY$.
	
	Assume $\alpha\in\ASupp\mcY$. Then there exists a monoform object $H$ in $\mcA$ with $\overline{H}=\alpha$ such that $H$ belongs to $\mcY$. By \autoref{LocAndPrelocSubcat}, the object $G_{\alpha}F_{\alpha}(H)$ belongs to $\mcY$. By the proof of \autoref{SimpleObjAndAtomicObj} \autoref{AtomicObjIsSimpleObjInLocalizedCat}, the object $G_{\alpha}F_{\alpha}(H)$ is isomorphic to $H(\alpha)$.
\end{proof}

We show the main result in this section.

\begin{Theorem}\label{LocOfPrelocSubcat}
	Assume \autoref{SettingOfLocOfPrelocSubcat}. Then there exist bijections between the following sets.
	\begin{enumerate}
		\item\label{LocOfPrelocSubcat:PrelocSubcats} The set of prelocalizing subcategories of $\mcA$.
		\item\label{LocOfPrelocSubcat:PrelocSubcatsOfQuotCats} The set of families $\set{\mcY_{\lambda}\subset\mcU_{\lambda}}_{\lambda\in\Lambda}$ of prelocalizing subcategories such that $F^{\lambda}_{\mu}(\mcY_{\lambda})=\mcY_{\mu}$ for each $\lambda,\mu\in\Lambda$ with $\ASpec\,\mcU_{\mu}\subset\ASpec\,\mcU_{\lambda}$.
		\item\label{LocOfPrelocSubcat:LocalizedPrelocSubcats} The set of families $\set{\mcY(\alpha)\subset\mcA_{\alpha}}_{\alpha\in\ASpec\mcA}$ of prelocalizing subcategories such that $\mcY(\beta)_{\alpha}=\mcY(\alpha)$ for each $\alpha,\beta\in\ASpec\mcA$ with $\alpha\leq\beta$.
	\end{enumerate}
	The correspondences are given as follows.
	\begin{align*}
		&\autoref{LocOfPrelocSubcat:PrelocSubcats}\ \mcY\mapsto
		\begin{cases}
			\autoref{LocOfPrelocSubcat:PrelocSubcatsOfQuotCats}\ \set{F_{\lambda}(\mcY)}_{\lambda\in\Lambda}\\
			\autoref{LocOfPrelocSubcat:LocalizedPrelocSubcats}\ \set{\mcY_{\alpha}}_{\alpha\in\ASpec\mcA},
		\end{cases}\\
		&\autoref{LocOfPrelocSubcat:PrelocSubcatsOfQuotCats}\ \set{\mcY_{\lambda}}_{\lambda\in\Lambda}\mapsto
		\begin{cases}
			\autoref{LocOfPrelocSubcat:PrelocSubcats}\ \displaystyle\bigcap_{\lambda\in\Lambda}F_{\lambda}^{-1}(\mcY_{\lambda})\\
			\autoref{LocOfPrelocSubcat:LocalizedPrelocSubcats}\ \set{\mcY(\alpha)}_{\alpha\in\ASpec\mcA},\\
			\phantom{\autoref{LocOfPrelocSubcat:LocalizedPrelocSubcats}}\text{ where }\mcY(\alpha)=(\mcY_{\lambda})_{\alpha}\text{ for }\lambda\in\Lambda\text{ with }\alpha\in\ASpec\,\mcU_{\lambda},
		\end{cases}\\
		&\autoref{LocOfPrelocSubcat:LocalizedPrelocSubcats}\ \set{\mcY(\alpha)}_{\alpha\in\ASpec\mcA}\mapsto
		\begin{cases}
			\autoref{LocOfPrelocSubcat:PrelocSubcats}\ \displaystyle\bigcap_{\alpha\in\ASpec\mcA}F_{\alpha}^{-1}(\mcY(\alpha))\\
			\autoref{LocOfPrelocSubcat:PrelocSubcatsOfQuotCats}\ \displaystyle\set[\Bigg]{\bigcap_{\alpha\in\ASpec\,\mcU_{\lambda}}(F^{\lambda}_{\alpha})^{-1}(\mcY(\alpha))}_{\lambda\in\Lambda}.
		\end{cases}
	\end{align*}
\end{Theorem}

\begin{proof}
	(\autoref{LocOfPrelocSubcat:PrelocSubcats}$\leftrightarrow$\autoref{LocOfPrelocSubcat:PrelocSubcatsOfQuotCats}) Let $\mcY$ be a prelocalizing subcategory of $\mcA$. It is obvious that $\mcY\subset\bigcap_{\lambda\in\Lambda}F_{\lambda}^{-1}F_{\lambda}(\mcY)$. Let $M$ be an object in $\mcA$ which belongs to $\bigcap_{\lambda\in\Lambda}F_{\lambda}^{-1}F_{\lambda}(\mcY)$. For each $\alpha\in\AAss M$, by \autoref{SettingOfLocOfPrelocSubcat} \autoref{SettingOfLocOfPrelocSubcat:Covering}, there exists $\lambda\in\Lambda$ such that $\alpha\in\ASpec\,\mcU_{\lambda}$. Then
	\begin{equation*}
		M_{\alpha}=F_{\lambda}(M)_{\alpha}\in F_{\lambda}(\mcY)_{\alpha}=\mcY_{\alpha}.
	\end{equation*}
	By \autoref{LocOfObjAtAAssAndPrelocSubcat}, the object $M$ belongs to $\mcY$. We obtain $\mcY=\bigcap_{\lambda\in\Lambda}F_{\lambda}^{-1}F_{\lambda}(\mcY)$.
	
	Let $\set{\mcY_{\lambda}}_{\lambda\in\Lambda}$ be an element of \autoref{LocOfPrelocSubcat:PrelocSubcatsOfQuotCats} and $\mcY:=\bigcap_{\lambda\in\Lambda}F_{\lambda}^{-1}(\mcY_{\lambda})$. It is obvious that $F_{\lambda}(\mcY)\subset\mcY_{\lambda}$. Let $M'$ be an object in $\mcU_{\lambda}$ which belongs to $\mcY_{\lambda}$. We show that $F_{\lambda'}G_{\lambda}(M')$ belongs to $\mcY_{\lambda'}$ for each $\lambda'\in\Lambda$. For each $\beta\in\AAss F_{\lambda'}G_{\lambda}(M')$, by \autoref{AAssOfObjInQuotCatInSetting} and \autoref{AAssAndASuppAndQuotCat} \autoref{AAssAndASuppAndSectionFunctor},
	\begin{equation*}
		\beta\in\AAss M'\cap\ASpec\,\mcU_{\lambda'}\subset\ASpec\,\mcU_{\lambda}\cap\ASpec\,\mcU_{\lambda'}.
	\end{equation*}
	Hence there exists $\mu\in\Lambda$ such that
	\begin{equation*}
		\beta\in\ASpec\,\mcU_{\mu}\subset\ASpec\,\mcU_{\lambda}\cap\ASpec\,\mcU_{\lambda'}.
	\end{equation*}
	Since the object
	\begin{equation*}
		F_{\lambda'}G_{\lambda}(M')_{\beta}=G_{\lambda}(M')_{\beta}=F_{\lambda}G_{\lambda}(M')_{\beta}=M'_{\beta}
	\end{equation*}
	belongs to
	\begin{equation*}
		(\mcY_{\lambda})_{\beta}=F^{\lambda}_{\mu}(\mcY_{\lambda})_{\beta}=(\mcY_{\mu})_{\beta}=F^{\lambda'}_{\mu}(\mcY_{\lambda'})_{\beta}=(\mcY_{\lambda'})_{\beta},
	\end{equation*}
	by \autoref{QuotCatSatisfiesSetting} and \autoref{LocOfObjAtAAssAndPrelocSubcat}, the object $F_{\lambda'}G_{\lambda}(M')$ belongs to $\mcY_{\lambda'}$. Hence $G_{\lambda}(M')$ belongs to $\mcY$, and $M'\cong F_{\lambda}G_{\lambda}(M')$ belongs to $F_{\lambda}(\mcY)$. We obtain $F_{\lambda}(\mcY)=\mcY_{\lambda}$.
	
	(Well-definedness of \autoref{LocOfPrelocSubcat:PrelocSubcatsOfQuotCats}$\to$\autoref{LocOfPrelocSubcat:LocalizedPrelocSubcats}) Let $\set{\mcY_{\lambda}}_{\lambda\in\Lambda}$ be an element of \autoref{LocOfPrelocSubcat:PrelocSubcatsOfQuotCats} and $\alpha\in\ASpec\mcA$. Let $\lambda_{1},\lambda_{2}\in\Lambda$ such that $\alpha\in\ASpec\,\mcU_{\lambda_{i}}$ for each $i=1,2$. Then by \autoref{SettingOfLocOfPrelocSubcat} \autoref{SettingOfLocOfPrelocSubcat:Covering}, there exists $\mu\in\Lambda$ such that
	\begin{equation*}
		\alpha\in\ASpec\,\mcU_{\mu}\subset\ASpec\,\mcU_{\lambda_{1}}\cap\ASpec\,\mcU_{\lambda_{2}}.
	\end{equation*}
	Hence
	\begin{equation*}
		(\mcY_{\lambda_{1}})_{\alpha}=F^{\lambda_{1}}_{\mu}(\mcY_{\lambda_{1}})_{\alpha}=(\mcY_{\mu})_{\alpha}=F^{\lambda_{2}}_{\mu}(\mcY_{\lambda_{2}})_{\alpha}=(\mcY_{\lambda_{2}})_{\alpha}.
	\end{equation*}
	
	(\autoref{LocOfPrelocSubcat:PrelocSubcatsOfQuotCats}$\leftrightarrow$\autoref{LocOfPrelocSubcat:LocalizedPrelocSubcats}) Let $\set{\mcY_{\lambda}}_{\lambda\in\Lambda}$ be an element of \autoref{LocOfPrelocSubcat:PrelocSubcatsOfQuotCats}. For each $\lambda\in\Lambda$, let
	\begin{equation*}
		\widetilde{\mcY}_{\lambda}:=\bigcap_{\alpha\in\ASpec\,\mcU_{\lambda}}(F^{\lambda}_{\alpha})^{-1}F^{\lambda}_{\alpha}(\mcY_{\lambda}).
	\end{equation*}
	Then $\mcY_{\lambda}\subset\widetilde{\mcY}_{\lambda}$. Let $M'$ be an object in $\mcU_{\lambda}$ which belongs to $\widetilde{\mcY}_{\lambda}$. For each $\alpha\in\AAss M'$, we have $M'_{\alpha}\in(\mcY_{\lambda})_{\alpha}$. Hence by \autoref{QuotCatSatisfiesSetting} and \autoref{LocOfObjAtAAssAndPrelocSubcat}, the object $M'$ belongs to $\mcY_{\lambda}$, and we obtain $\mcY_{\lambda}=\widetilde{\mcY}_{\lambda}$.
	
	Let $\set{\mcY(\alpha)}_{\alpha\in\ASpec\mcA}$ be an element of \autoref{LocOfPrelocSubcat:LocalizedPrelocSubcats}. For each $\lambda\in\Lambda$, let
	\begin{equation*}
		\mcY_{\lambda}:=\bigcap_{\alpha\in\ASpec\,\mcU_{\lambda}}(F^{\lambda}_{\alpha})^{-1}(\mcY(\alpha)).
	\end{equation*}
	For each $\alpha\in\ASpec\mcA$, by \autoref{SettingOfLocOfPrelocSubcat} \autoref{SettingOfLocOfPrelocSubcat:Covering}, there exists $\mu\in\Lambda$ such that $\alpha\in\ASpec\,\mcU_{\mu}$. It is obvious that $F^{\mu}_{\alpha}(\mcY_{\lambda})\subset\mcY(\alpha)$. Let $M''$ be an object in $\mcA_{\alpha}$ which belongs to $\mcY(\alpha)$. We show that $F^{\lambda}_{\beta}G^{\lambda}_{\alpha}(M'')$ belongs to $\mcY(\beta)$ for each $\beta\in\ASpec\,\mcU_{\lambda}$. For each $\gamma\in\AAss F^{\lambda}_{\beta}G^{\lambda}_{\alpha}(M'')$, by \autoref{AAssOfObjInQuotCatInSetting} and \autoref{AAssAndASuppAndQuotCat} \autoref{AAssAndASuppAndSectionFunctor},
	\begin{equation*}
		\gamma\in\AAss M''\cap\Lambda(\beta)\subset\Lambda(\alpha)\cap\Lambda(\beta).
	\end{equation*}
	Since the object
	\begin{equation*}
		F^{\lambda}_{\beta}G^{\lambda}_{\alpha}(M'')_{\gamma}=G^{\lambda}_{\alpha}(M'')_{\gamma}=F^{\lambda}_{\alpha}G^{\lambda}_{\alpha}(M'')_{\gamma}=M''_{\gamma}
	\end{equation*}
	belongs to
	\begin{equation*}
		\mcY(\alpha)_{\gamma}=\mcY(\gamma)=\mcY(\beta)_{\gamma},
	\end{equation*}
	the object $F^{\lambda}_{\beta}G^{\lambda}_{\alpha}(M'')$ belongs to $\mcY(\beta)$. Hence $G^{\lambda}_{\alpha}(M'')$ belongs to $\mcY_{\lambda}$, and $M''\cong F^{\lambda}_{\alpha}G^{\lambda}_{\alpha}(M'')\in F^{\lambda}_{\alpha}(\mcY_{\lambda})$. We obtain $F^{\mu}_{\beta}(\mcY_{\lambda})=\mcY(\beta)$.
\end{proof}

For a family $\set{\mcY^{\omega}}_{\omega\in\Omega}$ of prelocalizing subcategories of $\mcA$, we can consider the smallest prelocalizing subcategory $\preloc{\bigcup_{\omega\in\Omega}\mcY^{\omega}}$ containing $\mcY^{\omega}$ for every $\omega\in\Omega$ and the intersection $\bigcap_{\omega\in\Omega}\mcY^{\omega}$. These are described in terms of prelocalizing subcategories of quotient categories in the following ways.

\begin{Proposition}\label{LocOfJoinOfPrelocSubcats}
	Assume that the following elements correspond to each other by the bijections in \autoref{LocOfPrelocSubcat} for each $\omega\in\Omega$.
	\begin{enumerate}
		\item $\mcY^{\omega}$.
		\item $\set{\mcY^{\omega}_{\lambda}}_{\lambda\in\Lambda}$.
		\item $\set{\mcY^{\omega}(\alpha)}_{\alpha\in\ASpec\mcA}$.
	\end{enumerate}
	Then the following elements correspond to each other by the bijections.
	\begin{enumerate}
		\item $\preloc{\bigcup_{\omega\in\Omega}\mcY^{\omega}}$.
		\item $\set{\preloc{\bigcup_{\omega\in\Omega}\mcY^{\omega}_{\lambda}}}_{\lambda\in\Lambda}$.
		\item $\set{\preloc{\bigcup_{\omega\in\Omega}\mcY^{\omega}(\alpha)}}_{\alpha\in\ASpec\mcA}$.
	\end{enumerate}
\end{Proposition}

\begin{proof}
	For each $\lambda\in\Lambda$,
	\begin{align*}
		F_{\lambda}\Bigg(\generatedset[\Bigg]{\!\!preloc}{\bigcup_{\omega\in\Omega}\mcY^{\omega}}\Bigg)&=\generatedset[\Bigg]{\!\!preloc}{F_{\lambda}\Bigg(\bigcup_{\omega\in\Omega}\mcY^{\omega}\Bigg)}\\
		&=\generatedset[\Bigg]{\!\!preloc}{\bigcup_{\omega\in\Omega}F_{\lambda}(\mcY^{\omega})}\\
		&=\generatedset[\Bigg]{\!\!preloc}{\bigcup_{\omega\in\Omega}\mcY^{\omega}_{\lambda}}.
	\end{align*}
	It is shown similarly that $(\preloc{\bigcup_{\omega\in\Omega}\mcY^{\omega}})_{\alpha}=\preloc{\bigcup_{\omega\in\Omega}\mcY^{\omega}(\alpha)}$ for each $\alpha\in\ASpec\mcA$.
\end{proof}

\begin{Proposition}\label{LocOfIntersectionOfPrelocSubcats}
	Assume that the following elements correspond to each other by the bijections in \autoref{LocOfPrelocSubcat} for each $\omega\in\Omega$.
	\begin{enumerate}
		\item $\mcY^{\omega}$.
		\item $\set{\mcY^{\omega}_{\lambda}}_{\lambda\in\Lambda}$.
		\item $\set{\mcY^{\omega}(\alpha)}_{\alpha\in\ASpec\mcA}$.
	\end{enumerate}
	Then the following elements correspond to each other by the bijections.
	\begin{enumerate}
		\item $\bigcap_{\omega\in\Omega}\mcY^{\omega}$.
		\item $\set{\bigcap_{\omega\in\Omega}\mcY^{\omega}_{\lambda}}_{\lambda\in\Lambda}$.
		\item $\set{\bigcap_{\omega\in\Omega}\mcY^{\omega}(\alpha)}_{\alpha\in\ASpec\mcA}$.
	\end{enumerate}
\end{Proposition}

\begin{proof}
	Let $\lambda\in\Lambda$. It is obvious that $F_{\lambda}(\bigcap_{\omega\in\Omega}\mcY^{\omega})\subset\bigcap_{\omega\in\Omega}F_{\lambda}(\mcY^{\omega})=\bigcap_{\omega\in\Omega}\mcY^{\omega}_{\lambda}$. Let $M'$ be an object in $\mcU_{\lambda}$ which belongs to $\bigcap_{\omega\in\Omega}\mcY^{\omega}_{\lambda}$. Then for each $\omega\in\Omega$, there exists an object $M_{\omega}$ in $\mcA$ which belongs to $\mcY^{\omega}$ such that $F_{\lambda}(M_{\omega})\cong M'$. By \autoref{LocAndPrelocSubcat}, the object $G_{\lambda}(M')\cong G_{\lambda}F_{\lambda}(M_{\omega})$ belongs to $\preloc{M_{\omega}}$. Hence $G_{\lambda}(M')$ belongs to $\bigcap_{\omega\in\Omega}\mcY^{\omega}$, and $M'\cong F_{\lambda}G_{\lambda}(M')\in F_{\lambda}(\bigcap_{\omega\in\Omega}\mcY^{\omega})$. This shows that $F_{\lambda}(\bigcap_{\omega\in\Omega}\mcY^{\omega})=\bigcap_{\omega\in\Omega}\mcY^{\omega}_{\lambda}$.
	
	It is shown similarly that $(\bigcap_{\omega\in\Omega}\mcY^{\omega})_{\alpha}=\bigcap_{\omega\in\Omega}\mcY^{\omega}(\alpha)$ for each $\alpha\in\ASpec\mcA$.
\end{proof}

Families in \autoref{LocOfPrelocSubcat} \autoref{LocOfPrelocSubcat:LocalizedPrelocSubcats} have the following characterization.

\begin{Proposition}\label{DescriptionOfPrelocSubcatsOnLocCats}
	For each family $\set{\mcY(\alpha)\subset\mcA_{\alpha}}_{\alpha\in\ASpec\mcA}$ of prelocalizing subcategories, the following assertions are equivalent.
	\begin{enumerate}
		\item\label{DescriptionOfPrelocSubcatsOnLocCats:WholeCats} There exists a prelocalizing subcategory $\mcY$ of $\mcA$ satisfying $\mcY_{\alpha}=\mcY(\alpha)$ for each $\alpha\in\ASpec\mcA$.
		\item\label{DescriptionOfPrelocSubcatsOnLocCats:QuotCats} For each $\alpha\in\ASpec\mcA$, there exist $\lambda\in\Lambda$ with $\alpha\in\ASpec\,\mcU_{\lambda}$ and a prelocalizing subcategory $\mcY'$ of $\mcU_{\lambda}$ satisfying $\mcY'_{\beta}=\mcY(\beta)$ for each $\beta\in\ASpec\,\mcU_{\lambda}$.
		\item\label{DescriptionOfPrelocSubcatsOnLocCats:LocCats} For each $\alpha,\beta\in\ASpec\mcA$ with $\alpha\leq\beta$, it holds that $\mcY(\beta)_{\alpha}=\mcY(\alpha)$.
	\end{enumerate}
\end{Proposition}

\begin{proof}
	This can be shown straightforwardly by using \autoref{LocOfPrelocSubcat}.
\end{proof}

In order to investigate the localizing subcategories of $\QCoh X$, we improve \autoref{QuotOfPrelocSubcat} under the assumptions of \autoref{SettingOfLocOfPrelocSubcat}.

\begin{Proposition}\label{CommutativityAndImageOfExtOfPrelocSubcats}
	Let $\mcX$ be a localizing subcategory of $\mcA$. Denote the canonical functor by $F\colon\mcA\to\mcA/\mcX$ and its right adjoint by $G\colon\mcA/\mcX\to\mcA$.
	\begin{enumerate}
		\item\label{CommutativityBetweenLocSubcatAndPrelocSubcat} Let $\mcY$ be a prelocalizing subcategory of $\mcA$. Then $\mcY*\mcX\subset\mcX*\mcY$.
		\item\label{ImageOfExtOfPrelocSubcats} Let $\mcY_{1}$ and $\mcY_{2}$ be prelocalizing subcategories of $\mcA$. Then
		\begin{equation*}
			F(\mcY_{1}*\mcY_{2})=F(\mcY_{1})*F(\mcY_{2}).
		\end{equation*}
		\item\label{ImageOfLocSubcat} Let $\mcY$ be a localizing subcategory of $\mcA$. Then $F(\mcY)$ is a localizing subcategory of $\mcA/\mcX$.
	\end{enumerate}
\end{Proposition}

\begin{proof}
	\autoref{CommutativityBetweenLocSubcatAndPrelocSubcat} Let $M$ be an object in $\mcA$ which belongs to $\mcY*\mcX$. Then there exists an exact sequence
	\begin{equation*}
		0\to L\to M\to N\to 0
	\end{equation*}
	in $\mcA$, where $L$ belongs to $\mcY$, and $N$ belongs to $\mcX$. Since $F(L)\cong F(M)$, the object $GF(M)\cong GF(L)$ belongs to $\mcY$ by \autoref{LocAndPrelocSubcat}. Let $\eta\colon 1_{\mcA}\to GF$ be the unit morphism. There is an exact sequence
	\begin{equation*}
		0\to\Ker\eta_{M}\to M\to\Im\eta_{M}\to 0.
	\end{equation*}
	By \autoref{PropertiesOfQuotCat} \autoref{DescriptionOfUnit}, the object $\Ker\eta_{M}$ belongs to $\mcX$. The subobject $\Im\eta_{M}$ of $GF(M)$ belongs to $\mcY$. Therefore $M$ belongs to $\mcX*\mcY$.
	
	\autoref{ImageOfExtOfPrelocSubcats} By \autoref{QuotOfPrelocSubcat} \autoref{QuotOfExtOfSubcatsIsExtOfQuotOfSubcats},
	\begin{equation*}
		F(\mcY_{1}*\mcY_{2})\subset F(\mcY_{1}*\mcX*\mcY_{2})=F(\mcY_{1})*F(\mcY_{2}).
	\end{equation*}
	By \autoref{CommutativityBetweenLocSubcatAndPrelocSubcat},
	\begin{equation*}
		F(\mcY_{1}*\mcX*\mcY_{2})\subset F(\mcX*\mcY_{1}*\mcY_{2})\subset F(\mcX)*F(\mcY_{1}*\mcY_{2})=F(\mcY_{1}*\mcY_{2}).
	\end{equation*}
	
	\autoref{ImageOfLocSubcat} By \autoref{ImageOfExtOfPrelocSubcats},
	\begin{equation*}
		F(\mcY)*F(\mcY)=F(\mcY*\mcY)=F(\mcY).\qedhere
	\end{equation*}
\end{proof}

\begin{Theorem}\label{LocOfExtOfPrelocSubcats}
	Assume that the following elements correspond to each other by the bijections in \autoref{LocOfPrelocSubcat} for each $i=1,2$.
	\begin{enumerate}
		\item $\mcY^{i}$.
		\item $\set{\mcY^{i}_{\lambda}}_{\lambda\in\Lambda}$.
		\item $\set{\mcY^{i}(\alpha)}_{\alpha\in\ASpec\mcA}$.
	\end{enumerate}
	Then the following elements correspond to each other by the bijections.
	\begin{enumerate}
		\item $\mcY^{1}*\mcY^{2}$.
		\item $\set{\mcY^{1}_{\lambda}*\mcY^{2}_{\lambda}}_{\lambda\in\Lambda}$.
		\item $\set{\mcY^{1}(\alpha)*\mcY^{2}(\alpha)}_{\alpha\in\ASpec\mcA}$.
	\end{enumerate}
\end{Theorem}

\begin{proof}
	This follows from \autoref{CommutativityAndImageOfExtOfPrelocSubcats}.
\end{proof}

\begin{Corollary}\label{LocOfLocSubcats}
	The bijections in \autoref{LocOfPrelocSubcat} restrict to bijections between following sets.
	\begin{enumerate}
		\item\label{LocOfLocSubcats:LocSubcat} The set of localizing subcategories of $\mcA$.
		\item\label{LocOfLocSubcats:LocSubcatsOfQuotCats} The set of families $\set{\mcY_{\lambda}\subset\mcU_{\lambda}}_{\lambda\in\Lambda}$ of localizing subcategories such that $F^{\lambda}_{\mu}(\mcY_{\lambda})=\mcY_{\mu}$ for each $\lambda,\mu\in\Lambda$ with $\ASpec\,\mcU_{\mu}\subset\ASpec\,\mcU_{\lambda}$.
		\item\label{LocOfLocSubcats:LocSubcatsOfLocalizedCats} The set of families $\set{\mcY(\alpha)\subset\mcA_{\alpha}}_{\alpha\in\ASpec\mcA}$ of localizing subcategories such that $\mcY(\beta)_{\alpha}=\mcY(\alpha)$ for each $\alpha,\beta\in\ASpec\mcA$ with $\alpha\leq\beta$.
	\end{enumerate}
\end{Corollary}

\begin{proof}
	This follows from \autoref{LocOfExtOfPrelocSubcats}.
\end{proof}

Prime localizing subcategories of $\mcA$ are characterized as follows.

\begin{Theorem}\label{CharacterizationOfPrimeLocSubcat}
	Assume \autoref{SettingOfLocOfPrelocSubcat}, and let $\mcX$ be a localizing subcategory of $\mcA$. Then the following assertions are equivalent.
	\begin{enumerate}
		\item\label{CharacterizationOfPrimeLocSubcat:PrimeLocSubcat} $\mcX$ is a prime localizing subcategory of $\mcA$.
		\item\label{CharacterizationOfPrimeLocSubcat:LocSubcatAssociatedToAtom} There exists $\alpha\in\ASpec\mcA$ such that $\mcX=\mcX(\alpha)$.
		\item\label{CharacterizationOfPrimeLocSubcat:CompletelyMeetIrredLocSubcat} For each family $\{\mcX_{\omega}\}_{\omega\in\Omega}$ of localizing subcategories of $\mcA$ satisfying $\mcX=\bigcap_{\omega\in\Omega}\mcX_{\omega}$, there exists $\omega\in\Omega$ such that $\mcX=\mcX_{\omega}$.
		\item\label{CharacterizationOfPrimeLocSubcat:CompletelyMeetIrredLocSubcatInPrelocSubcats} For each family $\{\mcY_{\omega}\}_{\omega\in\Omega}$ of prelocalizing subcategories of $\mcA$ satisfying $\mcX=\bigcap_{\omega\in\Omega}\mcY_{\omega}$, there exists $\omega\in\Omega$ such that $\mcX=\mcY_{\omega}$.
	\end{enumerate}
\end{Theorem}

\begin{proof}
	The equivalence \autoref{CharacterizationOfPrimeLocSubcat:PrimeLocSubcat}$\Leftrightarrow$\autoref{CharacterizationOfPrimeLocSubcat:LocSubcatAssociatedToAtom} follows from \autoref{BijectionBetweenAtomsAndPrimeLocSubcat}.
	
	Let $\{\mcY_{\omega}\}_{\omega\in\Omega}$ be a family of prelocalizing subcategories of $\mcA$ satisfying $\mcX(\alpha)=\bigcap_{\omega\in\Omega}\mcY_{\omega}$. Since $H(\alpha)$ does not belong to $\mcX(\alpha)$, there exists $\omega\in\Omega$ such that $H(\alpha)$ does not belong to $\mcY_{\omega}$. By \autoref{AtomAndAtomicObjBelongToPrelocSubcat}, $\alpha\notin\ASupp\mcY_{\omega}$, and hence $\mcY_{\omega}\subset\mcX(\alpha)$. This shows \autoref{CharacterizationOfPrimeLocSubcat:LocSubcatAssociatedToAtom}$\Rightarrow$\autoref{CharacterizationOfPrimeLocSubcat:CompletelyMeetIrredLocSubcatInPrelocSubcats}.
	
	The implication \autoref{CharacterizationOfPrimeLocSubcat:CompletelyMeetIrredLocSubcatInPrelocSubcats}$\Rightarrow$\autoref{CharacterizationOfPrimeLocSubcat:CompletelyMeetIrredLocSubcat} is obvious. The implication \autoref{CharacterizationOfPrimeLocSubcat:CompletelyMeetIrredLocSubcat}$\Rightarrow$\autoref{CharacterizationOfPrimeLocSubcat:LocSubcatAssociatedToAtom} follows from \autoref{LocSubcatIsIntersectionOfPrimeLocSubcat}.
\end{proof}

\section{Classification of prelocalizing subcategories}
\label{sec:ClassificationOfPrelocSubcats}

Let $X$ be a locally noetherian scheme with the structure sheaf $\OX$. In this section, we classify the prelocalizing subcategories of $\QCoh X$. Let $\set{U_{\lambda}}_{\lambda\in\Lambda}$ be an open affine basis of $X$. Let $i_{\lambda}\colon U_{\lambda}\into X$ be the immersion for each $\lambda\in\Lambda$, and let $i_{\lambda,\mu}\colon U_{\mu}\into U_{\lambda}$ be the immersion for each $\lambda,\mu\in\Lambda$ with $U_{\mu}\subset U_{\lambda}$.

We recall the notion of a filter. This is an essential tool to classify prelocalizing subcategories.

\begin{Definition}\label{FiltAndPrincipalFiltOfSubobjs}
	Let $\mcA$ be a Grothendieck category, and let $M$ be an object in $\mcA$.
	\begin{enumerate}
		\item\label{FiltOfSubobjs} A \emph{filter} of subobjects of $M$ in $\mcA$ is a set $\mcF$ of subobjects of $M$ satisfying the following conditions.
		\begin{enumerate}
			\item\label{FiltOfSubobjs:LargestElement} $M\in\mcF$.
			\item\label{FiltOfSubobjs:UpwardClosedness} If $L\subset L'$ are subobjects of $M$ with $L\in\mcF$, then $L'\in\mcF$.
			\item\label{FiltOfSubobjs:FiniteIntersection} If $L_{1},L_{2}\in\mcF$, then $L_{1}\cap L_{2}\in\mcF$.
		\end{enumerate}
		If there is no danger of confusion, we simply say that $\mcF$ is a filter of $M$.
		\item\label{PrincipalFiltOfSubobjs} For each subobject $L$ of $M$, denote by $\mcF(L)$ the filter consisting of all subobjects $L'$ of $M$ with $L\subset L'$. A filter of the form $\mcF(L)$ is called a \emph{principal filter}.
	\end{enumerate}
\end{Definition}

\begin{Remark}\label{CharacterizationOfPrincipalFilt}
	In \autoref{FiltAndPrincipalFiltOfSubobjs} \autoref{PrincipalFiltOfSubobjs}, the principal filter $\mcF(L)$ is closed under arbitrary intersection. Conversely, if a filter $\mcF$ of $M$ is closed under arbitrary intersection, then $\mcF=\mcF(L)$, where $L$ is the smallest element of $\mcF$.
	
	It is obvious that the map
	\begin{equation*}
		\setwithspace{\text{subobjects of }M}\to\setwithspace{\text{principal filters of }M}
	\end{equation*}
	given by $L\mapsto\mcF(L)$ is bijective.
\end{Remark}

For a ring $\Lambda$, we say that a filter $\mcF$ of right ideals of $\Lambda$ is \emph{prelocalizing} if for each $L\in\mcF$ and $a\in\Lambda$, the right ideal
\begin{equation*}
	a^{-1}L=\conditionalset{b\in\Lambda}{ab\in L}
\end{equation*}
of $\Lambda$ belongs to $\mcF$. For a ring $\Lambda$, Gabriel \cite{Gabriel} gave a classification of the prelocalizing subcategories of $\Mod\Lambda$.

\begin{Theorem}[{\cite[Lemma~V.2.1]{Gabriel}}]\label{BijectionBetweenPrelocSubcatsAndFiltsOfSubobjOfRing}
	Let $\Lambda$ be a ring. There is a bijection
	\begin{equation*}
		\setwithspace{\text{prelocalizing subcategories of }\Mod\Lambda}\to\setwithspace{\text{prelocalizing filters of right ideals of }\Lambda}
	\end{equation*}
	given by
	\begin{equation*}
		\mcY\mapsto\conditionalset[\bigg]{L\subset\Lambda\text{ in }\Mod\Lambda}{\frac{\Lambda}{L}\in\mcY}.
	\end{equation*}
	Its inverse is given by
	\begin{align*}
		\mcF&\mapsto\conditionalset{M\in\Mod\Lambda}{\Ann_{\Lambda}(x)\in\mcF\text{ for every }x\in M}\\
		&=\conditionalgeneratedset[\bigg]{\!\!preloc}{\frac{\Lambda}{L}\in\Mod\Lambda}{L\in\mcF}.
	\end{align*}
\end{Theorem}

\begin{proof}
	\cite[Theorem~4.9.1]{Popescu}.
\end{proof}

For a commutative ring $R$, every filter $\mcF$ of $R$ is prelocalizing. Indeed, for $L\in\mcF$ and $a\in R$, we have $L\subset a^{-1}L$, and hence $a^{-1}L\in\mcF$. Therefore the following assertion holds.

\begin{Corollary}\label{BijectionBetweenPrelocSubcatsAndFiltsOfSubobjOfCommRing}
	Let $R$ be a commutative ring. There is a bijection
	\begin{equation*}
		\setwithspace{\text{prelocalizing subcategories of }\Mod R}\to\setwithspace{\text{filters of ideals of }R}
	\end{equation*}
	given by
	\begin{equation*}
		\mcY\mapsto\conditionalset[\bigg]{I\subset R\text{ in }\Mod R}{\frac{R}{I}\in\mcY}.
	\end{equation*}
	Its inverse is given by
	\begin{align*}
		\mcF&\mapsto\conditionalset{M\in\Mod R}{\Ann_{R}(x)\in\mcF\text{ for every }x\in M}\\
		&=\conditionalgeneratedset[\bigg]{\!\!preloc}{\frac{R}{I}\in\Mod R}{I\in\mcF}.
	\end{align*}
\end{Corollary}

\begin{proof}
	This is immediate from \autoref{BijectionBetweenPrelocSubcatsAndFiltsOfSubobjOfRing}.
\end{proof}

In the case of a locally noetherian scheme $X$, we need to use the notion of a local filter instead of a filter (see \autoref{ClassficationOfPrelocSubcatsForLocNoethSch} and \autoref{ExOfCoprodOfSch}).

\begin{Definition}\label{LocFiltOfSubobjs}
	Let $X$ be a locally noetherian scheme. We say that a filter $\mcF$ of subobjects of $\OX$ in $\QCoh X$ is a \emph{local filter} of $\OX$ if it satisfies the following condition: let $I$ be a subobject of $\OX$, and assume that for each $x\in X$, there exist an open affine neighborhood $U$ of $x$ in $X$ and $I'\in\mcF$ such that $I'|_{U}\subset I|_{U}$ as a subobject of $\mcO_{U}$. Then $I\in\mcF$.
\end{Definition}

\begin{Proposition}\label{PrincipalFiltIsLocFilt}
	Every principal filter of $\OX$ is a local filter.
\end{Proposition}

\begin{proof}
	For every subobject $I$ of $\OX$, we show that $\mcF(I)$ is a local filter. Let $I'$ be a subobject of $\OX$ such that for each $x\in X$, there exist an open affine neighborhood $U(x)$ of $x$ in $X$ and $J(x)\in\mcF(I)$ such that $J(x)|_{U(x)}\subset I'|_{U(x)}$. Let $J:=\bigcap_{x\in X}J(x)$ in $\QCoh X$. Then $J\in\mcF(I)$, and $J|_{U(x)}\subset I'|_{U(x)}$ for each $x\in X$. For each open subset $U$ of $X$,
	\begin{align*}
		J(U)&=\conditionalset{s\in\OX(U)}{s|_{U(x)\cap U}\in J(U(x)\cap U)\text{ for each }x\in X}\\
		&\subset\conditionalset{s\in\OX(U)}{s|_{U(x)\cap U}\in I'(U(x)\cap U)\text{ for each }x\in X}\\
		&=I'(U),
	\end{align*}
	and hence $J\subset I'$ follows. This implies that $I'\in\mcF$.
\end{proof}

The next result shows that the local filters of $\OX$ are exactly the same as the filters of $\OX$ in the case where $X$ is quasi-compact. This is the reason that we do not need to consider a local filter in the case of a commutative ring.

\begin{Proposition}\label{FiltOfSubobjsForNoethSchIsLocFilt}
	If $X$ is a noetherian scheme, then every filter of $\OX$ is a local filter.
\end{Proposition}

\begin{proof}
	Let $\mcF$ be a filter of $\OX$. Let $I$ be a subobject of $\OX$, and assume that for each $x\in X$, there exist an open affine neighborhood $U(x)$ of $x$ in $X$ and $I'(x)\in\mcF$ such that $I'(x)|_{U(x)}\subset I|_{U(x)}$. Since $X$ is quasi-compact, there exists $x_{1},\ldots,x_{n}\in X$ such that $X=\bigcup_{j=1}^{n}U(x_{j})$. Let $I':=\bigcap_{j=1}^{n}I'(x_{j})$. Then $I'$ belongs to $\mcF$. Since $I'|_{U(x_{j})}\subset I|_{U(x_{j})}$ for each $j=1,\ldots,n$, we have $I'\subset I$, and hence $I$ also belongs to $\mcF$.
\end{proof}

The following result describes the local filter generated by a set of subobjects of $\OX$.

\begin{Proposition}\label{LocFiltGeneratedBySet}
	Let $\mcS$ be a set of subobjects of $\OX$. Let $\mcF$ be the set consisting of all subobjects $I$ of $\OX$ satisfying the following condition: for each $x\in X$, there exist an open affine neighborhood $U$ of $x$ in $X$ and $n\in\mbZ_{\geq 1}$ and $I_{1},\ldots,I_{n}\in\mcS$ such that
	\begin{equation*}
		(I_{1}\cap\cdots\cap I_{n})|_{U}\subset I|_{U}.
	\end{equation*}
	Then $\mcF$ is the smallest local filter of $\OX$ including $\mcS$.
\end{Proposition}

\begin{proof}
	It is obvious that $\mcF$ satisfies the conditions \autoref{FiltOfSubobjs:LargestElement} and \autoref{FiltOfSubobjs:UpwardClosedness} in \autoref{FiltAndPrincipalFiltOfSubobjs} \autoref{FiltOfSubobjs}. We show that \autoref{FiltOfSubobjs:FiniteIntersection} is satisfied. Let $I^{(1)},I^{(2)}\in\mcF$. Then for each $j=1,2$ and $x\in X$, there exist an open affine neighborhood $U^{(j)}$ of $x$ in $X$ and $n_{j}\in\mbZ_{\geq 1}$ and $I^{(j)}_{1},\ldots,I^{(j)}_{n_{j}}\in\mcS$ such that
	\begin{equation*}
		(I^{(j)}_{1}\cap\cdots\cap I^{(j)}_{n_{j}})|_{U^{(j)}}\subset I^{(j)}|_{U^{(j)}}.
	\end{equation*}
	Then
	\begin{equation*}
		(I^{(1)}_{1}\cap\cdots\cap I^{(1)}_{n_{1}}\cap I^{(2)}_{1}\cap\cdots\cap I^{(2)}_{n_{2}})|_{U^{(1)}\cap U^{(2)}}\subset (I^{(1)}\cap I^{(2)})|_{U^{(1)}\cap U^{(2)}}.
	\end{equation*}
	This shows $I^{(1)}\cap I^{(2)}\in\mcF$. Hence $\mcF$ is a filter of $\OX$.
	
	Let $I$ be a subobject of $\OX$ such that for each $x\in X$, there exist an open affine neighborhood $U$ of $x$ in $X$ and $I'\in\mcF$ such that $I'|_{U}\subset I|_{U}$. Let $x\in X$, and take such $U$ and $I'$. Then there exists an open affine neighborhood $U'$ of $x$ in $X$ and $n\in\mbZ_{\geq 1}$ and $I'_{1},\ldots,I'_{n}\in\mcS$ such that
	\begin{equation*}
		(I'_{1}\cap\cdots\cap I'_{n})|_{U'}\subset I'|_{U'}.
	\end{equation*}
	Since
	\begin{equation*}
		(I'_{1}\cap\cdots\cap I'_{n})|_{U\cap U'}\subset I'|_{U\cap U'}\subset I|_{U\cap U'},
	\end{equation*}
	it holds that $I\in\mcF$. This shows that $\mcF$ is a local filter. It is obvious that $\mcF$ is the smallest local filter of $\OX$ including $\mcS$.
\end{proof}

In the setting of \autoref{LocFiltGeneratedBySet}, the local filter $\mcF$ is denoted by $\locfilt{\mcS}$.

We investigate the restriction of a filter to an open affine subscheme and the localization at a point.

\begin{Proposition}\label{ImageOfFiltIsFilt}
	Let $\mcF$ be a filter of $\OX$.
	\begin{enumerate}
		\item\label{ImageOfFiltToOpenAffineSubschIsFilt} For every $\lambda\in\Lambda$, the set
		\begin{equation*}
			\mcF|_{U_{\lambda}}:=\conditionalset{I|_{U_{\lambda}}\subset\mcO_{U_{\lambda}}\text{ in }\QCoh U_{\lambda}}{I\in\mcF}
		\end{equation*}
		is a filter of $\mcO_{U_{\lambda}}$.
		\item\label{ImageOfFiltToLocalRingIsFilt} For every $x\in X$, the set
		\begin{equation*}
			\mcF_{x}:=\conditionalset{I_{x}\subset\OXx\text{ in }\Mod\OXx}{I\in\mcF}
		\end{equation*}
		is a filter of $\OXx$.
	\end{enumerate}
\end{Proposition}

\begin{proof}
	\autoref{ImageOfFiltToOpenAffineSubschIsFilt} Since $\OX\in\mcF$, we have $\mcO_{U_{\lambda}}\in\mcF|_{U_{\lambda}}$.
	
	Let $\widetilde{I}\subset\widetilde{I}'$ be subobjects of $\mcO_{U_{\lambda}}$ with $\widetilde{I}\in\mcF|_{U_{\lambda}}$. By \autoref{SubobjOfObjInQuotCat}, there exists a largest subobject $I$ (resp.\ $I'$) of $\OX$ satisfying $I|_{U_{\lambda}}\subset\widetilde{I}$ (resp.\ $I'|_{U_{\lambda}}\subset\widetilde{I}'$), and it holds that $I|_{U_{\lambda}}=\widetilde{I}$ (resp.\ $I'|_{U_{\lambda}}=\widetilde{I}'$). Then $I\in\mcF$ and $I\subset I'$ imply $I'\in\mcF$. We deduce that $\widetilde{I}'=I'|_{U_{\lambda}}\in\mcF|_{U_{\lambda}}$.
	
	Let $\widetilde{I}_{1},\widetilde{I}_{2}\in\mcF|_{U_{\lambda}}$. Then for each $i=1,2$, there exists $I_{i}\in\mcF$ such that $I_{i}|_{U_{\lambda}}=\widetilde{I}_{i}$. It holds that $I_{1}\cap I_{2}\in\mcF$. Since $(-)|_{U_{\lambda}}\colon\QCoh X\to\QCoh U_{\lambda}$ is an exact functor, $I_{1}|_{U_{\lambda}}\cap I_{2}|_{U_{\lambda}}=(I_{1}\cap I_{2})|_{U_{\lambda}}\in\mcF|_{U_{\lambda}}$.
	
	\autoref{ImageOfFiltToLocalRingIsFilt} This is shown similarly to \autoref{ImageOfFiltToOpenAffineSubschIsFilt}.
\end{proof}

We give a characterization of a local filter.

\begin{Proposition}\label{CharacterizationOfLocFilt}
	Let $\mcF$ be a filter of $\OX$. Then the following assertions are equivalent.
	\begin{enumerate}
		\item\label{CharacterizationOfLocFilt:LocFilt} $\mcF$ is a local filter.
		\item\label{CharacterizationOfLocFilt:CharacterizationOfLocFilt} Let $I$ be a subobject of $\OX$ such that for each $x\in X$, there exists an open affine neighborhood $U$ of $x$ in $X$ satisfying $I|_{U}\in\mcF|_{U}$. Then $I\in\mcF$.
	\end{enumerate}
\end{Proposition}

\begin{proof}
	It is obvious that \autoref{CharacterizationOfLocFilt:LocFilt} implies \autoref{CharacterizationOfLocFilt:CharacterizationOfLocFilt}.
	
	Assume \autoref{CharacterizationOfLocFilt:CharacterizationOfLocFilt}. Let $I$ be a subobject of $\OX$ such that for each $x\in X$, there exist an open affine neighborhood $U$ of $x$ in $X$ and $I'\in\mcF$ satisfying $I'|_{U}\subset I|_{U}$. Since $\mcF|_{U}$ is a filter of $\mcO_{U}$ by \autoref{ImageOfFiltIsFilt} \autoref{ImageOfFiltToOpenAffineSubschIsFilt}, we have $I|_{U}\in\mcF|_{U}$. Hence $I\in\mcF$. This shows \autoref{CharacterizationOfLocFilt:LocFilt}.
\end{proof}

The following lemmas show that the bijection in \autoref{BijectionBetweenPrelocSubcatsAndFiltsOfSubobjOfCommRing} commutes with the restriction to an open affine subscheme and the localization at a point.

\begin{Lemma}\label{PrelocSubcatCorrespondingToImageOfFilt}
	Let $\lambda,\mu\in\Lambda$ with $U_{\mu}\subset U_{\lambda}$. Let $\mcY$ be a prelocalizing subcategory of $\QCoh U_{\lambda}$, and let $\mcF$ be the corresponding filter of $\mcO_{U_{\lambda}}$ by the bijection in \autoref{BijectionBetweenPrelocSubcatsAndFiltsOfSubobjOfCommRing}. Then the filter $\mcF|_{U_{\mu}}$ of $\mcO_{U_{\mu}}$ corresponds to the prelocalizing subcategory $\mcY|_{U_{\mu}}$ of $\QCoh U_{\mu}$ by the bijection.
\end{Lemma}

\begin{proof}
	Let $\mcF'$ be the filter of $\mcO_{U_{\mu}}$ corresponding to $\mcY|_{U_{\mu}}$, that is,
	\begin{equation*}
		\mcF'=\conditionalset[\bigg]{\widetilde{I}\subset\mcO_{U_{\mu}}\text{ in }\QCoh U_{\mu}}{\frac{\mcO_{U_{\mu}}}{\widetilde{I}}\in\mcY|_{U_{\mu}}}.
	\end{equation*}
	It is obvious that $\mcF|_{U_{\mu}}\subset\mcF'$. Let $\widetilde{I}\in\mcF'$. Then there exists an object $M$ in $\QCoh X$ which belongs to $\mcY$ such that $\mcO_{U_{\mu}}/\widetilde{I}\cong M|_{U_{\mu}}$. By \autoref{SubobjOfObjInQuotCat}, there exists a subobject $I$ of $\mcO_{U_{\lambda}}$ such that $I|_{U_{\mu}}=\widetilde{I}$, and $\mcO_{U_{\lambda}}/I$ is $\mcX$-torsionfree, where
	\begin{equation*}
		\mcX=\conditionalset{M'\in\QCoh U_{\lambda}}{M'|_{U_{\mu}}=0}.
	\end{equation*}
	By \autoref{PropertiesOfQuotCat} \autoref{DescriptionOfUnit}, the canonical morphism $\mcO_{U_{\lambda}}/I\to(i_{\lambda,\mu})_{*}i_{\lambda,\mu}^{*}(\mcO_{U_{\lambda}}/I)$ is a monomorphism. By \autoref{LocAndPrelocSubcat}, the object
	\begin{equation*}
		(i_{\lambda,\mu})_{*}i_{\lambda,\mu}^{*}\left(\frac{\mcO_{U_{\lambda}}}{I}\right)\cong(i_{\lambda,\mu})_{*}\left(\frac{\mcO_{U_{\mu}}}{\widetilde{I}}\right)\cong(i_{\lambda,\mu})_{*}i_{\lambda,\mu}^{*}M
	\end{equation*}
	belongs to $\mcY$. Hence $\mcO_{U_{\lambda}}/I$ also belongs to $\mcY$. This shows that $I\in\mcF$ and that $\widetilde{I}=I|_{U_{\mu}}\in\mcF|_{U_{\mu}}$. Therefore $\mcF|_{U_{\mu}}=\mcF'$.
\end{proof}

\begin{Lemma}\label{PrelocSubcatCorrespondingToLocOfFilt}
	Let $x,y\in X$ with $y\in\overline{\set{x}}$. Let $\mcY$ be a prelocalizing subcategory of $\Mod\OXy$, and let $\mcF$ be the corresponding filter of $\OXy$ by the bijection in \autoref{BijectionBetweenPrelocSubcatsAndFiltsOfSubobjOfCommRing}. Then the filter $\mcF_{x}$ of $\OXx$ corresponds to the prelocalizing subcategory $\mcY_{x}$ of $\Mod\OXx$ by the bijection.
\end{Lemma}

\begin{proof}
	This is shown similarly to \autoref{PrelocSubcatCorrespondingToImageOfFilt}.
\end{proof}

We show a lemma to glue filters on open affine basis to a local filter of $\OX$.

\begin{Lemma}\label{GluingOfFiltsAndDescriptionOfLocFilt}\leavevmode
	\begin{enumerate}
		\item\label{DescriptionOfLocFilt} For every local filter $\mcF$ of $\OX$,
		\begin{equation*}
			\mcF=\conditionalset{I\subset\OX\text{ in }\QCoh X}{I|_{U_{\lambda}}\in\mcF|_{U_{\lambda}}\text{ for each }\lambda\in\Lambda}.
		\end{equation*}
		\item\label{GluingOfFilts} Let $\mcF_{\lambda}$ be a filter of $\mcO_{U_{\lambda}}$ for each $\lambda\in\Lambda$, and assume that $\mcF_{\lambda}|_{U_{\mu}}=\mcF_{\mu}$ for each $\lambda,\mu\in\Lambda$ with $U_{\mu}\subset U_{\lambda}$. Then there exists a unique local filter $\mcF$ of $\OX$ satisfying $\mcF|_{U_{\lambda}}=\mcF_{\lambda}$ for each $\lambda\in\Lambda$.
	\end{enumerate}
\end{Lemma}

\begin{proof}
	\autoref{DescriptionOfLocFilt} This follows from \autoref{CharacterizationOfLocFilt}.
	
	\autoref{GluingOfFilts} The uniqueness follows from \autoref{DescriptionOfLocFilt}. Let
	\begin{equation*}
		\mcF:=\conditionalset{I\subset\OX\text{ in }\QCoh X}{I|_{U_{\lambda}}\in\mcF_{\lambda}\text{ for each }\lambda\in\Lambda}.
	\end{equation*}
	It is straightforward to show that $\mcF$ is a filter of $\OX$ satisfying $\mcF|_{U_{\lambda}}\subset\mcF_{\lambda}$ for each $\lambda\in\Lambda$.
	
	Let $I$ be a subobject of $\OX$ such that for each $x\in X$, there exists an open affine neighborhood $U$ of $x$ in $X$ satisfying $I|_{U}\in\mcF|_{U}$. For each $\lambda\in\Lambda$ and $y\in U_{\lambda}$, there exists an open affine neighborhood $U'$ of $y$ in $X$ satisfying $I|_{U'}\in\mcF|_{U'}$. Take $\mu\in\Lambda$ satisfying $y\in U_{\mu}\subset U_{\lambda}\cap U'$. Then $(I|_{U_{\lambda}})|_{U_{\mu}}=(I|_{U'})|_{U_{\mu}}\in(\mcF|_{U'})|_{U_{\mu}}=(\mcF|_{U_{\lambda}})|_{U_{\mu}}$. Since $\mcF|_{U_{\lambda}}$ is a local filter by \autoref{ImageOfFiltIsFilt} \autoref{ImageOfFiltToOpenAffineSubschIsFilt} and \autoref{FiltOfSubobjsForNoethSchIsLocFilt}, we have $I|_{U_{\lambda}}\in\mcF|_{U_{\lambda}}\subset\mcF_{\lambda}$. This shows that $I\in\mcF$. By \autoref{CharacterizationOfLocFilt}, the filter $\mcF$ is a local filter.
	
	We show that $\mcF_{\lambda}\subset\mcF|_{U_{\lambda}}$. Let $\widetilde{J}\in\mcF_{\lambda}$. By \autoref{SubobjOfObjInQuotCat}, there exists a subobject $J$ of $\OX$ such that $J|_{U_{\lambda}}=\widetilde{J}$, and $\OX/J$ is $\mcX_{U_{\lambda}}$-torsionfree (see \autoref{OpenSubschIsQuotCat}). It suffices to show that $J\in\mcF$, that is, $J|_{U_{\mu}}\in\mcF_{\mu}$ for each $\mu\in\Lambda$. Denote by $\mcY_{\lambda}$ and $\mcY_{\mu}$ the prelocalizing subcategories of $\QCoh U_{\lambda}$ and $\QCoh U_{\mu}$ corresponding to $\mcF_{\lambda}$ and $\mcF_{\mu}$ by \autoref{BijectionBetweenPrelocSubcatsAndFiltsOfSubobjOfCommRing}, respectively. We show that the object $\mcO_{U_{\mu}}/J|_{U_{\mu}}$ belongs to $\mcY_{\mu}$. Let $x\in\Ass_{U_{\mu}}(\mcO_{U_{\mu}}/J|_{U_{\mu}})$. By \autoref{AAssOfObjInQuotCatInSetting},
	\begin{equation*}
		x\in\Ass_{U_{\mu}}\frac{\OX}{J}\bigg|_{U_{\mu}}=\Ass_{X}\frac{\OX}{J}\cap U_{\mu}\subset U_{\lambda}\cap U_{\mu}.
	\end{equation*}
	Hence
	\begin{equation*}
		\left(\frac{\mcO_{U_{\mu}}}{J|_{U_{\mu}}}\right)_{x}=\frac{\OXx}{J_{x}}
		=\left(\frac{\mcO_{U_{\lambda}}}{J|_{U_{\lambda}}}\right)_{x}
		=\left(\frac{\mcO_{U_{\lambda}}}{\widetilde{J}}\right)_{x}
		\in(\mcY_{\lambda})_{x}.
	\end{equation*}
	Take $\nu\in\Lambda$ such that $x\in U_{\nu}\subset U_{\lambda}\cap U_{\mu}$. Then $(\mcY_{\lambda})_{x}=(\mcY_{\lambda}|_{U_{\nu}})_{x}=(\mcY_{\nu})_{x}=(\mcY_{\mu}|_{U_{\nu}})_{x}=(\mcY_{\mu})_{x}$. Hence by \autoref{LocOfObjAtAAssAndPrelocSubcat}, the object $\mcO_{U_{\mu}}/J|_{U_{\mu}}$ belongs to $\mcY_{\mu}$. This shows that $J|_{U_{\mu}}\in\mcF_{\mu}$.
\end{proof}

The following theorem is the main result in this section, which gives a classification of the prelocalizing subcategory of $\QCoh X$.

\begin{Theorem}\label{ClassficationOfPrelocSubcatsForLocNoethSch}
	Let $X$ be a locally noetherian scheme, and let $\set{U_{\lambda}}_{\lambda\in\Lambda}$ be an open affine basis of $X$. Then there exist bijections between the following sets.
	\begin{enumerate}
		\item\label{ClassficationOfPrelocSubcatsForLocNoethSch:PrelocSubcat} The set of prelocalizing subcategories of $\QCoh X$.
		\item\label{ClassficationOfPrelocSubcatsForLocNoethSch:PrelocSubcatsOnOpenAffineSubschs} The set of families $\set{\mcY_{\lambda}\subset\QCoh U_{\lambda}}_{\lambda\in\Lambda}$ of prelocalizing subcategories such that $\mcY_{\lambda}|_{U_{\mu}}=\mcY_{\mu}$ for each $\lambda,\mu\in\Lambda$ with $U_{\mu}\subset U_{\lambda}$.
		\item\label{ClassficationOfPrelocSubcatsForLocNoethSch:PrelocSubcatsOnLocRings} The set of families $\set{\mcY(x)\subset\Mod\OXx}_{x\in X}$ of prelocalizing subcategories such that $\mcY(y)_{x}=\mcY(x)$ for each $x,y\in X$ with $y\in\overline{\set{x}}$.
		\item\label{ClassficationOfPrelocSubcatsForLocNoethSch:LocFilt} The set of local filters of $\OX$.
		\item\label{ClassficationOfPrelocSubcatsForLocNoethSch:FiltsOnOpenAffineSubschs} The set of families $\set{\mcF_{\lambda}}_{\lambda\in\Lambda}$, where $\mcF_{\lambda}$ is a filter of $\mcO_{U_{\lambda}}$ for each $\lambda\in\Lambda$, such that $\mcF_{\lambda}|_{U_{\mu}}=\mcF_{\mu}$ for each $\lambda,\mu\in\Lambda$ with $U_{\mu}\subset U_{\lambda}$.
		\item\label{ClassficationOfPrelocSubcatsForLocNoethSch:FiltsOnLocRings} The set of families $\set{\mcF(x)}_{x\in X}$, where $\mcF(x)$ is a filter of $\OXx$ for each $x\in X$, such that $\mcF(y)_{x}=\mcF(x)$ for each $x,y\in X$ with $y\in\overline{\set{x}}$.
	\end{enumerate}
	The correspondences are given as follows.
	\begin{align*}
		&\autoref{ClassficationOfPrelocSubcatsForLocNoethSch:PrelocSubcat}\ \mcY\mapsto
		\begin{cases}
			\autoref{ClassficationOfPrelocSubcatsForLocNoethSch:LocFilt}\ \displaystyle\conditionalset[\bigg]{I\subset\OX\text{ in }\QCoh X}{\frac{\OX}{I}\in\mcY}\\
			\autoref{ClassficationOfPrelocSubcatsForLocNoethSch:PrelocSubcatsOnOpenAffineSubschs}\ \set{\mcY|_{U_{\lambda}}}_{\lambda\in\Lambda}\\
			\autoref{ClassficationOfPrelocSubcatsForLocNoethSch:PrelocSubcatsOnLocRings}\ \set{\mcY_{x}}_{x\in X}
		\end{cases}\\
		&\autoref{ClassficationOfPrelocSubcatsForLocNoethSch:LocFilt}\ \mcF\mapsto
		\begin{cases}
			\autoref{ClassficationOfPrelocSubcatsForLocNoethSch:PrelocSubcat}\ \displaystyle\conditionalgeneratedset[\bigg]{\!\!preloc}{\frac{\OX}{I}}{I\in\mcF}.\\
			\autoref{ClassficationOfPrelocSubcatsForLocNoethSch:FiltsOnOpenAffineSubschs}\ \set{\mcF|_{U_{\lambda}}}_{\lambda\in\Lambda}\\
			\autoref{ClassficationOfPrelocSubcatsForLocNoethSch:FiltsOnLocRings}\ \set{\mcF_{x}}_{x\in X}
		\end{cases}\\
		&\autoref{ClassficationOfPrelocSubcatsForLocNoethSch:PrelocSubcatsOnOpenAffineSubschs}\ \set{\mcY_{\lambda}}_{\lambda\in\Lambda}\mapsto\autoref{ClassficationOfPrelocSubcatsForLocNoethSch:PrelocSubcat}\ \conditionalset{M\in\QCoh X}{M|_{U_{\lambda}}\in\mcY_{\lambda}\text{ for each }\lambda\in\Lambda}\\
		&\autoref{ClassficationOfPrelocSubcatsForLocNoethSch:PrelocSubcatsOnLocRings}\ \set{\mcY(x)}_{x\in X}\mapsto\autoref{ClassficationOfPrelocSubcatsForLocNoethSch:PrelocSubcat}\ \conditionalset{M\in\QCoh X}{M_{x}\in\mcY(x)\text{ for each }x\in X}\\
		&\autoref{ClassficationOfPrelocSubcatsForLocNoethSch:FiltsOnOpenAffineSubschs}\ \set{\mcF_{\lambda}}_{\lambda\in\Lambda}\mapsto\autoref{ClassficationOfPrelocSubcatsForLocNoethSch:LocFilt}\ \conditionalset{I\subset\OX\text{ in }\QCoh X}{I|_{U_{\lambda}}\in\mcF_{\lambda}\text{ for each }\lambda\in\Lambda}\\
		&\autoref{ClassficationOfPrelocSubcatsForLocNoethSch:FiltsOnLocRings}\ \set{\mcF(x)}_{x\in X}\mapsto\autoref{ClassficationOfPrelocSubcatsForLocNoethSch:LocFilt}\ \conditionalset{I\subset\OX\text{ in }\QCoh X}{I_{x}\in\mcF(x)\text{ for each }x\in X}
	\end{align*}
\end{Theorem}

\begin{proof}
	\autoref{LocOfPrelocSubcat} gives bijections between \autoref{ClassficationOfPrelocSubcatsForLocNoethSch:PrelocSubcat}, \autoref{ClassficationOfPrelocSubcatsForLocNoethSch:PrelocSubcatsOnOpenAffineSubschs}, and \autoref{ClassficationOfPrelocSubcatsForLocNoethSch:PrelocSubcatsOnLocRings}. \autoref{BijectionBetweenPrelocSubcatsAndFiltsOfSubobjOfCommRing} and \autoref{PrelocSubcatCorrespondingToImageOfFilt} (resp.\ \autoref{PrelocSubcatCorrespondingToLocOfFilt}) give a bijection between \autoref{ClassficationOfPrelocSubcatsForLocNoethSch:PrelocSubcatsOnOpenAffineSubschs} and \autoref{ClassficationOfPrelocSubcatsForLocNoethSch:FiltsOnOpenAffineSubschs} (resp.\ \autoref{ClassficationOfPrelocSubcatsForLocNoethSch:PrelocSubcatsOnLocRings} and \autoref{ClassficationOfPrelocSubcatsForLocNoethSch:FiltsOnLocRings}). \autoref{GluingOfFiltsAndDescriptionOfLocFilt} gives a bijection between \autoref{ClassficationOfPrelocSubcatsForLocNoethSch:LocFilt} and \autoref{ClassficationOfPrelocSubcatsForLocNoethSch:FiltsOnOpenAffineSubschs}.
\end{proof}

For a family of prelocalizing subcategories of $\QCoh X$, the supremum and the intersection are described in terms of local filters as follows.

\begin{Proposition}\label{JoinOfPrelocSubcatsAndJoinOfFilts}
	Assume that the following elements correspond to each other by the bijections in \autoref{ClassficationOfPrelocSubcatsForLocNoethSch} for each $\omega\in\Omega$.
	\begin{multicols}{2}
		\begin{enumerate}
			\item $\mcY^{\omega}$.
			\item $\set{\mcY^{\omega}_{\lambda}}_{\lambda\in\Lambda}$.
			\item $\set{\mcY^{\omega}(x)}_{x\in X}$.
			\item $\mcF^{\omega}$.
			\item $\set{\mcF^{\omega}_{\lambda}}_{\lambda\in\Lambda}$.
			\item $\set{\mcF^{\omega}(x)}_{x\in X}$.
		\end{enumerate}
	\end{multicols}
	\noindent Then the following elements correspond to each other by the bijections.
	\begin{multicols}{2}
		\begin{enumerate}
			\item $\preloc{\bigcup_{\omega\in\Omega}\mcY^{\omega}}$.
			\item $\set{\preloc{\bigcup_{\omega\in\Omega}\mcY^{\omega}_{\lambda}}}_{\lambda\in\Lambda}$.
			\item $\set{\preloc{\bigcup_{\omega\in\Omega}\mcY^{\omega}(x)}}_{x\in X}$.
			\item $\locfilt{\bigcup_{\omega\in\Omega}\mcF^{\omega}}$.
			\item $\set{\locfilt{\bigcup_{\omega\in\Omega}\mcF^{\omega}_{\lambda}}}_{\lambda\in\Lambda}$.
			\item $\set{\locfilt{\bigcup_{\omega\in\Omega}\mcF^{\omega}(x)}}_{x\in X}$.
		\end{enumerate}
	\end{multicols}
\end{Proposition}

\begin{proof}
	This follows from \autoref{LocOfJoinOfPrelocSubcats}.
\end{proof}

\begin{Proposition}\label{IntersectionOfPrelocSubcatsAndIntersectionOfFilts}
	Assume that the following elements correspond to each other by the bijections in \autoref{ClassficationOfPrelocSubcatsForLocNoethSch} for each $\omega\in\Omega$.
	\begin{multicols}{2}
		\begin{enumerate}
			\item $\mcY^{\omega}$.
			\item $\set{\mcY^{\omega}_{\lambda}}_{\lambda\in\Lambda}$.
			\item $\set{\mcY^{\omega}(x)}_{x\in X}$.
			\item $\mcF^{\omega}$.
			\item $\set{\mcF^{\omega}_{\lambda}}_{\lambda\in\Lambda}$.
			\item $\set{\mcF^{\omega}(x)}_{x\in X}$.
		\end{enumerate}
	\end{multicols}
	\noindent Then the following elements correspond to each other by the bijections.
	\begin{multicols}{2}
		\begin{enumerate}
			\item $\bigcap_{\omega\in\Omega}\mcY^{\omega}$.
			\item $\set{\bigcap_{\omega\in\Omega}\mcY^{\omega}_{\lambda}}_{\lambda\in\Lambda}$.
			\item $\set{\bigcap_{\omega\in\Omega}\mcY^{\omega}(x)}_{x\in X}$.
			\item $\bigcap_{\omega\in\Omega}\mcF^{\omega}$.
			\item $\set{\bigcap_{\omega\in\Omega}\mcF^{\omega}_{\lambda}}_{\lambda\in\Lambda}$.
			\item $\set{\bigcap_{\omega\in\Omega}\mcF^{\omega}(x)}_{x\in X}$.
		\end{enumerate}
	\end{multicols}
\end{Proposition}

\begin{proof}
	This follows from \autoref{LocOfIntersectionOfPrelocSubcats}.
\end{proof}

We demonstrate a calculation of the prelocalizing subcategories by using \autoref{ClassficationOfPrelocSubcatsForLocNoethSch}.

\begin{Example}\label{ExOfPrelocSubcatsOfAffineSch}
	Let $k$ be an algebraically closed field, and consider the polynomial ring $k[x]$ with a variable $x$. For each $a\in k$, let $\mfp_{a}:=(x-a)\subset k[x]$ and $\mfm_{a}:=\mfp_{a}k[x]_{\mfp_{a}}$. Then
	\begin{equation*}
		\Spec k[x]=\conditionalset{\mfp_{a}}{a\in k}\cup\set{0}.
	\end{equation*}
	Since $k[x]_{\mfp_{a}}$ is a discrete valuation ring, the set of ideals of $k[x]_{\mfp_{a}}$ is
	\begin{equation*}
		\conditionalset{\mfm_{a}^{i}}{i\in\mbZ_{\geq 0}}\cup\set{0},
	\end{equation*}
	where $\mfm_{a}^{0}=k[x]_{\mfp_{a}}$. For each $n\in\mbZ_{\geq 0}$, define the filter $\mcF_{a}^{n}$ of $k[x]_{\mfp_{a}}$ by
	\begin{equation*}
		\mcF_{a}^{n}=\conditionalset{\mfm_{a}^{i}}{0\leq i\leq n},
	\end{equation*}
	and let
	\begin{equation*}
		\mcF_{a}^{\infty}:=\conditionalset{\mfm_{a}^{i}}{i\in\mbZ_{\geq 0}},\ \mcF_{a}:=\conditionalset{\mfm_{a}^{i}}{i\in\mbZ_{\geq 0}}\cup\set{0}.
	\end{equation*}
	Then the set of filters of $k[x]_{\mfp_{a}}$ is
	\begin{equation*}
		\conditionalset{\mcF_{a}^{n}}{n\in\mbZ_{\geq 0}\cup\set{\infty}}\cup\set{\mcF_{a}}.
	\end{equation*}
	Since $k[x]_{0}=k(x)$ is a field, the set of the filters of $k(x)$ consists of $\mcF^{\infty}=\set{k(x)}$ and $\mcF=\set{0,\,k(x)}$. For each $a\in k$ and $n\in\mbZ_{\geq 0}\cup\set{\infty}$, $(\mcF_{a}^{n})_{0}=\mcF^{\infty}$, and $(\mcF_{a})_{0}=\mcF$. Hence the set
	\begin{equation*}
		\conditionalset[\Bigg]{(\set{\mcF_{a}^{r(a)}}_{a\in k},\,\mcF^{\infty})}{r=\set{r(a)}_{a\in k}\in\prod_{a\in k}(\mbZ_{\geq 0}\cup\set{\infty})}\cup\set{(\set{\mcF_{a}}_{a\in k},\,\mcF)}
	\end{equation*}
	is the set of families of filters which are compatible with localizations. By \autoref{ClassficationOfPrelocSubcatsForLocNoethSch}, the set of prelocalizing subcategories of $\Mod k[x]$ is
	\begin{equation*}
		\conditionalset[\Bigg]{\mcY_{r}}{r\in\prod_{a\in k}(\mbZ_{\geq 0}\cup\set{\infty})}\cup\set{\Mod k[x]},
	\end{equation*}
	where
	\begin{equation*}
		\mcY_{r}=\conditionalset{M\in\Mod k[x]}{ M_{\mfp_{a}}\mfm_{a}^{r(a)}=0\text{ for each }a\in k\text{ with }r(a)\neq\infty}
	\end{equation*}
	for each $r\in\prod_{a\in k}(\mbZ_{\geq 0}\cup\set{\infty})$.
\end{Example}

\section{Classification of localizing subcategories}
\label{sec:ClassificationOfLocSubcats}

In this section, we investigate extensions of prelocalizing subcategories (\autoref{SubcatClosedUnderExtAndLocSubcat} \autoref{ExtOfFullSubcats}) in terms of local filters and classify the localizing subcategories of $\QCoh X$ for a locally noetherian scheme $X$. The classification is given as a restriction of \autoref{ClassficationOfPrelocSubcatsForLocNoethSch}. We start with recalling Gabriel's classification of the localizing subcategories of $\Mod\Lambda$ for a ring $\Lambda$.

\begin{Definition}\label{ProductOfPrelocFiltAndGabrielFilt}
	Let $\Lambda$ be a ring.
	\begin{enumerate}
		\item\label{ProductOfPrelocFilt} For prelocalizing filters $\mcF_{1}$ and $\mcF_{2}$ of right ideals of $\Lambda$, define the \emph{product} $\mcF_{1}*\mcF_{2}$ as follows: $L\in\mcF_{1}*\mcF_{2}$ if and only if there exists $L_{1}\in\mcF_{1}$ satisfying $a^{-1}L\in\mcF_{2}$ for every $a\in L_{1}$.
		\item\label{GabrielFilt} A prelocalizing filter $\mcF$ of right ideals of $\Lambda$ is called a \emph{Gabriel filter} if $\mcF*\mcF\subset\mcF$ holds.
	\end{enumerate}
\end{Definition}

\begin{Proposition}\label{ProductOfPrelocFiltIsLarger}
	Let $\Lambda$ be a ring. If $\mcF_{1}$ and $\mcF_{2}$ are prelocalizing filters of right ideals of $\Lambda$, then $\mcF_{1}\subset\mcF_{1}*\mcF_{2}$, and $\mcF_{2}\subset\mcF_{1}*\mcF_{2}$.
\end{Proposition}

\begin{proof}
	Let $L_{1}\in\mcF_{1}$. Then $a^{-1}L_{1}=\Lambda\in\mcF_{2}$ for each $a\in L_{1}$. This shows that $\mcF_{1}\subset\mcF_{1}*\mcF_{2}$.
	
	Let $L_{2}\in\mcF_{2}$. Then $\Lambda\in\mcF_{1}$, and $a^{-1}L_{2}\in\mcF_{2}$ for each $a\in\Lambda$. This shows that $\mcF_{2}\subset\mcF_{1}*\mcF_{2}$.
\end{proof}

\begin{Theorem}[{\cite[p.~412]{Gabriel}}]\label{ExtOfPrelocSubcatsAndProdOfPrelocFiltsAndBijectionBetweenLocSubcatsAndGabrielFilts}
	Let $\Lambda$ be a ring.
	\begin{enumerate}
		\item\label{ExtOfPrelocSubcatsAndProdOfPrelocFilts} For each $i=1,2$, let $\mcY_{i}$ be a prelocalizing subcategory of $\Mod\Lambda$, and let $\mcF_{i}$ be the prelocalizing filter of right ideals of $\Lambda$ corresponding to $\mcF_{i}$ by the bijection in \autoref{BijectionBetweenPrelocSubcatsAndFiltsOfSubobjOfRing}. Then $\mcY_{1}*\mcY_{2}$ corresponds to $\mcF_{2}*\mcF_{1}$ by the bijection.
		\item\label{BijectionBetweenLocSubcatsAndGabrielFilts} The bijection in \autoref{BijectionBetweenPrelocSubcatsAndFiltsOfSubobjOfRing} restricts to a bijection
		\begin{equation*}
			\setwithspace{\text{localizing subcategories of }\Mod\Lambda}\to\setwithspace{\text{Gabriel filters of right ideals of }\Lambda}.
		\end{equation*}
	\end{enumerate}
\end{Theorem}

\begin{proof}
	\cite[Theorem~VI.5.1]{Stenstrom}.
\end{proof}

For a commutative ring $R$, we say that a filter $\mcF$ of $R$ is \emph{closed under products} if $I_{1},I_{2}\in\mcF$ implies $I_{1}I_{2}\in\mcF$. In the case of a commutative noetherian ring, products of filters and Gabriel filters are characterized as follows.

\begin{Proposition}\label{ProductOfFiltsAndGabrielFiltOfCommNoethRing}
	Let $R$ be a commutative noetherian ring.
	\begin{enumerate}
		\item\label{ProductOfFiltsOfCommNoethRing} Let $\mcF_{1}$ and $\mcF_{2}$ be filters of $R$. Then
		\begin{equation*}
			\mcF_{1}*\mcF_{2}=\conditionalset{I\subset R\text{ in }\Mod R}{I_{1}I_{2}\subset I\text{ for some }I_{1}\in\mcF_{1},\,I_{2}\in\mcF_{2}}.
		\end{equation*}
		\item\label{GabrielFiltOfCommNoethRing} Let $\mcF$ be a filter of $R$. Then $\mcF$ is a Gabriel filter if and only if $\mcF$ is closed under products.
	\end{enumerate}
\end{Proposition}

\begin{proof}
	\autoref{ProductOfFiltsOfCommNoethRing} Let $I\in\mcF_{1}*\mcF_{2}$. Then there exists $I_{1}\in\mcF_{1}$ such that $a^{-1}I\in\mcF_{2}$ for each $a\in I_{1}$. Since $R$ is noetherian, there exist $b_{1},\ldots,b_{n}\in R$ such that $I_{1}=b_{1}R+\cdots+b_{n}R$. Let $I_{2}:=b_{1}^{-1}I\cap\cdots\cap b_{n}^{-1}I$. Then $I_{2}\in\mcF_{2}$, and
	\begin{equation*}
		I_{1}I_{2}=b_{1}I_{2}+\cdots+b_{n}I_{2}\subset b_{1}(b_{1}^{-1}I)+\cdots+b_{n}(b_{n}^{-1}I)\subset I.
	\end{equation*}
	
	Conversely, let $J_{1}\in\mcF_{1}$ and $J_{2}\in\mcF_{2}$. For each $a\in J_{1}$, we have $J_{2}\subset a^{-1}J_{1}J_{2}$, and hence $a^{-1}J_{1}J_{2}\in\mcF_{2}$. This implies that $J_{1}J_{2}\in\mcF_{1}*\mcF_{2}$.
	
	\autoref{GabrielFiltOfCommNoethRing} This follows from \autoref{ProductOfFiltsOfCommNoethRing}.
\end{proof}

For a commutative noetherian ring $R$, the classification of the localizing subcategories of $\Mod R$ is stated as follows.

\begin{Corollary}\label{BijectionBetweenLocSubcatsAndFiltsOfSubobjOfCommNoethRing}
	Let $R$ be a commutative noetherian ring. Then the bijection in \autoref{BijectionBetweenPrelocSubcatsAndFiltsOfSubobjOfCommRing} restricts to a bijection
	\begin{equation*}
		\setwithspace{\text{localizing subcategories of }\Mod R}\to\setwithspace{\text{filters of }R\text{ closed under products}}.
	\end{equation*}
\end{Corollary}

\begin{proof}
	This follows from \autoref{ExtOfPrelocSubcatsAndProdOfPrelocFiltsAndBijectionBetweenLocSubcatsAndGabrielFilts} \autoref{BijectionBetweenLocSubcatsAndGabrielFilts} and \autoref{ProductOfFiltsAndGabrielFiltOfCommNoethRing} \autoref{GabrielFiltOfCommNoethRing}.
\end{proof}

In the rest of this section, let $X$ be a locally noetherian scheme, and let $\set{U_{\lambda}}_{\lambda\in\Lambda}$ be an open affine basis of $X$. For an object $M$ in $\QCoh X$ and a subobject $I$ of $\OX$, the subobject $MI$ of $M$ is defined as the image of the canonical morphism $M\otimes_{\OX} I\to M$ in $\QCoh X$.

\begin{Definition}\label{ProductOfFiltAndClosednessUnderProducts}\leavevmode
	\begin{enumerate}
		\item\label{ProductOfFiltAndClosednessUnderProducts:ProductOfFilt} Let $\mcF_{1}$ and $\mcF_{2}$ be local filters of $\OX$. We define the \emph{product} $\mcF_{1}*\mcF_{2}$ by
		\begin{equation*}
			\mcF_{1}*\mcF_{2}=\conditionalgeneratedset{locfilt}{I_{1}I_{2}\subset\OX\text{ in }\QCoh X}{I_{i}\in\mcF_{i}\text{ for each }i=1,2}.
		\end{equation*}
		\item\label{ProductOfFiltAndClosednessUnderProducts:ClosednessUnderProducts} We say that a local filter $\mcF$ is \emph{closed under products} if $\mcF*\mcF\subset\mcF$ holds.
	\end{enumerate}
\end{Definition}

Note that a local filter $\mcF$ is closed under products if and only if $I_{1},I_{2}\in\mcF$ implies $I_{1}I_{2}\in\mcF$.

Products of local filters of $\OX$ commute with the restriction to an open affine subscheme and the localization at a point.

\begin{Lemma}\label{LocOfProductOfLocFilts}
	Let $\mcF_{i}$ be a local filter of $\OX$ for each $i=1,2$.
	\begin{enumerate}
		\item\label{LocOfProductOfLocFiltsToOpenAffineSubsch} For every $\lambda\in\Lambda$,
		\begin{equation*}
			(\mcF_{1}*\mcF_{2})|_{U_{\lambda}}=\mcF_{1}|_{U_{\lambda}}*\mcF_{2}|_{U_{\lambda}}.
		\end{equation*}
		\item\label{LocOfProductOfLocFiltsToLocRing} For every $x\in X$,
		\begin{equation*}
			(\mcF_{1}*\mcF_{2})_{x}=(\mcF_{1})_{x}*(\mcF_{2})_{x}.
		\end{equation*}
	\end{enumerate}
\end{Lemma}

\begin{proof}
	\autoref{LocOfProductOfLocFiltsToOpenAffineSubsch} Let $J\in(\mcF_{1}*\mcF_{2})|_{U_{\lambda}}$. Then there exists $I\in\mcF_{1}*\mcF_{2}$ such that $I|_{U_{\lambda}}=J$. For each $x\in U_{\lambda}$, there exist an open affine neighborhood $U$ of $x$ in $X$ and $I_{1}\in\mcF_{1}$ and $I_{2}\in\mcF_{2}$ such that $(I_{1}I_{2})|_{U}\subset I|_{U}$. Hence
	\begin{equation*}
		(I_{1}|_{U_{\lambda}}I_{2}|_{U_{\lambda}})|_{U_{\lambda}\cap U}=(I_{1}I_{2})|_{U_{\lambda}\cap U}\subset I|_{U_{\lambda}\cap U}=J|_{U_{\lambda}\cap U}.
	\end{equation*}
	This shows that $J\in\mcF_{1}|_{U_{\lambda}}*\mcF_{2}|_{U_{\lambda}}$.
	
	Conversely, assume $J\in\mcF_{1}|_{U_{\lambda}}*\mcF_{2}|_{U_{\lambda}}$. Then for each $x\in U_{\lambda}$, there exist an open affine neighborhood $V$ of $x$ in $U_{\lambda}$ and $J_{1}\in\mcF_{1}|_{U_{\lambda}}$ and $J_{2}\in\mcF_{2}|_{U_{\lambda}}$ such that $(J_{1}J_{2})|_{V}\subset J|_{V}$. For each $i=1,2$, there exists $I_{i}\in\mcF_{i}$ such that $I_{i}|_{U_{\lambda}}=J_{i}$. Then $(I_{1}I_{2})|_{U_{\lambda}}\in(\mcF_{1}*\mcF_{2})|_{U_{\lambda}}$, and
	\begin{equation*}
		((I_{1}I_{2})|_{U_{\lambda}})|_{V}=(J_{1}J_{2})|_{V}\subset J|_{V}.
	\end{equation*}
	Since $(\mcF_{1}*\mcF_{2})|_{U_{\lambda}}$ is a local filter by \autoref{ImageOfFiltIsFilt} \autoref{ImageOfFiltToOpenAffineSubschIsFilt} and \autoref{FiltOfSubobjsForNoethSchIsLocFilt}, we obtain $J\in(\mcF_{1}*\mcF_{2})|_{U_{\lambda}}$.
	
	\autoref{LocOfProductOfLocFiltsToLocRing} This can be shown similarly to \autoref{LocOfProductOfLocFiltsToOpenAffineSubsch}.
\end{proof}

We describe extensions of prelocalizing subcategories of $\QCoh X$ in terms of products of local filters.

\begin{Theorem}\label{ExtOfPrelocSubcatsAndProductOfFilts}
	Assume that the following elements correspond to each other by the bijections in \autoref{ClassficationOfPrelocSubcatsForLocNoethSch} for each $i=1,2$.
	\begin{multicols}{2}
		\begin{enumerate}
			\item $\mcY^{i}$.
			\item $\set{\mcY^{i}_{\lambda}}_{\lambda\in\Lambda}$.
			\item $\set{\mcY^{i}(x)}_{x\in X}$.
			\item $\mcF^{i}$.
			\item $\set{\mcF^{i}_{\lambda}}_{\lambda\in\Lambda}$.
			\item $\set{\mcF^{i}(x)}_{x\in X}$.
		\end{enumerate}
	\end{multicols}
	\noindent Then the following elements correspond to each other by the bijections.
	\begin{multicols}{2}
		\begin{enumerate}
			\item $\mcY^{1}*\mcY^{2}$.
			\item $\set{\mcY^{1}_{\lambda}*\mcY^{2}_{\lambda}}_{\lambda\in\Lambda}$.
			\item $\set{\mcY^{1}(x)*\mcY^{2}(x)}_{x\in X}$.
			\item $\mcF^{1}*\mcF^{2}$.
			\item $\set{\mcF^{1}_{\lambda}*\mcF^{2}_{\lambda}}_{\lambda\in\Lambda}$.
			\item $\set{\mcF^{1}(x)*\mcF^{2}(x)}_{x\in X}$.
		\end{enumerate}
	\end{multicols}
\end{Theorem}

\begin{proof}
	This follows from \autoref{LocOfExtOfPrelocSubcats}, \autoref{ExtOfPrelocSubcatsAndProdOfPrelocFiltsAndBijectionBetweenLocSubcatsAndGabrielFilts} \autoref{ExtOfPrelocSubcatsAndProdOfPrelocFilts}, and \autoref{LocOfProductOfLocFilts}.
\end{proof}

\begin{Corollary}\label{ClassificationOfLocSubcatsForLocNoethSch}
	The bijections in \autoref{ClassficationOfPrelocSubcatsForLocNoethSch} restrict to bijections between following sets.
	\begin{enumerate}
		\item\label{ClassficationOfLocSubcatsForLocNoethSch:LocSubcat} The set of localizing subcategories of $\QCoh X$.
		\item\label{ClassficationOfLocSubcatsForLocNoethSch:LocSubcatsOnOpenAffineSubschs} The set of families $\set{\mcX_{\lambda}\subset\QCoh U_{\lambda}}_{\lambda\in\Lambda}$ of localizing subcategories such that $\mcX_{\lambda}|_{U_{\mu}}=\mcX_{\mu}$ for each $\lambda,\mu\in\Lambda$ with $U_{\mu}\subset U_{\lambda}$.
		\item\label{ClassficationOfLocSubcatsForLocNoethSch:LocSubcatsOnLocRings} The set of families $\set{\mcX(x)\subset\Mod\OXx}_{x\in X}$ of localizing subcategories such that $\mcX(y)_{x}=\mcX(x)$ for each $x,y\in X$ with $y\in\overline{\set{x}}$.
		\item\label{ClassficationOfLocSubcatsForLocNoethSch:LocFilt} The set of local filters of $\OX$ closed under products.
		\item\label{ClassficationOfLocSubcatsForLocNoethSch:FiltsOnOpenAffineSubschs} The set of families $\set{\mcF_{\lambda}}_{\lambda\in\Lambda}$, where $\mcF_{\lambda}$ is a filter of $\mcO_{U_{\lambda}}$ closed under products for each $\lambda\in\Lambda$, such that $\mcF_{\lambda}|_{U_{\mu}}=\mcF_{\mu}$ for each $\lambda,\mu\in\Lambda$ with $U_{\mu}\subset U_{\lambda}$.
		\item\label{ClassficationOfLocSubcatsForLocNoethSch:FiltsOnLocRings} The set of families $\set{\mcF(x)}_{x\in X}$, where $\mcF(x)$ is a filter of $\OXx$ closed under products for each $x\in X$, such that $\mcF(y)_{x}=\mcF(x)$ for each $x,y\in X$ with $y\in\overline{\set{x}}$.
	\end{enumerate}
\end{Corollary}

\begin{proof}
	This follows from \autoref{ExtOfPrelocSubcatsAndProductOfFilts}.
\end{proof}

We apply \autoref{ClassificationOfLocSubcatsForLocNoethSch} to \autoref{ExOfPrelocSubcatsOfAffineSch}.

\begin{Example}\label{ExOfLocSubcatsOfAffineSch}
	In the setting of \autoref{ExOfPrelocSubcatsOfAffineSch},
	\begin{equation*}
		\mcF_{a}^{m}*\mcF_{a}^{n}=\mcF_{a}^{m+n}
	\end{equation*}
	for each $a\in k$ and $m,n\in\mbZ_{\geq 0}\cup\set{\infty}$. Hence by \autoref{ClassificationOfLocSubcatsForLocNoethSch}, the set of localizing subcategories of $\Mod k[x]$ is
	\begin{equation*}
		\conditionalset[\Bigg]{\mcY_{r}}{r\in\prod_{a\in k}\set{0,\infty}}\cup\set{\Mod k[x]}.
	\end{equation*}
\end{Example}

In \autoref{BijectionBetweenLocSubcatsAndSpecializationClosedSubsetsForLocNoethSch}, we showed that there exists a bijection between the localizing subcategories of $\QCoh X$ and the specialization-closed subsets of $X$. For a local filter $\mcF$ of $\OX$ closed under products, the corresponding specialization-closed subset of $X$ is $\conditionalset{x\in X}{\mcF_{x}\neq\set{\OXx}}$.

Prime localizing subcategories of $\QCoh X$ are characterized in terms of local filters as follows.

\begin{Theorem}\label{CharacterizationOfLocFiltCorrespondingToPrimeLocSubcat}
	Let $\mcF$ be a local filter of $\OX$ closed under products. Then the following assertions are equivalent.
	\begin{enumerate}
		\item\label{CharacterizationOfLocFiltCorrespondingToPrimeLocSubcat:Prime} By the bijection in \autoref{ClassficationOfPrelocSubcatsForLocNoethSch}, the local filter $\mcF$ corresponds to a prime localizing subcategory of $\QCoh X$.
		\item\label{CharacterizationOfLocFiltCorrespondingToPrimeLocSubcat:Description} There exists $x\in X$ such that
		\begin{equation*}
			\mcF=\conditionalset{I\subset\OX\text{ in }\QCoh X}{I_{x}=\OXx}.
		\end{equation*}
		\item For each family $\{\mcF_{\omega}\}_{\omega\in\Omega}$ of local filters of $\OX$ closed under products satisfying $\mcF=\bigcap_{\omega\in\Omega}\mcF_{\omega}$, there exists $\omega\in\Omega$ such that $\mcF=\mcF_{\omega}$.
		\item For each family $\{\mcF_{\omega}\}_{\omega\in\Omega}$ of local filters of $\OX$ satisfying $\mcF=\bigcap_{\omega\in\Omega}\mcF_{\omega}$, there exists $\omega\in\Omega$ such that $\mcF=\mcF_{\omega}$.
	\end{enumerate}
\end{Theorem}

\begin{proof}
	This follows from \autoref{CharacterizationOfPrimeLocSubcat}.
\end{proof}

\section{Classification of closed subcategories}
\label{sec:ClassificationOfClosedSubcats}

In this section, we investigate the closed subcategories of $\QCoh X$ for a locally noetherian scheme $X$, whose definition is as follows.

\begin{Definition}\label{ClosedSubcat}
	Let $\mcA$ be a Grothendieck category. A prelocalizing subcategory $\mcX$ of $\mcA$ is called a \emph{closed} subcategory of $\mcA$ if $\mcX$ is closed under arbitrary direct products.
\end{Definition}

Note that every Grothendieck category has arbitrary direct products (\cite[Corollary~3.7.10]{Popescu}).

Closed subcategories are characterized by \autoref{CharacterizationOfPrelocSubcat} and the following result.

\begin{Proposition}\label{CharacterizationOfClosedSubcat}
	Let $\mcA$ be a Grothendieck category (or more generally, an abelian category admitting arbitrary direct products), and let $\mcY$ be a full subcategory of $\mcA$ closed under subobjects and quotient objects. Then the following assertions are equivalent.
	\begin{enumerate}
		\item\label{CharacterizationOfClosedSubcat:ClosedSubcat} $\mcY$ is closed under arbitrary direct products.
		\item\label{CharacterizationOfClosedSubcat:LeftAdjoint} The inclusion functor $\mcY\into\mcA$ has a left adjoint.
		\item\label{CharacterizationOfClosedSubcat:SmallestSubobj} For each object $M$ in $\mcA$, there exists a smallest subobject $L$ of $M$ satisfying $M/L\in\mcY$.
	\end{enumerate}
\end{Proposition}

\begin{proof}
	Apply \autoref{CharacterizationOfPrelocSubcat} to the opposite category of $\mcA$.
\end{proof}

Note that for a Grothendieck category, every full subcategory which is closed under subobjects and arbitrary direct products is also closed under arbitrary direct sums.

For a ring $\Lambda$, Rosenberg \cite{Rosenberg1} showed that there exists a bijection between the closed subcategories of $\Mod\Lambda$ and the two-sided ideals of $\Lambda$. This result can be unified into \autoref{BijectionBetweenPrelocSubcatsAndFiltsOfSubobjOfRing} as follows.

\begin{Theorem}[{Gabriel \cite[Lemma~V.2.1]{Gabriel} and Rosenberg \cite[Proposition~III.6.4.1]{Rosenberg1}}]\label{BijectionBetweenClosedSubcatsAndIdealsOfRing}
	Let $\Lambda$ be a ring. Then there exist bijections between the following sets.
	\begin{enumerate}
		\item\label{BijectionBetweenClosedSubcatsAndIdealsOfRing:ClosedSubcats} The set of closed subcategories of $\Mod\Lambda$.
		\item\label{BijectionBetweenClosedSubcatsAndIdealsOfRing:PrincipalPrelocFilts} The set of principal prelocalizing filters of right ideals of $\Lambda$.
		\item\label{BijectionBetweenClosedSubcatsAndIdealsOfRing:TwoSidedIdeals} The set of two-sided ideals of $\Lambda$.
	\end{enumerate}
	
	The bijection between \autoref{BijectionBetweenClosedSubcatsAndIdealsOfRing:ClosedSubcats} and \autoref{BijectionBetweenClosedSubcatsAndIdealsOfRing:PrincipalPrelocFilts} is induced by the bijection in \autoref{BijectionBetweenPrelocSubcatsAndFiltsOfSubobjOfRing}.
	
	The bijection between \autoref{BijectionBetweenClosedSubcatsAndIdealsOfRing:ClosedSubcats} and \autoref{BijectionBetweenClosedSubcatsAndIdealsOfRing:TwoSidedIdeals} is given by
	\begin{align*}
		\autoref{BijectionBetweenClosedSubcatsAndIdealsOfRing:ClosedSubcats}\to\autoref{BijectionBetweenClosedSubcatsAndIdealsOfRing:TwoSidedIdeals}&:\ \mcY\mapsto\bigcap_{M\in\mcY}\Ann_{\Lambda}(M),\\
		\autoref{BijectionBetweenClosedSubcatsAndIdealsOfRing:TwoSidedIdeals}\to\autoref{BijectionBetweenClosedSubcatsAndIdealsOfRing:ClosedSubcats}&:\ I\mapsto\conditionalset{M\in\Mod\Lambda}{MI=0}=\generatedset[\bigg]{\!\!preloc}{\frac{\Lambda}{I}}.
	\end{align*}
\end{Theorem}

\begin{proof}
	We show that for each right ideal $L$ of $\Lambda$, the principal filter $\mcF(L)$ of right ideals of $\Lambda$ is prelocalizing if and only if $L$ is a two-sided ideal of $\Lambda$. Assume that $\mcF(L)$ is prelocalizing. Then for each $a\in\Lambda$, we have $a^{-1}L\in\mcF(L)$. This implies $L\subset a^{-1}L$, and hence $aL\subset L$. Therefore $L$ is a two-sided ideal of $\Lambda$. The converse is obvious. The bijection between \autoref{BijectionBetweenClosedSubcatsAndIdealsOfRing:PrincipalPrelocFilts} and \autoref{BijectionBetweenClosedSubcatsAndIdealsOfRing:TwoSidedIdeals} follows from \autoref{CharacterizationOfPrincipalFilt}.
	
	Let $\mcY$ be a prelocalizing subcategory of $\mcA$, and let $\mcF$ be the corresponding prelocalizing filter of right ideals of $\Lambda$. If $\mcY$ is a closed subcategory of $\mcA$, then by \autoref{CharacterizationOfClosedSubcat}, there exists a smallest element of $\mcF$. Hence $\mcF$ is principal.
	
	Conversely, assume that $\mcF$ is principal. Then $\mcF=\mcF(I)$ for some two-sided ideal $I$ of $\Lambda$. Since
	\begin{align*}
		\mcY&=\conditionalset{M\in\Mod\Lambda}{I\subset\Ann_{\Lambda}(x)\text{ for each }x\in M}\\
		&=\conditionalset{M\in\Mod\Lambda}{MI=0},
	\end{align*}
	the prelocalizing subcategory $\mcY$ of $\mcA$ is also closed under arbitrary direct products.
\end{proof}

The aim of this section is to generalize \autoref{BijectionBetweenClosedSubcatsAndIdealsOfRing} to a locally noetherian scheme $X$. Let $\set{U_{\lambda}}_{\lambda\in\Lambda}$ be an open affine basis of $X$.

We show a lemma on gluing of subobjects on open affine subschemes.

\begin{Lemma}\label{GluingOfSubobjForSch}
	Let $M$ be an object in $\QCoh X$, and let $L_{\lambda}$ be a subobject of $M|_{U_{\lambda}}$ for each $\lambda\in\Lambda$. Assume that $L_{\lambda}|_{U_{\mu}}=L_{\mu}$ for each $\lambda,\mu\in\Lambda$ with $U_{\mu}\subset U_{\lambda}$. Then there exists a unique subobject $L$ of $M$ such that $L|_{U_{\lambda}}=L_{\lambda}$ for each $\lambda\in\Lambda$.
\end{Lemma}

\begin{proof}
	(Existence) Define a subsheaf $L$ of $M$ by
	\begin{equation*}
		L(U)=\conditionalset{s\in M(U)}{s|_{U_{\lambda}}\in L_{\lambda}(U_{\lambda})\text{ for each }\lambda\in\Lambda\text{ with }U_{\lambda}\subset U}
	\end{equation*}
	for each open subset $U$ of $X$. It is straightforward to show that $L$ is a subsheaf of $M$ satisfying $L|_{U_{\lambda}}=L_{\lambda}$ for each $\lambda\in\Lambda$. In particular, the sheaf $L$ is quasi-coherent.
	
	(Uniqueness) Let $L'$ be a subobject of $M$ in $\QCoh X$ such that $L'|_{U_{\lambda}}=L_{\lambda}$ for each $\lambda\in\Lambda$. Then
	\begin{align*}
		L'(U)&=\conditionalset{s\in M(U)}{s|_{U_{\lambda}}\in L'(U_{\lambda})\text{ for each }\lambda\in\Lambda\text{ with }U_{\lambda}\subset U}\\
		&=\conditionalset{s\in M(U)}{s|_{U_{\lambda}}\in L_{\lambda}(U_{\lambda})\text{ for each }\lambda\in\Lambda\text{ with }U_{\lambda}\subset U}
	\end{align*}
	for each open subset $U$ of $X$.
\end{proof}

The following lemma shows that for a principal filter of $\OX$, its restriction to an open affine subscheme and its localization at a point are also principal filters.

\begin{Lemma}\label{ImageOfPrincipalFiltIsPrincipalFilt}
	Let $I$ be a subobject of $\OX$.
	\begin{enumerate}
		\item\label{ImageOfPrincipalFiltToOpenAffineSubschIsPrincipalFilt} For every $\lambda\in\Lambda$, we have $\mcF(I)|_{U_{\lambda}}=\mcF(I|_{U_{\lambda}})$.
		\item\label{ImageOfPrincipalFiltToLocalRingIsPrincipalFilt} For every $x\in X$, we have $\mcF(I)_{x}=\mcF(I_{x})$.
	\end{enumerate}
\end{Lemma}

\begin{proof}
	\autoref{ImageOfPrincipalFiltToOpenAffineSubschIsPrincipalFilt} For each $J'\in\mcF(I)|_{U_{\lambda}}$, there exists $J\in\mcF(I)$ such that $J|_{U_{\lambda}}=J'$. Since $I\subset J$, it holds that $I|_{U_{\lambda}}\subset J|_{U_{\lambda}}=J'$. This shows that $\mcF(I)|_{U_{\lambda}}\subset\mcF(I|_{U_{\lambda}})$.
	
	It is follows from $I\in\mcF(I)$ that $I|_{U_{\lambda}}\in\mcF(I)|_{U_{\lambda}}$. Since $\mcF(I)|_{U_{\lambda}}$ is a filter of $\mcO_{U_{\lambda}}$ by \autoref{ImageOfFiltIsFilt} \autoref{ImageOfFiltToOpenAffineSubschIsFilt}, we have $\mcF(I)|_{U_{\lambda}}\supset\mcF(I|_{U_{\lambda}})$.
	
	\autoref{ImageOfPrincipalFiltToLocalRingIsPrincipalFilt} This is shown similarly by using \autoref{ImageOfFiltIsFilt} \autoref{ImageOfFiltToLocalRingIsFilt}.
\end{proof}

Conversely, if the restriction of a local filter of $\OX$ to each open affine subscheme $U_{\lambda}$ is principal, then the local filter is principal.

\begin{Lemma}\label{LocOfPrincipalFilt}
	Let $\mcF$ be a local filter of $\OX$. Then $\mcF$ is a principal filter if and only if the filter $\mcF|_{U_{\lambda}}$ of $\mcO_{U_{\lambda}}$ is principal for every $\lambda\in\Lambda$.
\end{Lemma}

\begin{proof}
	If $\mcF$ is a principal filter, then $\mcF|_{U_{\lambda}}$ is a principal filter for every $\lambda\in\Lambda$ by \autoref{ImageOfPrincipalFiltIsPrincipalFilt} \autoref{ImageOfPrincipalFiltToOpenAffineSubschIsPrincipalFilt}.
	
	Assume that there exists a subobject $I_{\lambda}$ of $\mcO_{U_{\lambda}}$ such that $\mcF|_{U_{\lambda}}=\mcF(I_{\lambda})$ for each $\lambda\in\Lambda$. For each $\lambda,\mu\in\Lambda$ with $U_{\mu}\subset U_{\lambda}$,
	\begin{equation*}
		\mcF(I_{\lambda}|_{U_{\mu}})=\mcF(I_{\lambda})|_{U_{\mu}}=(\mcF|_{U_{\lambda}})|_{U_{\mu}}=\mcF|_{U_{\mu}}=\mcF(I_{\mu}).
	\end{equation*}
	Hence $I_{\lambda}|_{U_{\mu}}=I_{\mu}$. By \autoref{GluingOfSubobjForSch}, there exists a subobject $I$ of $\OX$ such that $I|_{U_{\lambda}}=I_{\lambda}$ for each $\lambda\in\Lambda$. Since $\mcF(I)|_{U_{\lambda}}=\mcF(I|_{U_{\lambda}})=\mcF(I_{\lambda})=\mcF_{\lambda}$ for each $\lambda\in\Lambda$, it follows from \autoref{GluingOfFiltsAndDescriptionOfLocFilt} \autoref{GluingOfFilts} that $\mcF(I)=\mcF$.
\end{proof}

\begin{Remark}\label{BeingPrincipalIsNotLocalProperty}
	Let $\mcF$ be a local filter of $\OX$. Even if $\mcF_{x}$ is a principal filter of $\OXx$ for each $x\in X$, the local filter $\mcF$ is not necessarily a principal filter. A counter-example is given in \autoref{ExOfClosedSubcatsOfAffineSch}.
\end{Remark}

We characterize closed subcategories of $\QCoh X$ in terms of local filters.

\begin{Lemma}\label{ClosedSubcatAndPrincipalFilt}
	Let $\mcY$ be a prelocalizing subcategory of $\QCoh X$, and let $\mcF$ be the corresponding local filter of $\OX$ by the bijection in \autoref{ClassficationOfPrelocSubcatsForLocNoethSch}. Then $\mcY$ is a closed subcategory of $\QCoh X$ if and only if $\mcF$ is a principal filter. If $\mcF=\mcF(I)$ for a subobject $I$ of $\OX$, then $I$ is the smallest subobject of $\OX$ satisfying $\OX/I\in\mcY$, and
	\begin{equation*}
		\mcY=\conditionalset{M\in\QCoh X}{MI=0}.
	\end{equation*}
\end{Lemma}

\begin{proof}
	Assume that $\mcY$ is a closed subcategory of $\QCoh X$. Then by \autoref{CharacterizationOfClosedSubcat}, there exists a smallest subobject $I$ of $\OX$ satisfying $\OX/I\in\mcY$. Hence $\mcF=\mcF(I)$.
	
	Conversely, assume that $\mcF=\mcF(I)$ for some subobject $I$ of $\OX$. Then for each $\lambda\in\Lambda$, we have $\mcF|_{U_{\lambda}}=\mcF(I|_{U_{\lambda}})$ by \autoref{ImageOfPrincipalFiltIsPrincipalFilt} \autoref{ImageOfPrincipalFiltToOpenAffineSubschIsPrincipalFilt}, and hence
	\begin{equation*}
		\mcY|_{U_{\lambda}}=\conditionalset{M'\in\QCoh U_{\lambda}}{M'(I|_{U_{\lambda}})=0}
	\end{equation*}
	by \autoref{BijectionBetweenClosedSubcatsAndIdealsOfRing}. By \autoref{ClassficationOfPrelocSubcatsForLocNoethSch},
	\begin{align*}
		\mcY&=\conditionalset{M\in\QCoh X}{M|_{U_{\lambda}}\in\mcY|_{U_{\lambda}}\text{ for every }\lambda\in\Lambda}\\
		&=\conditionalset{M\in\QCoh X}{M|_{U_{\lambda}}I|_{U_{\lambda}}=0\text{ for every }\lambda\in\Lambda}\\
		&=\conditionalset{M\in\QCoh X}{MI=0}.
	\end{align*}
	For each object $M$ in $\QCoh X$, the subobject $MI$ of $M$ is the smallest among the subobjects $L$ of $M$ satisfying $(M/L)I=0$. Therefore $\mcY$ is a closed subcategory of $\QCoh X$.
\end{proof}

As in \autoref{BeingPrincipalIsNotLocalProperty}, the same type of theorem as \autoref{ClassificationOfLocSubcatsForLocNoethSch} does not hold for the closed subcategories. For this reason, we use the characterization in \autoref{DescriptionOfPrelocSubcatsOnLocCats} in order to obtain a generalization to the closed subcategories.

\begin{Theorem}\label{ClassificationOfClosedSubcatsForLocNoethSch}
	Let $X$ be a locally noetherian scheme, and let $\set{U_{\lambda}}_{\lambda\in\Lambda}$ be an open affine basis of $X$. Then there exist bijections between the following sets.
	\begin{enumerate}
		\item\label{ClassficationOfClosedSubcatsForLocNoethSch:ClosedSubcat} The set of closed subcategories of $\QCoh X$.
		\item\label{ClassficationOfClosedSubcatsForLocNoethSch:ClosedSubcatsOnOpenAffineSubschs} The set of families $\set{\mcY_{\lambda}\subset\QCoh U_{\lambda}}_{\lambda\in\Lambda}$ of closed subcategories such that $\mcY_{\lambda}|_{U_{\mu}}=\mcY_{\mu}$ for each $\lambda,\mu\in\Lambda$ with $U_{\mu}\subset U_{\lambda}$.
		\item\label{ClassficationOfClosedSubcatsForLocNoethSch:ClosedSubcatsOnLocRings} The set of families $\set{\mcY(x)\subset\Mod\OXx}_{x\in X}$ of closed subcategories such that for each $x\in X$, there exist $\lambda\in\Lambda$ with $x\in U_{\lambda}$ and a closed subcategory $\mcY'$ of $\QCoh U_{\lambda}$ satisfying $\mcY'_{y}=\mcY(y)$ for each $y\in U_{\lambda}$.
		\item\label{ClassficationOfClosedSubcatsForLocNoethSch:PrincipalFilt} The set of principal filters of $\OX$.
		\item\label{ClassficationOfClosedSubcatsForLocNoethSch:PrincipalFiltsOnOpenAffineSubschs} The set of families $\set{\mcF_{\lambda}}_{\lambda\in\Lambda}$, where $\mcF_{\lambda}$ is a principal filter of $\mcO_{U_{\lambda}}$ for each $\lambda\in\Lambda$, such that $\mcF_{\lambda}|_{U_{\mu}}=\mcF_{\mu}$ for each $\lambda,\mu\in\Lambda$ with $U_{\mu}\subset U_{\lambda}$.
		\item\label{ClassficationOfClosedSubcatsForLocNoethSch:PrincipalFiltsOnLocRings} The set of families $\set{\mcF(x)}_{x\in X}$, where $\mcF(x)$ is a principal filter of $\OXx$ for each $x\in X$, such that for each $x\in X$, there exist $\lambda\in\Lambda$ with $x\in U_{\lambda}$ and a principal filter $\mcF'$ of $\mcO_{U_{\lambda}}$ satisfying $\mcF'_{y}=\mcF(y)$ for each $y\in U_{\lambda}$.
		\item\label{ClassficationOfClosedSubcatsForLocNoethSch:Ideal} The set of subobjects of $\OX$.
		\item\label{ClassficationOfClosedSubcatsForLocNoethSch:IdealsOfOpenAffineSubschs} The set of families $\set{I_{\lambda}}_{\lambda\in\Lambda}$, where $I_{\lambda}$ is a subobject of $\mcO_{U_{\lambda}}$ for each $\lambda\in\Lambda$, such that $I_{\lambda}|_{U_{\mu}}=I_{\mu}$ for each $\lambda,\mu\in\Lambda$ with $U_{\mu}\subset U_{\lambda}$.
		\item\label{ClassficationOfClosedSubcatsForLocNoethSch:IdealsOfLocRings} The set of families $\set{I(x)}_{x\in X}$, where $I(x)$ is an ideal of $\OXx$ for each $x\in X$, such that for each $x\in X$, there exist $\lambda\in\Lambda$ with $x\in U_{\lambda}$ and a subobject $I'$ of $\mcO_{U_{\lambda}}$ satisfying $I'_{y}=I(y)$ for each $y\in U_{\lambda}$.
	\end{enumerate}
	
	The bijections between the sets \autoref{ClassficationOfClosedSubcatsForLocNoethSch:ClosedSubcat}, \ldots, \autoref{ClassficationOfClosedSubcatsForLocNoethSch:PrincipalFiltsOnLocRings} are induced by \autoref{ClassficationOfPrelocSubcatsForLocNoethSch}.
	
	The bijections \autoref{ClassficationOfClosedSubcatsForLocNoethSch:PrincipalFilt}$\leftrightarrow$\autoref{ClassficationOfClosedSubcatsForLocNoethSch:Ideal}, \autoref{ClassficationOfClosedSubcatsForLocNoethSch:PrincipalFiltsOnOpenAffineSubschs}$\leftrightarrow$\autoref{ClassficationOfClosedSubcatsForLocNoethSch:IdealsOfOpenAffineSubschs}, and \autoref{ClassficationOfClosedSubcatsForLocNoethSch:PrincipalFiltsOnLocRings}$\leftrightarrow$\autoref{ClassficationOfClosedSubcatsForLocNoethSch:IdealsOfLocRings} are defined by the bijection $L\mapsto\mcF(L)$ in \autoref{CharacterizationOfPrincipalFilt}.
\end{Theorem}

\begin{proof}
	This follows from \autoref{ClassficationOfPrelocSubcatsForLocNoethSch}, \autoref{BijectionBetweenClosedSubcatsAndIdealsOfRing}, \autoref{LocOfPrincipalFilt}, and \autoref{ClosedSubcatAndPrincipalFilt}.
\end{proof}

We establish a bijection between the closed subcategories of $\QCoh X$ and the closed subschemes of $X$ by using the following fact.

\begin{Proposition}\label{BijectionBetweenQCohSubsheavesAndClosedSubschs}
	There is a bijection
	\begin{equation*}
		\setwithspace{\text{subobjects of }\OX}\to\setwithspace{\text{closed subschemes of }X}
	\end{equation*}
	given by $I\mapsto (\Supp(\OX/I),i^{-1}(\OX/I))$, where $i\colon\Supp(\OX/I)\into X$ is the immersion. For each closed subscheme $Y$ of $X$, the corresponding subobject $I$ of $\OX$ is given by the exact sequence
	\begin{equation*}
		0\to I\to\OX\to i_{*}\mcO_{Y}\to 0,
	\end{equation*}
	where $i\colon Y\into X$ is the immersion, and $\OX\to i_{*}\mcO_{Y}$ is the canonical morphism.
\end{Proposition}

\begin{proof}
	\cite[Proposition~II.5.9]{Hartshorne2}.
\end{proof}

\begin{Theorem}\label{BijectionBetweenClosedSubcatsAndClosedSubschs}
	Let $X$ be a locally noetherian scheme. Then there exists a bijection between
	\begin{enumerate}
		\item\label{BijectionBetweenClosedSubcatsAndClosedSubschs:ClosedSubcats} The set of closed subcategories of $\QCoh X$.
		\item\label{BijectionBetweenClosedSubcatsAndClosedSubschs:ClosedSubschs} The set of closed subschemes of $X$.
	\end{enumerate}
	
	For each closed subscheme $Y$ of $X$ with the immersion $i\colon Y\into X$, the functor $i_{*}\colon\QCoh Y\to\QCoh X$ is fully faithful and induces an equivalence between $\QCoh Y$ and the closed subcategory of $\QCoh X$ corresponding to $Y$.
\end{Theorem}

\begin{proof}
	The bijection is obtained by \autoref{ClassificationOfClosedSubcatsForLocNoethSch} and \autoref{BijectionBetweenQCohSubsheavesAndClosedSubschs}. By \cite[0.5.1.4]{Grothendieck1}, \cite[Proposition~I.5.5.1 (i)]{Grothendieck1}, and \cite[Corollary~I.9.2.2 (a)]{Grothendieck1}, we have the functor $i^{*}\colon\QCoh X\to\QCoh Y$ and its right adjoint $i_{*}\colon\QCoh Y\to\QCoh X$. It is straightforward to show that the counit morphism $i^{*}i_{*}\to 1_{\QCoh Y}$ is an isomorphism. Hence $i_{*}$ is fully faithful. An object $M$ in $\QCoh X$ is isomorphic to the image of an object in $\QCoh Y$ by $i_{*}$ if and only if the canonical morphism $M\to i_{*}i^{*}M$ is an isomorphism. Let $I$ be the subobject of $\OX$ corresponding to $Y$. Since we have the exact sequence
	\begin{equation*}
		0\to MI\to M\to i_{*}i^{*}M\to 0,
	\end{equation*}
	$M\to i_{*}i^{*}M$ is an isomorphism if and only if $MI=0$. Therefore the claim follows.
\end{proof}

\begin{Example}\label{ExOfClosedSubcatsOfAffineSch}
	We follow the notations in \autoref{ExOfPrelocSubcatsOfAffineSch} and \autoref{ExOfLocSubcatsOfAffineSch}. Each nonzero proper ideal $I$ of $k[x]$ is generated by an element of the form $(x-a_{1})^{r_{1}}\cdots(x-a_{l})^{r_{l}}$, where $l\in\mbZ_{\geq 1}$, $a_{1},\ldots,a_{l}$ are distinct elements of $k$, and $r_{1},\ldots,r_{l}\in\mbZ_{\geq 1}$. We have
	\begin{equation*}
		I_{\mfp_{a}}=
		\begin{cases}
			\mfm_{a_{i}}^{r_{i}} & \text{if }a=a_{i}\text{ for some }i\in\set{1,\ldots,l}\\
			k[x]_{\mfp_{a}} & \text{if }a\in k\setminus\set{a_{1},\ldots,a_{l}}
		\end{cases}.
	\end{equation*}
	For each $r\in\prod_{a\in k}(\mbZ_{\geq 0}\cup\set{\infty})$, the object $k[x]/I$ belongs to $\mcY_{r}$ if and only if $r_{i}\leq r(a_{i})$ for every $i=1,\ldots,l$. Hence the corresponding filter of $k[x]$ to $\mcY_{r}$ is
	\begin{equation*}
		\conditionalset[\bigg]{(x-a_{1})^{r_{1}}\cdots(x-a_{l})^{r_{l}}\subset k[x]}{
			\begin{gathered}
				l\in\mbZ_{\geq 1},\ a_{1},\ldots,a_{l}\in k\text{ (distinct)},\ r_{1},\ldots,r_{l}\in\mbZ_{\geq 1}\\
				r_{i}\leq r(a_{i})\text{ for each }i=1,\ldots,l
			\end{gathered}
		}\cup\set{k[x]}.
	\end{equation*}
	This is equal to
	\begin{equation*}
		\conditionalgeneratedset{locfilt}{(x-a)^{r}\subset k[x]}{a\in k,\ r\in\mbZ_{\geq 1},\ r\leq r(a)},
	\end{equation*}
	and we have the description
	\begin{equation*}
		\mcY_{r}=\conditionalgeneratedset[\bigg]{\!\!preloc}{\frac{k[x]}{(x-a)^{r}}}{a\in k,\ r\in\mbZ_{\geq 1},\ r\leq r(a)}.
	\end{equation*}
	By \autoref{ClassificationOfClosedSubcatsForLocNoethSch}, the set of closed subcategories of $\Mod k[x]$ is
	\begin{equation*}
		\conditionalset[\Bigg]{\mcY_{r}}{r\in\bigoplus_{a\in k}\mbZ_{\geq 0}}\cup\set{\Mod k[x]}.
	\end{equation*}
	
	Let $r\in(\prod_{a\in k}\mbZ_{\geq 0})\setminus(\bigoplus_{a\in k}\mbZ_{\geq 0})$. Then for every $\mfp\in\Spec k[x]$, the prelocalizing subcategory $(\mcY_{r})_{\mfp}$ of $\Mod k[x]_{\mfp}$ is a closed subcategory, and $\mcF_{\mfp}$ is a principal filter of $k[x]_{\mfp}$. However, the prelocalizing subcategory $\mcY_{r}$ of $\Mod k[x]$ is not a closed subcategory, and the corresponding local filter of $k[x]$ is not a principal filter.
\end{Example}

\section{Classification of bilocalizing subcategories}
\label{sec:ClassificationOfBilocSubcats}

Let $X$ be a locally noetherian scheme. We investigate extensions of closed subcategories. The following lemma shows that products of principal filters are also principal.

\begin{Lemma}\label{ProdOfPrincipalFilts}
	Let $I_{1}$ and $I_{2}$ be subobjects of $\OX$. Then $\mcF(I_{1})*\mcF(I_{2})=\mcF(I_{1}I_{2})$.
\end{Lemma}

\begin{proof}
	This follows from \autoref{PrincipalFiltIsLocFilt}.
\end{proof}

Extensions of closed subcategories are described in terms of products of principal filters.

\begin{Theorem}\label{ExtOfClosedSubcatsAndProductOfPrincipalFilts}
	Assume that the following elements correspond to each other by the bijections in \autoref{ClassificationOfClosedSubcatsForLocNoethSch} for each $i=1,2$.
	\begin{multicols}{3}
		\begin{enumerate}
			\item $\mcY^{i}$.
			\item $\set{\mcY^{i}_{\lambda}}_{\lambda\in\Lambda}$.
			\item $\set{\mcY^{i}(x)}_{x\in X}$.
			\item $\mcF^{i}$.
			\item $\set{\mcF^{i}_{\lambda}}_{\lambda\in\Lambda}$.
			\item $\set{\mcF^{i}(x)}_{x\in X}$.
			\item $I^{i}$.
			\item $\set{I^{i}_{\lambda}}_{\lambda\in\Lambda}$.
			\item $\set{I^{i}(x)}_{x\in X}$.
		\end{enumerate}
	\end{multicols}
	\noindent Then the following elements correspond to each other by the bijections.
	\begin{multicols}{3}
		\begin{enumerate}
			\item $\mcY^{1}*\mcY^{2}$.
			\item $\set{\mcY^{1}_{\lambda}*\mcY^{2}_{\lambda}}_{\lambda\in\Lambda}$.
			\item $\set{\mcY^{1}(x)*\mcY^{2}(x)}_{x\in X}$.
			\item $\mcF^{1}*\mcF^{2}$.
			\item $\set{\mcF^{1}_{\lambda}*\mcF^{2}_{\lambda}}_{\lambda\in\Lambda}$.
			\item $\set{\mcF^{1}(x)*\mcF^{2}(x)}_{x\in X}$.
			\item $I^{1}I^{2}$.
			\item $\set{I^{1}_{\lambda}I^{2}_{\lambda}}_{\lambda\in\Lambda}$.
			\item $\set{I^{1}(x)I^{2}(x)}_{x\in X}$.
		\end{enumerate}
	\end{multicols}
\end{Theorem}

\begin{proof}
	This follows from \autoref{ExtOfPrelocSubcatsAndProductOfFilts} and \autoref{ProdOfPrincipalFilts}.
\end{proof}

As a corollary of \autoref{ExtOfClosedSubcatsAndProductOfPrincipalFilts}, we obtain a classification of the bilocalizing subcategories of $\QCoh X$. They are defined as follows.

\begin{Definition}\label{BilocSubcat}
	Let $\mcA$ be a Grothendieck category. A prelocalizing subcategory $\mcX$ of $\mcA$ is called a \emph{bilocalizing} subcategory of $\mcA$ if $\mcX$ is both localizing and closed.
\end{Definition}

Bilocalizing subcategories have the following characterization.

\begin{Proposition}\label{CharacterizationOfBilocSubcat}
	Let $\mcA$ be a Grothendieck category, and let $\mcX$ be a localizing subcategory of $\mcA$. Then $\mcX$ is a bilocalizing subcategory of $\mcA$ if and only if the canonical functor $\mcA\to\mcA/\mcX$ has a left adjoint.
\end{Proposition}

\begin{proof}
	\cite[Theorem~4.21.1]{Popescu}.
\end{proof}

For a ring $\Lambda$, the bilocalizing subcategories of $\Mod\Lambda$ are classified by the idempotent two-sided ideals of $\Lambda$.

\begin{Definition}\label{IdempIdeal}
	Let $\Lambda$ be a ring. A two-sided ideal $I$ of $\Lambda$ is called \emph{idempotent} if $I^{2}=I$ holds.
\end{Definition}

\begin{Proposition}\label{ProductOfIdealsAndExtOfClosedSubcatsAndBijectionBetweenBilocSubcatsAndIdempIdealsOfRing}
	Let $\Lambda$ be a ring.
	\begin{enumerate}
		\item\label{ProductOfIdealsAndExtOfClosedSubcats} For each $i=1,2$, let $\mcY_{i}$ be a closed subcategory of $\Mod\Lambda$, and let $I_{i}$ be the corresponding two-sided ideal of $\Lambda$ by the bijection in \autoref{BijectionBetweenClosedSubcatsAndIdealsOfRing}. Then $\mcY_{1}*\mcY_{2}$ corresponds to $I_{2}I_{1}$ by the bijection.
		\item\label{BijectionBetweenBilocSubcatsAndIdempIdealsOfRing} The bijection in \autoref{BijectionBetweenClosedSubcatsAndIdealsOfRing} restricts to a bijection
		\begin{equation*}
			\setwithspace{\text{bilocalizing subcategories of }\Mod\Lambda}\to\setwithspace{\text{idempotent two-sided ideals of }\Lambda}.
		\end{equation*}
	\end{enumerate}
\end{Proposition}

\begin{proof}
	\autoref{ProductOfIdealsAndExtOfClosedSubcats} For two-sided ideals $I_{1}$ and $I_{2}$ of $\Lambda$, it is straightforward to show that $\mcF(I_{1})*\mcF(I_{2})=\mcF(I_{1}I_{2})$. Therefore the claim follows from \autoref{ExtOfPrelocSubcatsAndProdOfPrelocFiltsAndBijectionBetweenLocSubcatsAndGabrielFilts} \autoref{ExtOfPrelocSubcatsAndProdOfPrelocFilts}.
	
	\autoref{BijectionBetweenBilocSubcatsAndIdempIdealsOfRing} This follows from \autoref{ProductOfIdealsAndExtOfClosedSubcats}.
\end{proof}

A subobject $I$ of $\OX$ is called \emph{idempotent} if $I^{2}=I$ holds. We classify the bilocalizing subcategories of $\QCoh X$ as follows.

\begin{Corollary}\label{ClassificationOfBilocSubcatsForLocNoethSch}
	The bijections in \autoref{ClassificationOfClosedSubcatsForLocNoethSch} restrict to bijections between following sets.
	\begin{enumerate}
		\item\label{ClassificationOfBilocSubcatsForLocNoethSch:BilocSubcat} The set of bilocalizing subcategories of $\QCoh X$.
		\item\label{ClassificationOfBilocSubcatsForLocNoethSch:BilocSubcatsOnOpenAffineSubschs} The set of families $\set{\mcY_{\lambda}\subset\QCoh U_{\lambda}}_{\lambda\in\Lambda}$ of bilocalizing subcategories such that $\mcY_{\lambda}|_{U_{\mu}}=\mcY_{\mu}$ for each $\lambda,\mu\in\Lambda$ with $U_{\mu}\subset U_{\lambda}$.
		\item\label{ClassificationOfBilocSubcatsForLocNoethSch:BilocSubcatsOnLocRings} The set of families $\set{\mcY(x)\subset\Mod\OXx}_{x\in X}$ of bilocalizing subcategories such that for each $x\in X$, there exist $\lambda\in\Lambda$ with $x\in U_{\lambda}$ and a bilocalizing subcategory $\mcY'$ of $\QCoh U_{\lambda}$ satisfying $\mcY'_{y}=\mcY(y)$ for each $y\in U_{\lambda}$.
		\item\label{ClassificationOfBilocSubcatsForLocNoethSch:PrincipalFiltClosedUnderProducts} The set of principal filters of $\OX$ closed under products.
		\item\label{ClassificationOfBilocSubcatsForLocNoethSch:PrincipalFiltsClosedUnderProductsOnOpenAffineSubschs} The set of families $\set{\mcF_{\lambda}}_{\lambda\in\Lambda}$, where $\mcF_{\lambda}$ is a principal filter of $\mcO_{U_{\lambda}}$ closed under products for each $\lambda\in\Lambda$, such that $\mcF_{\lambda}|_{U_{\mu}}=\mcF_{\mu}$ for each $\lambda,\mu\in\Lambda$ with $U_{\mu}\subset U_{\lambda}$.
		\item\label{ClassificationOfBilocSubcatsForLocNoethSch:PrincipalFiltsClosedUnderProductsOnLocRings} The set of families $\set{\mcF(x)}_{x\in X}$, where $\mcF(x)$ is a principal filter of $\OXx$ closed under products for each $x\in X$, such that for each $x\in X$, there exist $\lambda\in\Lambda$ with $x\in U_{\lambda}$ and a principal filter of subobjects $\mcF'$ of $\mcO_{U_{\lambda}}$ which is closed under products and satisfies $\mcF'_{y}=\mcF(y)$ for each $y\in U_{\lambda}$.
		\item\label{ClassificationOfBilocSubcatsForLocNoethSch:IdempIdeal} The set of idempotent subobjects of $\OX$.
		\item\label{ClassificationOfBilocSubcatsForLocNoethSch:IdempIdealsOfOpenAffineSubschs} The set of families $\set{I_{\lambda}}_{\lambda\in\Lambda}$, where $I_{\lambda}$ is an idempotent subobjects of $\mcO_{U_{\lambda}}$ for each $\lambda\in\Lambda$, such that $I_{\lambda}|_{U_{\mu}}=I_{\mu}$ for each $\lambda,\mu\in\Lambda$ with $U_{\mu}\subset U_{\lambda}$.
		\item\label{ClassificationOfBilocSubcatsForLocNoethSch:IdempIdealsOfLocRings} The set of families $\set{I(x)}_{x\in X}$, where $I(x)$ is an idempotent ideal of $\OXx$ for each $x\in X$, such that for each $x\in X$, there exist $\lambda\in\Lambda$ with $x\in U_{\lambda}$ and an idempotent subobject $I'$ of $\mcO_{U_{\lambda}}$ satisfying $I'_{y}=I(y)$ for each $y\in U_{\lambda}$.
	\end{enumerate}
\end{Corollary}

\begin{proof}
	This follows from \autoref{ExtOfClosedSubcatsAndProductOfPrincipalFilts}.
\end{proof}

\begin{Example}\label{ExOfBilocSubcatsOfAffineSch}
	In the setting of \autoref{ExOfClosedSubcatsOfAffineSch}, the set of bilocalizing subcategories of $\Mod k[x]$ is
	\begin{equation*}
		\conditionalset{\mcY_{r}}{r=\set{0}_{a\in k}}\cup\set{\Mod k[x]}=\set{0,\,\Mod k[x]}.
	\end{equation*}
\end{Example}

We show that the sets in \autoref{ClassificationOfBilocSubcatsForLocNoethSch} also bijectively correspond to the set of open closed subsets of $X$. We start with the following well-known fact on a commutative noetherian ring.

\begin{Lemma}\label{IdempIdealOfCommNoethRingSplits}
	Let $R$ be a commutative noetherian ring, and let $I$ be an idempotent ideal of $R$. Then there exists an ideal $J$ of $R$ such that $R=I\oplus J$ in $\Mod R$. In particular, the subset $\Supp(R/I)$ of $\Spec R$ is open and closed.
\end{Lemma}

\begin{proof}
	By Nakayama's lemma (\cite[Theorem~2.2]{Matsumura}), there exists $a\in R$ such that $aI=0$ and $1-a\in I$. Then $a^{2}=a$ and $(1-a)R=I$. By letting $J=aR$, we obtain $R=I\oplus J$, and $\Spec R$ is the disjoint union of the closed subsets $V(I)$ and $V(J)$ determined by $I$ and $J$, respectively.
\end{proof}

The idempotence of a subobject of $\OX$ is characterized in terms of the corresponding closed subscheme.

\begin{Lemma}\label{SubobjIsIdempIffCorrespondingClosedSubschIsOpenSubsch}
	Let $X$ be a locally noetherian scheme. Let $I$ be a subobject of $\OX$, and let $Y$ be the corresponding closed subscheme of $X$ by the bijection in \autoref{BijectionBetweenQCohSubsheavesAndClosedSubschs}. Then $I$ is idempotent if and only if $Y$ is also an open subscheme of $X$.
\end{Lemma}

\begin{proof}
	Assume that $I$ is idempotent. For each open affine subscheme $U$ of $X$, the subobject $I|_{U}$ of $\mcO_{U}$ is idempotent. By \autoref{IdempIdealOfCommNoethRingSplits}, the subset $\Supp(\mcO_{U}/I|_{U})$ of $U$ is open and closed. Since
	\begin{equation*}
		\Supp\frac{\mcO_{U}}{I|_{U}}=U\cap\Supp\frac{\OX}{I},
	\end{equation*}
	the underlying space $\Supp(\OX/I)$ of $Y$ is an open subset of $X$. For each $y\in Y$, the ideal $I_{y}$ of $\OXy$ is idempotent, and $(\OX/I)_{y}\neq 0$. Hence $I_{y}=0$. It follows that $\mcO_{Y}=(\OX/I)|_{Y}=\OX|_{Y}$.
	
	Conversely, assume that $Y$ is also an open subscheme. Let $i\colon Y\into X$ be the immersion. There is an exact sequence
	\begin{equation*}
		0\to I\to\OX\to i_{*}(\OX|_{Y})\to 0.
	\end{equation*}
	For each $x\in X$,
	\begin{equation*}
		I_{x}=
		\begin{cases}
			0 & \text{if }x\in Y\\
			\OXx & \text{if }x\notin Y
		\end{cases},
	\end{equation*}
	and hence $I_{x}$ is idempotent. It follows that $I$ is idempotent.
\end{proof}

\begin{Corollary}\label{BijectionBetweenBilocSubcatsAndIdempSubobjAndClopenSubsch}
	Let $X$ be a locally noetherian scheme. Then there exist bijections between the following sets.
	\begin{enumerate}
		\item\label{BijectionBetweenBilocSubcatsAndIdempSubobjAndClopenSubsch:BilocSubcats} The set of bilocalizing subcategories of $\QCoh X$.
		\item\label{BijectionBetweenBilocSubcatsAndIdempSubobjAndClopenSubsch:IdempSubobjs} The set of idempotent subobjects of $\OX$.
		\item\label{BijectionBetweenBilocSubcatsAndIdempSubobjAndClopenSubsch:ClopenSubschs} The set of closed subschemes of $X$ which are also open subschemes.
		\item\label{BijectionBetweenBilocSubcatsAndIdempSubobjAndClopenSubsch:ClopenSubsets} The set of subsets of $X$ which are open and closed.
	\end{enumerate}
	
	The bijection \autoref{BijectionBetweenBilocSubcatsAndIdempSubobjAndClopenSubsch:BilocSubcats}$\leftrightarrow$\autoref{BijectionBetweenBilocSubcatsAndIdempSubobjAndClopenSubsch:IdempSubobjs} is in \autoref{ClassificationOfBilocSubcatsForLocNoethSch}. The bijection \autoref{BijectionBetweenBilocSubcatsAndIdempSubobjAndClopenSubsch:IdempSubobjs}$\leftrightarrow$\autoref{BijectionBetweenBilocSubcatsAndIdempSubobjAndClopenSubsch:ClopenSubschs} is induced by the bijection in \autoref{BijectionBetweenQCohSubsheavesAndClosedSubschs}. For each element $Y$ of \autoref{BijectionBetweenBilocSubcatsAndIdempSubobjAndClopenSubsch:ClopenSubschs}, the corresponding element of \autoref{BijectionBetweenBilocSubcatsAndIdempSubobjAndClopenSubsch:ClopenSubsets} is the underlying space of $Y$.
\end{Corollary}

\begin{proof}
	This follows from \autoref{ClassificationOfBilocSubcatsForLocNoethSch} and \autoref{SubobjIsIdempIffCorrespondingClosedSubschIsOpenSubsch}.
\end{proof}

By using the classification of the prelocalizing (resp.\ localizing, closed) subcategories of $\Mod k[x]$, we can obtain a classification of the prelocalizing (resp.\ localizing, closed) subcategories for the projective line.

\begin{Example}\label{ExOfProjSch}
	Let $k$ be an algebraically closed field, and consider the projective line $X=\mbP_{k}^{1}$. Denote by $\Phi$ the set of closed points in $X$. For each $r\in\prod_{x\in\Phi}(\mbZ_{\geq 0}\cup\set{\infty})$, we define a prelocalizing subcategory $\mcY_{r}$ of $\QCoh X$ by
	\begin{equation*}
		\mcY_{r}=\conditionalset{M\in\QCoh X}{M_{x}\mfm_{x}^{r(x)}=0\text{ for each }x\in\Phi\text{ with }r(x)\neq\infty}.
	\end{equation*}
	Then by the main results (\autoref{ClassficationOfPrelocSubcatsForLocNoethSch}, \autoref{ClassificationOfLocSubcatsForLocNoethSch}, \autoref{ClassificationOfClosedSubcatsForLocNoethSch}, and \autoref{ClassificationOfBilocSubcatsForLocNoethSch}) and the examples on $\Spec k[x]$ (\autoref{ExOfPrelocSubcatsOfAffineSch}, \autoref{ExOfLocSubcatsOfAffineSch}, \autoref{ExOfClosedSubcatsOfAffineSch}, and \autoref{ExOfBilocSubcatsOfAffineSch}), the set of prelocalizing subcategories of $\QCoh X$ is
	\begin{equation*}
		\conditionalset[\Bigg]{\mcY_{r}}{r\in\prod_{x\in\Phi}(\mbZ_{\geq 0}\cup\set{\infty})}\cup\set{\QCoh X},
	\end{equation*}
	the set of localizing subcategories of $\QCoh X$ is
	\begin{equation*}
		\conditionalset[\Bigg]{\mcY_{r}}{r\in\prod_{x\in\Phi}\set{0,\infty}}\cup\set{\QCoh X},
	\end{equation*}
	the set of closed subcategories of $\QCoh X$ is
	\begin{equation*}
		\conditionalset[\Bigg]{\mcY_{r}}{r\in\bigoplus_{x\in\Phi}\mbZ_{\geq 0}}\cup\set{\QCoh X},
	\end{equation*}
	and the set of bilocalizing subcategories of $\QCoh X$ is
	\begin{equation*}
		\conditionalset{\mcY_{r}}{r=\set{0}_{x\in\Phi}}\cup\set{\QCoh X}=\set{0,\,\QCoh X}.
	\end{equation*}
\end{Example}

\begin{Example}\label{ExOfCoprodOfSch}
	For each $i\in\mbZ$, let $k_{i}$ be a field, and let $U_{i}:=\Spec k_{i}$. Consider the disjoint union $X:=\coprod_{i\in\mbZ}U_{i}$. For each subset $B$ of $\mbZ$, define a prelocalizing subcategory $\mcY_{B}$ of $\QCoh X$ by
	\begin{equation*}
		\mcY_{B}=\conditionalset{M\in\QCoh X}{M|_{U_{i}}=0\text{ for each }i\in\mbZ\setminus B}.
	\end{equation*}
	Then by \autoref{ClassficationOfPrelocSubcatsForLocNoethSch}, \autoref{ClassificationOfLocSubcatsForLocNoethSch}, \autoref{ClassificationOfClosedSubcatsForLocNoethSch}, and \autoref{ClassificationOfBilocSubcatsForLocNoethSch}, the set
	\begin{equation*}
		\conditionalset{\mcY_{B}}{B\subset\mbZ}
	\end{equation*}
	is the set of prelocalizing subcategories of $\QCoh X$, and every prelocalizing subcategory of $\QCoh X$ is bilocalizing. Therefore every local filter of $\OX$ is a principal filter. For each subset $B$ of $\mbZ$, let $I_{B}$ be the idempotent subobject of $\OX$ corresponding to the bilocalizing subcategory $\mcY_{B}$. Then the filter
	\begin{equation*}
		\mcF=\conditionalset{I_{B}}{\mbZ\setminus B\text{ is a finite set}}
	\end{equation*}
	of $\OX$ is not a local filter since $\mcF$ is not a principal filter.
\end{Example}



\end{document}